\documentclass{amsart}
\usepackage{amssymb}

\newtheorem{theorem}{Theorem}[section]
\newtheorem{corollary}[theorem]{Corollary}
\newtheorem{lemma}[theorem]{Lemma}
\newtheorem{proposition}[theorem]{Proposition}
\theoremstyle{definition}

\newtheorem{remark}[theorem]{Remark}
\newtheorem{scholium}[theorem]{Scholium}

\newtheorem{fubbyloofer}{Description}[section]

\theoremstyle{remark}
\newtheorem{acknowledgements}{Acknowledgements}

\numberwithin{equation}{section}

\newcommand{\at}{{@}}
\newcommand{\frogleg}{\ }
\newcommand{\frogkern}{\kern0.5em}
\newlength{\customskipamount}
\newlength{\displayboxwidth}
\newlength{\qedskip}
\newlength{\Jwidth}
\DeclareMathSymbol{\bftdiagup}       {\mathord}{AMSb}{"1E}
\DeclareMathSymbol{\bftdiagdown}     {\mathord}{AMSb}{"1F}
\def\openone{\hbox{\upshape \small1\kern-3.3pt\normalsize1}}
\newbox\crossbox
\def\bigcross{\mathop{\mathchoice
{\setbox\crossbox=\hbox{\Huge \bfseries $\times$}
\raise\fontdimen22\textfont2%
\hbox{\raise0.5\dp\crossbox%
\hbox{\lower0.5\ht\crossbox%
\box\crossbox}}}
{\setbox\crossbox=\hbox{\huge \bfseries $\times$}
\raise\fontdimen22\textfont2%
\hbox{\raise0.5\dp\crossbox%
\hbox{\lower0.5\ht\crossbox%
\box\crossbox}}}
{\setbox\crossbox=\hbox{\LARGE \bfseries $\times$}
\raise\fontdimen22\scriptfont2%
\hbox{\raise0.5\dp\crossbox%
\hbox{\lower0.5\ht\crossbox%
\box\crossbox}}}
{\setbox\crossbox=\hbox{\Large \bfseries $\times$}
\raise\fontdimen22\scriptscriptfont2%
\hbox{\raise0.5\dp\crossbox%
\hbox{\lower0.5\ht\crossbox%
\box\crossbox}}}}}
\newbox\ipbox
\newcommand{\ip}[2]{\left\langle #1\mathrel{\mathchoice
{\setbox\ipbox=\hbox{$\displaystyle \mathstrut #1#2$}
\vrule height\ht\ipbox width0.25pt depth\dp\ipbox}
{\setbox\ipbox=\hbox{$\textstyle \mathstrut #1#2$}
\vrule height\ht\ipbox width0.25pt depth\dp\ipbox}
{\setbox\ipbox=\hbox{$\scriptstyle \mathstrut #1#2$}
\vrule height\ht\ipbox width0.25pt depth\dp\ipbox}
{\setbox\ipbox=\hbox{$\scriptscriptstyle \mathstrut #1#2$}
\vrule height\ht\ipbox width0.25pt depth\dp\ipbox}
} #2\right\rangle}
\newbox\frogdown
\long\def\hookdownarrow{\setbox\frogdown=\hbox to 0pt{\hss
$\displaystyle \downarrow $\hss }\vrule height0pt width0pt depth1.5\ht\frogdown
\setlength{\unitlength}{0.4pt}
\begin{picture}(10,5)
\put(5,0){\oval(10,10)[t]}
\end{picture}
\lower\ht\frogdown\box\frogdown}
\def\func#1{\mathop{\rm #1}}%
\def\limfunc#1{\mathop{\rm #1}}%

\begin{document}
\title[Cuntz algebras and multiresolution wavelet analysis]{Isometries, shifts,
Cuntz algebras and multiresolution wavelet analysis of
scale $N$}
\author{Ola Bratteli}
\address{Mathematics Institute\\
University of Oslo\\
BP 1053 -- Blindern\\
N-0316 Oslo\\
Norway}
\email{bratteli%
%TCIMACRO{\TeXButton{@}{\at} }
%BeginExpansion
\at%
%EndExpansion
math.uio.no}
\thanks{Work supported in part by the U.S. National Science Foundation and the
Norwegian Research Council.}
\author{Palle E.T. Jorgensen}
\address{Department of Mathematics\\
The University of Iowa\\
14 MacLean Hall\\
Iowa City, IA 52242-1419\\
U.S.A.}
\email{jorgen%
%TCIMACRO{\TeXButton{@}{\at} }
%BeginExpansion
\at%
%EndExpansion
math.uiowa.edu}

\begin{abstract}
In this paper we show how wavelets originating from multiresolution analysis
of scale $N$ give rise to certain representations of the Cuntz algebras $%
\mathcal{O}_{N}$, and conversely how the wavelets can be recovered from
these representations. The representations are given on the Hilbert space $%
L^{2}\left( \mathbb{T}\right) $ by $\left( S_{i}\xi \right) \left( z\right)
=m_{i}\left( z\right) \xi \left( z^{N}\right) $. We characterize the Wold
decomposition of such operators. If the operators come from wavelets they
are shifts, and this can be used to realize the representation on a certain
Hardy space over $L^{2}\left( \mathbb{T}\right) $. This is used to compare
the usual scale-$2$ theory of wavelets with the scale-$N$ theory. Also some
other representations of $\mathcal{O}_{N}$ of the above form called diagonal
representations are characterized and classified up to unitary equivalence
by a homological invariant.
\end{abstract}

%TCIMACRO{\TeXButton{maketitle}{\maketitle}}
%BeginExpansion
\maketitle%
%EndExpansion
%TCIMACRO{
%\TeXButton{tableofcontents}%
%{\tableofcontents\renewcommand{\frogleg}{\\}\renewcommand{\frogkern}{\relax}}}
%BeginExpansion
\tableofcontents\renewcommand{\frogleg}{\\}\renewcommand{\frogkern}{\relax}%
%EndExpansion
\setlength{\displayboxwidth}{\textwidth}\addtolength{\displayboxwidth}{-2%
\leftmargini}

\section{\label{Intr}%
%TCIMACRO{\TeXButton{frogkern}{\protect\frogkern } }
%BeginExpansion
\protect\frogkern %
%EndExpansion
Introduction}

Continuing \cite{BJP96}, \cite{BrJo96b}, \cite{Jor95} we will consider some
representations $\pi $ of the Cuntz algebra $\mathcal{O}_{N}$ coming from
wavelet theory. Our ultimate goal is to establish connections between
certain representations of $\mathcal{O}_{N}$, and their decompositions, and
wavelet decompositions for the wavelets arising from multiresolutions with
scaling $N $. The map from wavelets into representations is described in
detail in Section \ref{Comp} (when $N=2$), Section \ref{Scal} (when the
translates of the father function are orthogonal), and in Section \ref{Fath}
(in more general cases). Unfortunately we have only partial results on how
to go the other way, from representations to wavelets, and for the moment
the path in both directions leads past certain functions from the circle $%
\mathbb{T}$ into unitary $N\times N$ matrices given by (\ref{Eq1.11}). We
will discuss further the connection between representations and wavelets at
the end of this introduction.

Recall from \cite{Cun77} that $\mathcal{O}_{N}$ is the $C^{*}$-algebra
generated by $N\in \mathbb{N}$ isometries $s_{0},s_{1},\dots ,s_{N-1}$
satisfying 
\begin{equation}
s_{i}^{*}s_{j}^{{}}=\delta _{ij}^{{}}\openone  \label{Eq1.1}
\end{equation}
and 
\begin{equation}
\sum_{i=0}^{N-1}s_{i}^{{}}s_{i}^{*}=\openone.  \label{Eq1.2}
\end{equation}
The representations we will consider are realized on Hilbert spaces $%
\mathcal{H}=L^{2}\left( \Omega ,\mu \right) $ where $\Omega $ is a measure
space and $\mu $ is a probability measure on $\Omega $. We define the
representations in terms of certain maps $\sigma _{i}:\Omega \rightarrow
\Omega $ with the property that $\mu \left( \sigma _{i}\left( \Omega \right)
\cap \sigma _{j}\left( \Omega \right) \right) =0$ for all $i\ne j$, and if $%
\rho _{i}=\mu \left( \sigma _{i}\left( \Omega \right) \right) $ then $\rho
_{i}>0$ and $\sum_{i=0}^{N-1}\rho _{i}=1$, i.e., $\left\{ \sigma _{0}\left(
\Omega \right) ,\dots ,\sigma _{N-1}\left( \Omega \right) \right\} $ is a
partition of $\Omega $ up to measure zero.

We further assume that 
\begin{equation}
\int_{\Omega }f\,d\mu =\sum_{i=0}^{N-1}\rho _{i}\int_{\Omega }f\circ \sigma
_{i}\,d\mu  \label{Eq1.3}
\end{equation}
for all $f\in L^{\infty }\left( \Omega ,\mu \right) $, or, equivalently, 
\begin{equation}
\mu \left( \sigma _{i}\left( Y\right) \right) =\rho _{i}\mu \left( Y\right)
\label{Eq1.4}
\end{equation}
for $i\in \mathbb{Z}_{N}=\left\{ 0,1,\dots ,N-1\right\} $ and all measurable 
$Y\subset \Omega $. Since $\rho _{i}>0$, this entails that the $\sigma _{i}$%
's are injections up to measure zero, and hence we may define an $N$-to-$1$
map $\sigma :\Omega \rightarrow \Omega $, well defined up to measure zero,
by $\sigma \circ \sigma _{i}=\func{id}$ for $i\in \mathbb{Z}_{N}$. Finally,
we assume that the sets $\sigma _{i_{1}}\sigma _{i_{2}}\cdots \sigma
_{i_{k}}\left( \Omega \right) $ generate the $\sigma $-algebra of measurable
sets of $\Omega $ up to sets of measure zero. Thus $\left( \Omega ,\mu
,\sigma _{i}\right) $ is canonically isomorphic by a coding map to $\left( %
\bigcross_{k=1}^{\infty }\mathbb{Z}_{N}%
%TCIMACRO{\TeXButton{vphantom}{\vphantom{\sigma _{i}^{\left( 0\right) }}}}
%BeginExpansion
\vphantom{\sigma _{i}^{\left( 0\right) }}%
%EndExpansion
\right. $, the product measure of measure on $\mathbb{Z}_{N}$ with weights $%
\left. \rho _{0},\rho _{1},\dots ,\rho _{N-1},\sigma _{i}^{\left( 0\right) }%
%TCIMACRO{
%\TeXButton{vphantom}{\vphantom{\bigcross _{k=1}^{\infty }\mathbb{Z}_{N}}}}
%BeginExpansion
\vphantom{\bigcross _{k=1}^{\infty }\mathbb{Z}_{N}}%
%EndExpansion
\right) $, where 
\begin{equation}
\sigma _{i}^{\left( 0\right) }\left( x_{1},x_{2},\dots \right) =\left(
i,x_{1},x_{2},\dots \right) ,  \label{Eq1.5}
\end{equation}
and then $\sigma $ is defined by 
\begin{equation}
\sigma \left( x_{1},x_{2},\dots \right) =\left( x_{2},x_{3},\dots \right) .
\label{Eq1.6}
\end{equation}
We will list many other realizations of $\left( \Omega ,\mu ,\sigma
_{i}\right) $ below.

The announced representations $s_{i}\rightarrow S_{i}$ of $\mathcal{O}_{N}$
on $L^{2}\left( \Omega ,\mu \right) $ are defined in terms of measurable
functions $m_{0},m_{1},\dots ,m_{N-1}$ from $\Omega $ into $\mathbb{C}$ with
the property that the $N\times N$ matrix 
\begin{equation}
%TCIMACRO{
%\TeXButton{pmatrix}{\begin{pmatrix}
% \sqrt{\rho _{0}}m_{0}(\sigma _{0}(x)) & \sqrt{\rho _{1}}m_{0}(\sigma _{1}(x))
% & \dots  & \sqrt{\rho _{N-1}}m_{0}(\sigma _{N-1}(x)) \\ 
% \sqrt{\rho _{0}}m_{1}(\sigma _{0}(x)) & \sqrt{\rho _{1}}m_{1}(\sigma _{1}(x))
% & \dots  & \sqrt{\rho _{N-1}}m_{1}(\sigma _{N-1}(x)) \\ 
% \vdots  & \vdots  & \ddots  & \vdots  \\ 
% \sqrt{\rho _{0}}m_{N-1}(\sigma _{0}(x)) & \sqrt{\rho _{1}}m_{N-1}(\sigma
% _{1}(x)) & \dots  & \sqrt{\rho _{N-1}}m_{N-1}(\sigma _{N-1}(x))
% \end{pmatrix}} }
%BeginExpansion
\begin{pmatrix}
 \sqrt{\rho _{0}}m_{0}(\sigma _{0}(x)) & \sqrt{\rho _{1}}m_{0}(\sigma _{1}(x))
 & \dots  & \sqrt{\rho _{N-1}}m_{0}(\sigma _{N-1}(x)) \\ 
 \sqrt{\rho _{0}}m_{1}(\sigma _{0}(x)) & \sqrt{\rho _{1}}m_{1}(\sigma _{1}(x))
 & \dots  & \sqrt{\rho _{N-1}}m_{1}(\sigma _{N-1}(x)) \\ 
 \vdots  & \vdots  & \ddots  & \vdots  \\ 
 \sqrt{\rho _{0}}m_{N-1}(\sigma _{0}(x)) & \sqrt{\rho _{1}}m_{N-1}(\sigma
 _{1}(x)) & \dots  & \sqrt{\rho _{N-1}}m_{N-1}(\sigma _{N-1}(x))
 \end{pmatrix}%
%EndExpansion
\label{Eq1.7}
\end{equation}
is unitary for almost all $x\in \Omega $. The representation is given by 
\begin{equation}
\left( S_{i}\xi \right) \left( x\right) =m_{i}\left( x\right) \xi \left(
\sigma \left( x\right) \right)  \label{Eq1.8}
\end{equation}
and one computes that this is really a representation of the Cuntz algebra
and 
\begin{equation}
\left( S_{i}^{*}\xi \right) \left( x\right) =\sum_{k\in \mathbb{Z}_{N}}\rho
_{k}\bar{m}_{i}\left( \sigma _{k}\left( x\right) \right) \xi \left( \sigma
_{k}\left( x\right) \right) ;  \label{Eq1.9}
\end{equation}
see \cite{Jor95}. The computations are (for $\eta ,\xi \in L^{2}\left(
\Omega ,\mu \right) $): 
\begin{align*}
\ip{\eta }{S_{i}^{*}\xi }& =\ip{S_{i}^{}\eta }{\xi } \\
& =\int_{\Omega }\bar{m}_{i}\left( x\right) \bar{\eta}\left( \sigma \left(
x\right) \right) \xi \left( x\right) \,d\mu \left( x\right) \\
& =\sum_{k=0}^{N-1}\rho _{k}\int_{\Omega }\bar{m}_{i}\left( \sigma
_{k}\left( x\right) \right) \bar{\eta}\left( x\right) \xi \left( \sigma
_{k}\left( x\right) \right) \,d\mu \left( x\right) ,
\end{align*}
which used (\ref{Eq1.3}) and gives (\ref{Eq1.9}), and thus 
\begin{align*}
\left( S_{i}^{*}S_{j}^{{}}\xi \right) \left( x\right) & =\sum_{k\in %
\mathbb{Z}_{N}}\rho _{k}\bar{m}_{i}\left( \sigma _{k}\left( x\right) \right)
m_{j}\left( \sigma _{k}\left( x\right) \right) \xi \left( \sigma \left(
\sigma _{k}\left( x\right) \right) \right) \\
& =\delta _{ij}\xi \left( x\right)
\end{align*}
by unitarity of (\ref{Eq1.7}), which gives (\ref{Eq1.1}). The formula (\ref
{Eq1.2}) follows similarly from unitarity of (\ref{Eq1.7}): 
\begin{equation*}
\smash[t]{\sum_{i}\left\| S_{i}^{*}\xi \right\| ^{2}=\sum_{i}\sum_{k\,l}\rho
_{k}\rho _{l}\int_{\Omega }m_{i}\left( \sigma _{k}\left( x\right) \right)
\bar{\xi}\left( \sigma _{k}\left( x\right) \right) \bar{m}_{i}\left( \sigma
_{l}\left( x\right) \right) \xi \left( \sigma _{l}\left( x\right) \right)
\,d\mu \left( x\right) }.
\end{equation*}
By unitarity of (\ref{Eq1.7}), 
\begin{equation*}
\sum_{i}\sqrt{\rho _{k}}\bar{m}_{i}\left( \sigma _{k}\left( x\right) \right) 
\sqrt{\rho _{l}}m_{i}\left( \sigma _{l}\left( x\right) \right) =\delta _{lk},
\end{equation*}
so 
\begin{align*}
\sum_{i}\left\| S_{i}^{*}\xi \right\| ^{2}& =\sum_{k\,l}\sqrt{\rho _{k}}%
\sqrt{\rho _{l}}\delta _{lk}\int_{\Omega }\bar{\xi}\left( \sigma _{k}\left(
x\right) \right) \xi \left( \sigma _{l}\left( x\right) \right) \,d\mu \left(
x\right) \\
& =\sum_{k}\rho _{k}\int_{\Omega }\left| \xi \left( \sigma _{k}\left(
x\right) \right) \right| ^{2}\,d\mu \left( x\right) \\
& =\left\| \xi \right\| ^{2}
\end{align*}
and (\ref{Eq1.2}) follows.

Before surveying the by now rather rich theory of representations of the
form (\ref{Eq1.8}), we will give some alternative descriptions of the system 
$\left( \Omega ,\mu ,\sigma _{i}\right) $ which are convenient to use in
special circumstances. From now, and through the rest of the paper, we make
the simplifying assumption 
\begin{equation}
\rho _{k}=\frac{1}{N}  \label{Eq1.10}
\end{equation}
for $k\in \mathbb{Z}_{N}$, although this assumption is easy to remove in
many cases. Thus the condition of unitarity is that the $N\times N$ matrix 
\begin{equation}
\frac{1}{\sqrt{N}} 
%TCIMACRO{
%\TeXButton{pmatrix}{\begin{pmatrix}
% m_{0}(\sigma _{0}(x)) & m_{0}(\sigma _{1}(x))
% & \dots  & m_{0}(\sigma _{N-1}(x)) \\ 
% m_{1}(\sigma _{0}(x)) & m_{1}(\sigma _{1}(x))
% & \dots  & m_{1}(\sigma _{N-1}(x)) \\ 
% \vdots  & \vdots  & \ddots  & \vdots  \\ 
% m_{N-1}(\sigma _{0}(x)) & m_{N-1}(\sigma
% _{1}(x)) & \dots  & m_{N-1}(\sigma _{N-1}(x))
% \end{pmatrix}} }
%BeginExpansion
\begin{pmatrix}
 m_{0}(\sigma _{0}(x)) & m_{0}(\sigma _{1}(x))
 & \dots  & m_{0}(\sigma _{N-1}(x)) \\ 
 m_{1}(\sigma _{0}(x)) & m_{1}(\sigma _{1}(x))
 & \dots  & m_{1}(\sigma _{N-1}(x)) \\ 
 \vdots  & \vdots  & \ddots  & \vdots  \\ 
 m_{N-1}(\sigma _{0}(x)) & m_{N-1}(\sigma
 _{1}(x)) & \dots  & m_{N-1}(\sigma _{N-1}(x))
 \end{pmatrix}%
%EndExpansion
\label{Eq1.11}
\end{equation}
is unitary for almost all $x\in \Omega $, and 
\begin{align}
\left( S_{i}\xi \right) \left( x\right) & =m_{i}\left( x\right) \xi \left(
\sigma \left( x\right) \right) ,  \label{Eq1.12} \\
\left( S_{i}^{*}\xi \right) \left( x\right) & =\frac{1}{N}\sum_{k\in %
\mathbb{Z}_{N}}\bar{m}_{i}\left( \sigma _{k}\left( x\right) \right) \xi
\left( \sigma _{k}\left( x\right) \right) .  \label{Eq1.13}
\end{align}

Note conversely that if $S_{i}$ is given by (\ref{Eq1.12}) (respectively (%
\ref{Eq1.8})) and if the $S_{i}$'s define a representation of $\mathcal{O}%
_{N}$, then the matrix (\ref{Eq1.11}) (respectively (\ref{Eq1.7})) is
unitary. Also, as the ranges of the maps $\sigma _{0},\dots ,\sigma _{N-1}$
are disjoint, any function from $\mathbb{T}$ into unitary matrices has the
form (\ref{Eq1.11}) (respectively (\ref{Eq1.7})). Compare this with the
well-known fact that if $S_{0},\dots ,S_{N-1}$ and $T_{0},\dots ,T_{N-1}$
are any two realizations of $\mathcal{O}_{N}$ on a Hilbert space $\mathcal{H}
$, there is a unique unitary $U$ on $\mathcal{H}$ such that $S_{k}=UT_{k}$,
namely $U=\sum_{k}^{{}}S_{k}^{{}}T_{k}^{*}$. Alternatively, if $%
M_{ij}^{{}}=T_{i}^{*}S_{j}^{{}}$, then $S_{k}=\sum_{j}T_{j}M_{jk}$, and $%
\left[ M_{ij}\right] $ is a unitary matrix on $\mathcal{H}\otimes \mathbb{C}%
^{N}$. Our representations correspond to the special case that 
\begin{equation*}
\left( T_{i}\xi \right) \left( x\right) =\sqrt{N}\chi _{\sigma _{i}\left(
\Omega \right) }\left( x\right) \xi \left( \sigma \left( x\right) \right)
\end{equation*}
and the $M_{ij}$ are multiplication operators defined by $m_{i}\left( \sigma
_{j}\left( \,\cdot \,\right) \right) \in L^{\infty }\left( \Omega \right) $.

Here are some equivalent descriptions of $\left( \Omega ,\mu ,\sigma
_{i}\right) $. Description \ref{Des2} will be particularly convenient in
connection with the examples coming from wavelets.

\begin{fubbyloofer}
\label{Des1} 
\begin{equation}
\begin{aligned} \Omega & =\bigcross_{k=1}^{\infty }\mathbb{Z}_{N}, \\ \mu &
=\text{Normalized Haar measure,} \\ \sigma _{i}\left( x_{1},x_{2},\dots
\right) & =\left( i,x_{1},x_{2},\dots \right) , \\ \sigma \left(
x_{1},x_{2},\dots \right) & =\left( x_{2},x_{3},\dots \right) . \end{aligned}
\label{Eq1.14}
\end{equation}
\end{fubbyloofer}

\begin{fubbyloofer}
\label{Des2} 
\begin{equation}
\begin{aligned} \Omega & =\mathbb{T}=\text{the unit circle in }\mathbb{C},
\\ \mu & =\text{Normalized Haar measure,} \\ \sigma _{k}\left( e^{2\pi
i\theta }\right) & =\exp \left( 2\pi i(\theta +k)/N\right) \text{ when
}0\leq \theta <1, \\ \sigma \left( z\right) & =z^{N}; \end{aligned}
\label{Eq1.15}
\end{equation}
so (\ref{Eq1.12})--(\ref{Eq1.13}) take the form 
%TCIMACRO{
%\TeXButton{Custom Aligned Equation}{\begin{align}
%\left( S_{i}\xi \right) \left( z\right) &=m_{i}\left( z\right) \xi \left(
%z^{N}\right) ,  \label{Eq1.16} \\
%\left( S_{i}^{*}\xi \right) \left( z\right) 
%&=\frac{1}{N}\sum_{w^{N}=z}\bar{m}_{i}\left( w\right) \xi \left( w\right) .  
%\label{Eq1.17}
%\end{align}}}
%BeginExpansion
\begin{align}
\left( S_{i}\xi \right) \left( z\right) &=m_{i}\left( z\right) \xi \left(
z^{N}\right) ,  \label{Eq1.16} \\
\left( S_{i}^{*}\xi \right) \left( z\right) 
&=\frac{1}{N}\sum_{w^{N}=z}\bar{m}_{i}\left( w\right) \xi \left( w\right) .  
\label{Eq1.17}
\end{align}%
%EndExpansion
\end{fubbyloofer}

\begin{fubbyloofer}
\label{Des3}This example has an obvious $\nu $-dimensional analogue,
replacing $N$ by a $\nu \times \nu $ matrix $\mathbf{N}$ with integer
entries such that $\left| \det \left( \mathbf{N}\right) \right| =N$. Then $%
\Omega =\mathbb{T}^{\nu }$, $\mu =$ Normalized Haar measure, and $\sigma
\left( x%
%TCIMACRO{\TeXButton{mod ZZ^\nu}{\pmod{\mathbb{Z}^\nu }}}
%BeginExpansion
\pmod{\mathbb{Z}^\nu }%
%EndExpansion
\right) =\mathbf{N}x%
%TCIMACRO{\TeXButton{mod ZZ^\nu}{\pmod{\mathbb{Z}^\nu }}}
%BeginExpansion
\pmod{\mathbb{Z}^\nu }%
%EndExpansion
$ for $x\in \mathbb{R}^{\nu }$. A somewhat different turn on this idea is
the following: let $\mathbf{N}$ be a $\nu \times \nu $ matrix with integer
entries such that all the (complex) eigenvalues of $\mathbf{N}$ have modulus
greater than $1$, and assume $N=\left| \det \left( \mathbf{N}\right) \right| 
$. Let $D\subset \mathbb{Z}^{\nu }$ be a set of $N$ points in $\mathbb{Z}%
^{\nu }$ which are incongruent modulo $\mathbf{N}\left( \mathbb{Z}^{\nu
}\right) $, i.e., such that each point $m\in \mathbb{Z}^{\nu }$ has a unique
expansion $m=d+\mathbf{N}m$ for $d\in D$, $m\in \mathbb{Z}^{\nu }$. It
follows easily \cite{BrJo96b} that there is then a unique compact subset $%
\mathbf{T}\subset \mathbb{R}^{\nu }$ such that 
\begin{equation}
\mathbf{T}=\mathbf{N}^{-1}\bigcup_{d\in D}\left( d+\mathbf{T}\right) 
\label{Eq1.19}
\end{equation}
If $\mu $ is Lebesgue measure on $\mathbb{R}^{\nu }$, we have $\mu \left( d+%
\mathbf{T}\right) =\mu \left( \mathbf{T}\right) $ for each $d\in D$, while $%
\mu \left( \mathbf{N}\left( \mathbf{T}\right) \right) =\left| \det \left( 
\mathbf{N}\right) \right| \mu \left( \mathbf{T}\right) =N\mu \left( \mathbf{T%
}\right) $, and hence the sets $\mathbf{N}^{-1}\left( d+\mathbf{T}\right) $
must be mutually disjoint up to sets of measure zero, since they have union $%
\mathbf{T}$. Hence we may define $\Omega =\mathbf{T}$, $\mu =$ Lebesgue
measure$|_{\mathbf{T}}$ (except for normalization) and 
\begin{equation}
\sigma _{i}\left( x\right) =\mathbf{N}^{-1}\left( d_{i}+x\right) 
\label{Eq1.20}
\end{equation}
where $D=\left\{ d_{0},d_{1},\dots ,d_{N-1}\right\} $ is an enumeration of $D
$, and 
\begin{equation}
\begin{minipage}[t]{\displayboxwidth}\raggedright $\sigma \left( x\right)
=y\in \mathbf{T}$ (the point such that there is a $d\in D$ with
$x=N^{-1}\left( d+y\right) $). \end{minipage}  \label{Eq1.21}
\end{equation}
Note $\mathbf{T}$ may or may not be a $\mathbb{Z}^{\nu }$ tiling of $%
\mathbb{R}^{\nu }$, and it may be a union of tiles. An exhaustive discussion
of the rich possibilities is given in \cite{BrJo96b}, based on \cite{Hut81}, 
\cite{BrJo96a}, \cite{JoPe94}, \cite{JoPe96}.
\end{fubbyloofer}

\begin{fubbyloofer}
\label{Des4}Let $\mathbb{C}$ be the Riemann sphere and let $R\left( z\right)
=P\left( z\right) /Q\left( z\right) $ be a rational function, where the
polynomials $P\left( z\right) $ and $Q\left( z\right) $ have no common
linear factor. If $N=\max \left\{ \deg P,\deg Q\right\} $, then $R$ defines
an $N$-fold cover of the Riemann sphere. Now, let $\Omega $ be the Julia
set, i.e., $\Omega $ is the set of $z_{0}\in \mathbb{C}$ such that the
sequence of iterations $R^{n}\left( z\right) $ is not a normal family near $%
z_{0}$, i.e., there is no neighborhood of $z_{0}\in \mathbb{C}$ such that
the sequence $R^{n}\left( z\right) $ is uniformly bounded for $z$ in the
neighborhood. It is known that if $z_{0}$ is an attracting periodic point,
then the boundary of the region of attraction of $z_{0}$ under $R$ is equal
to $\Omega $ \cite[Theorem 2.1]{CaGa93}, and also that $\Omega $ is the
closure of the repelling periodic points under $R$, \cite[Theorem 3.1]
{CaGa93}. Following \cite{Bro65}, \cite{CaGa93}, if $\nu $ is any
probability measure on $\Omega $, define the energy integral 
\begin{equation}
I\left( \nu \right) =\int_{\Omega }\int_{\Omega }\log \left( \frac{1}{\left|
\zeta -\eta \right| }\right) \,d\nu \left( \eta \right) \,d\nu \left( \zeta
\right) .  \label{Eq1.23}
\end{equation}
Then $\inf_{\nu }I\left( \nu \right) =0$, and there is a unique probability
measure $\mu $ such that $I\left( \mu \right) =0$. Furthermore, for a
generic set of points $z_{0}\in \Omega $, the measure $\mu $ can be obtained
as the weak limit of the set of probability measures defined by 
\begin{equation}
\mu _{n}=\frac{1}{N^{n}}\sum_{%
\substack{w \\ \makebox[0pt]{\hss
{$\scriptstyle R^{n}\left( w\right) =z_{0}$}\hss }}}\delta _{w}.
\label{Eq1.24}
\end{equation}
Then $\left( \Omega ,\mu ,\sigma \right) $ satisfies all our requirements,
while $\sigma _{i}$ corresponds to an explicit choice of Riemann cover. In
the special case $R\left( z\right) =z^{N}$ we recover Description \ref{Des2}.
\end{fubbyloofer}

Let us now describe some known results on the representations defined by (%
\ref{Eq1.12})--(\ref{Eq1.13}), alias (\ref{Eq1.16})--(\ref{Eq1.17}) (recall
that we assume $\rho _{i}=\frac{1}{N}$ throughout). If 
\begin{equation}
m_{i}\left( x\right) =\eta _{i}^{-1}\chi _{\sigma _{i}\Omega }
\label{Eq1.25}
\end{equation}
where $\eta _{i}\in \mathbb{C}$ are nonzero complex numbers with $%
\sum_{i}\left| \eta _{i}\right| ^{2}=1$, then $\openone $ is cyclic for the
representation, and represents the so-called Cuntz state on $\mathcal{O}_{N}$%
; see, e.g., \cite[Section 8]{BJP96}. In particular these representations
are irreducible. In \cite[Section 8]{BJP96} we considered the particular
representation with 
\begin{equation}
m_{i}\left( x_{0},x_{1},x_{2},\dots \right) =N^{\frac{1}{2}}\delta
_{ix_{0}}\left\langle i,x_{1}\right\rangle ,  \label{Eq1.26}
\end{equation}
where $\left\langle \;,\;\right\rangle $ is the usual Pontrjagin duality $%
\mathbb{Z}_{N}\times \hat{\mathbb{Z}}_{N}\rightarrow \mathbb{T}$, and showed
that the resulting representation is irreducible and disjoint from all the
Cuntz state representations. In \cite{BrJo96a} we considered the
representations with 
\begin{equation}
m_{i}\left( x_{0},x_{1},\dots \right) =N^{\frac{1}{2}}\delta
_{ix_{0}}u\left( x_{0},x_{1},\dots \right) ,  \label{Eq1.27}
\end{equation}
where $u:\mathbb{T}\rightarrow \mathbb{T}$ is a measurable function (using
the identification $\Omega =\mathbb{T}$ of Description \ref{Des2}), and we
showed in Proposition 7.1 that the resulting representation of $\mathcal{O}%
_{N}$ is irreducible, and even the restriction to the canonical $\limfunc{UHF%
}$-subalgebra $\limfunc{UHF}\nolimits_{N}\subset \mathcal{O}_{N}$ is
irreducible. $\limfunc{UHF}\nolimits_{N}$ is the $C^{*}$-subalgebra
generated by the monomials $s_{I}^{{}}s_{J}^{*}$ with $\left| I\right|
=\left| J\right| $. (See also Remark \ref{Rem6.2} of the present paper.)
Here $I=\left( i_{1},i_{2},\dots ,i_{n}\right) $ is a finite sequence in $%
\mathbb{Z}_{N}$, and $\left| I\right| =n$. See \cite{BJP96}, \cite{BrJo96a}
for details. The significance of the subalgebra $\limfunc{UHF}\nolimits_{N}$
for our representations derives from the work of Powers on endomorphisms of
operator algebras \cite{Pow88}. His endomorphisms correspond to
representations of $\mathcal{O}_{N}$, and the endomorphisms are \emph{shifts}
in the sense of Powers iff the corresponding representation is irreducible
when restricted to $\limfunc{UHF}\nolimits_{N}$. One of the main results of
the present paper, Corollary \ref{Cor6.2}, states that the irreducible
representations obtained from two such functions $u_{1},u_{2}:\mathbb{T}%
\rightarrow \mathbb{T}$ are unitarily equivalent if and only if there is
another measurable function $\Delta :\mathbb{T}\rightarrow \mathbb{T}$ such
that 
\begin{equation}
u_{1}\left( z\right) \Delta \left( z^{N}\right) =\Delta \left( z\right)
u_{2}\left( z\right) .  \label{Eq1.28}
\end{equation}
In the product space language this relation states that 
\begin{equation}
u_{1}\left( x_{0},x_{1},\dots \right) \Delta \left( x_{1},x_{2},\dots
\right) =\Delta \left( x_{0},x_{1},\dots \right) u_{2}\left(
x_{0},x_{1},\dots \right) .  \label{Eq1.29}
\end{equation}
Thus, if for example $u_{2}=1$, we see that some function $u_{1}$ of the
form 
\begin{equation}
u_{1}\left( x_{0},x_{1},\dots \right) =u_{1}\left( x_{0},x_{1}\right) =\frac{%
\Delta \left( x_{0}\right) }{\Delta \left( x_{1}\right) }  \label{Eq1.30}
\end{equation}
will define representations unitarily equivalent to the representation
defined by the particular Cuntz state $S_{i}^{*}\openone=N^{-\frac{1}{2}}%
\openone$. There are of course functions $u_{1}\left( x_{0},x_{1}\right) $
that do not have this form, for example the function $u_{1}\left(
x_{0},x_{1},\dots \right) =\left\langle x_{0},x_{1}\right\rangle $ in 
\cite[Section 8]{BJP96}. This can be used to recover the result from that
paper. We use the notation $\left\langle x_{0},x_{1}\right\rangle =\exp
\left( i\frac{x_{0}x_{1}2\pi }{N}\right) $ for $x_{0},x_{1}\in \mathbb{Z}%
_{N}=\left\{ 0,1,\dots ,N-1\right\} $.

In Section \ref{CharPi} we will give an intrinsic characterization of the
representations $\pi $ of $\mathcal{O}_{N}$ which are given by (\ref{Eq1.27}%
). If $\mathcal{D}_{N}$ is the canonical diagonal $C^{*}$-subalgebra of $%
\limfunc{UHF}\nolimits_{N}$, i.e., $\mathcal{D}_{N}$ is the closure of the
linear span of elements of the form $s_{I}^{{}}s_{I}^{*}$, the
characterization up to a decoding of $\Omega $ is simply that $\pi \left( 
\mathcal{D}_{N}\right) ^{\prime \prime }\subset M_{L^{\infty }\left( %
\mathbb{T}\right) }$, where $M_{L^{\infty }\left( \mathbb{T}\right) }$ is
the image of $L^{\infty }\left( \mathbb{T}\right) $ acting as multiplication
operators on $L^{2}\left( \mathbb{T}\right) $.

In \cite{BrJo96b} and \cite{DaPi96}, representations of the form (\ref
{Eq1.16}) with 
\begin{equation}
m_{i}\left( z\right) =\lambda _{i}z^{d_{i}}  \label{Eq1.31}
\end{equation}
were considered, with $\lambda _{i}\in \mathbb{T}$ and $D=\left\{
d_{0},\dots ,d_{N-1}\right\} $ a set of $N$ integers incongruent modulo $N$.
These representations turn out not to be irreducible, but at least when $%
\lambda _{i}=1$ they decompose into a discrete direct sum of mutually
disjoint irreducible representations of $\mathcal{O}_{N}$; and the
restriction to $\limfunc{UHF}\nolimits_{N}$ decomposes similarly \cite
{BrJo96b}. When $\lambda _{i}\neq 1$, even continuous decompositions may
occur \cite{DaPi96}.

Let us give a more intrinsic characterization of the representations of $%
\mathcal{O}_{N}$ given by (\ref{Eq1.8}).

\begin{proposition}
\label{Pro1.1}Assume that $\left( \Omega ,\mu ,\sigma _{i}\right) $
satisfies the requirements around \textup{(\ref{Eq1.3})--(\ref{Eq1.4})} and
define an endomorphism $\bar{\sigma}$ of $L^{\infty }\left( \Omega ,\mu
\right) $ by $\left( \bar{\sigma}f\right) \left( x\right) =f\left( \sigma
\left( x\right) \right) $. Let $s_{i}\rightarrow S_{i}$ be a representation
of $\mathcal{O}_{N}$ on $L^{2}\left( \Omega ,\mu \right) $. Then the
following conditions are equivalent. 
%TCIMACRO{
%\TeXButton{Custom Aligned Equation}{\begin{align}
%& \begin{minipage}[t]{\displayboxwidth}\raggedright There are functions
%$m_{i}\in L^{\infty }\left( \Omega \right) $ such that $\left( S_{i} \xi
%\right) \left( x\right) =m_{i}\left( x\right) \xi \left( \sigma \left(
%x\right) \right) $ for all $\xi \in L^{2}\left( \Omega ,\mu \right) $, $x
%\in \Omega $. \end{minipage}  \label{Eq1.32} \\
%& \begin{minipage}[t]{\displayboxwidth}\raggedright $\sum
%_{i=0}^{N-1}S_{i}^{}M_{f}^{}S_{i}^{*}=M_{\bar{\sigma }\left( f\right) }^{}$
%for all $f\in L^{\infty }\left( \Omega \right) $, where $M_{f}$ is the
%multiplication operator defined by $f$ on $L^{2}\left( \Omega ,\mu \right)
%$. \end{minipage}  \label{Eq1.33}
%\end{align}} }
%BeginExpansion
\begin{align}
& \begin{minipage}[t]{\displayboxwidth}\raggedright There are functions
$m_{i}\in L^{\infty }\left( \Omega \right) $ such that $\left( S_{i} \xi
\right) \left( x\right) =m_{i}\left( x\right) \xi \left( \sigma \left(
x\right) \right) $ for all $\xi \in L^{2}\left( \Omega ,\mu \right) $, $x
\in \Omega $. \end{minipage}  \label{Eq1.32} \\
& \begin{minipage}[t]{\displayboxwidth}\raggedright $\sum
_{i=0}^{N-1}S_{i}^{}M_{f}^{}S_{i}^{*}=M_{\bar{\sigma }\left( f\right) }^{}$
for all $f\in L^{\infty }\left( \Omega \right) $, where $M_{f}$ is the
multiplication operator defined by $f$ on $L^{2}\left( \Omega ,\mu \right)
$. \end{minipage}  \label{Eq1.33}
\end{align}%
%EndExpansion
Furthermore, when these conditions are fulfilled, then 
\begin{equation}
m_{i}=S_{i}\openone.  \label{Eq1.34}
\end{equation}
\end{proposition}

%TCIMACRO{\TeXButton{Begin Proof}{\begin{proof}}}
%BeginExpansion
\begin{proof}%
%EndExpansion
(\ref{Eq1.32})$\Rightarrow $(\ref{Eq1.33}). By (\ref{Eq1.8}) and (\ref{Eq1.9}%
) we have, for $f\in L^{\infty }\left( \Omega \right) $, $\xi \in
L^{2}\left( \Omega ,\mu \right) $, 
\begin{align*}
\sum_{i=0}^{N-1}S_{i}^{{}}M_{f}^{}S_{i}^{*}\xi \left( x\right) & =\sum_{i\in %
\mathbb{Z}_{N}}m_{i}\left( x\right) \left( M_{f}^{}S_{i}^{*}\xi \right)
\left( \sigma \left( x\right) \right) \\
& =\sum_{i\in \mathbb{Z}_{N}}m_{i}\left( x\right) f\left( \sigma \left(
x\right) \right) \left( S_{i}^{*}\xi \right) \left( \sigma \left( x\right)
\right) \\
& =\sum_{i\in \mathbb{Z}_{N}}m_{i}\left( x\right) f\left( \sigma \left(
x\right) \right) \sum_{k\in \mathbb{Z}_{N}}\rho _{k}\bar{m}_{i}\left( \sigma
_{k}\sigma \left( x\right) \right) \xi \left( \sigma _{k}\sigma \left(
x\right) \right) .
\end{align*}
Now, let $k_{x}\in \mathbb{Z}_{N}$ be the unique (for almost all $x$) number
such that $x=\sigma _{k_{x}}\sigma \left( x\right) $. By unitarity of (\ref
{Eq1.7}) for $x:=\sigma \left( x\right) $ we have 
\begin{equation*}
\sum_{i\in \mathbb{Z}_{N}}m_{i}\left( x\right) \bar{m}_{i}\left( \sigma
_{k}\sigma \left( x\right) \right) =\begin{cases} \rho _{k_x}^{-1} &\text{if
$k=k_x$} \\ 0 &\text{otherwise} \end{cases}
\end{equation*}
and hence (\ref{Eq1.33}) follows.

(\ref{Eq1.33})$\Rightarrow $(\ref{Eq1.32}). Put 
\begin{equation*}
m_{i}=S_{i}\openone.
\end{equation*}
If $f\in L^{\infty }\left( \Omega \right) $, we have 
\begin{align*}
M_{\bar{\sigma}\left( f\right) }^{{}}S_{j}& =\sum_{i\in \mathbb{Z}%
_{N}}S_{i}^{{}}M_{f}^{{}}S_{i}^{*}S_{j}^{{}} \\
& =S_{j}M_{f}^{{}}
\end{align*}
and applying this to $\openone$ we have 
\begin{equation*}
f\left( \sigma \left( x\right) \right) m_{j}\left( x\right) =\left(
S_{j}f\right) \left( x\right) .
\end{equation*}
As $L^{\infty }\left( \Omega \right) $ is dense in $L^{2}\left( \Omega
\right) $, this implies (\ref{Eq1.32}).%
%TCIMACRO{\TeXButton{End Proof}{\end{proof}}}
%BeginExpansion
\end{proof}%
%EndExpansion

Let us remark that not all representations of $\mathcal{O}_{N}$ on a
separable Hilbert space $\mathcal{H}$ have the form (\ref{Eq1.12}) for a
suitable realization of $\mathcal{H}$ as $L^{2}\left( \Omega ,\mu \right) $.
We will for example establish in Theorem \ref{Thm3.1} that the unitary parts
of the Wold decompositions of the respective generators $S_{i}$ have to be
zero- or one-dimensional: and, in the case that the representation comes
from a wavelet, they have to be zero-dimensional by Lemma \ref{Lem7.3}. This
is already a severe restriction, which for example immediately implies that
none of the representations coming from monomials $m_{i}$ on $\mathbb{T}$
considered in \cite{BrJo96b} comes from a wavelet! Since we cannot really
characterize abstractly the representations of $\mathcal{O}_{N}$ coming from
wavelets, we can of course also not find a completely general way of going
the other way, from representations to wavelets. But let us mention some
connections from representations to wavelets which are as direct as possible
with our present technology: if $\ \varphi $ is a father wavelet in $%
L^{2}\left( \mathbb{R}\right) $ satisfying the standard requirements (\ref
{Eq8.1})--(\ref{Eq8.3}) in scale $N$, and $\psi _{1},\dots ,\psi _{N-1}$ are
corresponding mother wavelets as in Theorem \ref{Thm8.1}, then any $\xi \in
L^{2}\left( \mathbb{R}\right) $ has an orthonormal decomposition 
\begin{equation}
\xi \left( \,\cdot \,\right) =\sum_{i=1%
%TCIMACRO{\TeXButton{vphantom}{\vphantom{j,k\in \mathbb{Z}}}}
%BeginExpansion
\vphantom{j,k\in \mathbb{Z}}%
%EndExpansion
}^{N-1}\sum_{j,k\in \mathbb{Z}}a_{jk}^{\left( i\right) }\left( \xi \right)
N^{-\frac{j}{2}}\psi _{i}^{{}}\left( N^{-j}\,\cdot \,-k\right) 
\label{Eq1.35}
\end{equation}
in $L^{2}\left( \mathbb{R}\right) $. In particular $\xi $ is contained in
the closed subspace $\mathcal{V}_{0}$ of $L^{2}\left( \mathbb{R}\right) $
spanned by the translates $\varphi \left( \,\cdot \,-k\right) $, $k\in %
\mathbb{Z}$, of the father wavelets if and only if $a_{jk}^{\left( i\right)
}\left( \xi \right) =0$ for all $j\le 0$, and in that case $\varphi $ has
also a representation 
\begin{equation}
\hat{\xi}\left( t\right) =f\left( t\right) \hat{\varphi}\left( t\right) 
\label{Eq1.36}
\end{equation}
in terms of an $f\in L^{2}\left( \mathbb{T}\right) =L^{2}\left( \mathbb{R}%
\diagup 2\pi \mathbb{Z}\right) $, where $\hat{\ }$ denotes Fourier transform
(\ref{Eq7.6}). See Lemma \ref{Lem10.1} for this. The link between the
representation and the wavelet formulation is then provided by 
\begin{equation}
a_{jk}^{\left( i\right) }\left( \xi \right) =\left(
S_{i}^{*}S_{0}^{*\,j-1}f\right) ^{\sim }\left( k\right)   \label{Eq1.37}
\end{equation}
where $\left( 
%TCIMACRO{\TeXButton{kern}{\mkern9mu}}
%BeginExpansion
\mkern9mu%
%EndExpansion
\right) ^{\sim }$ refers to the Fourier transform on $L^{2}\left( \mathbb{T}%
\right) $: 
\begin{equation}
\tilde{g}\left( k\right) =\frac{1}{2\pi }\int_{-\pi }^{\pi }e^{-ikt}g\left(
t\right) \,dt.  \label{Eq1.38}
\end{equation}
The formula (\ref{Eq1.37}) follows from Corollary \ref{Cor8.3} and Theorem 
\ref{Thm4.2}, and the details of the proof will be given in Corollary \ref
{Cor8.4}. So at least given the father wavelet $\varphi $, the formula (\ref
{Eq1.36})--(\ref{Eq1.37}) give a path from the representation of $\mathcal{O}%
_{N}$ to the wavelet $\psi _{1},\dots ,\psi _{N-1}$. Furthermore, recall
that the father wavelet $\varphi $ under mild regularity assumptions can be
recovered from the function $m_{0}$ via the Mallat algorithm (\ref{Eq9.3})
(see \cite{Mal89}, \cite{Dau92}), i.e., 
\begin{equation}
\hat{\varphi}\left( t\right) =\left( 2\pi \right) ^{-\frac{1}{2}%
}\prod_{k=1}^{\infty }\left( N^{-\frac{1}{2}}m_{0}\left( tN^{-k}\right)
\right) .  \label{Eq1.39}
\end{equation}
But iterating (\ref{Eq1.16}) $n$ times, and applying the scaling operator 
\begin{equation*}
\left( U_{N}\xi \right) \left( t\right) =N^{-\frac{1}{2}}\xi \left(
N^{-1}t\right) ,
\end{equation*}
we see that 
\begin{equation}
\left( U_{N}^{n}S_{0}^{n}\xi \right) \left( t\right) =\prod_{k=1}^{n}\left(
N^{-\frac{1}{2}}m_{0}\left( tN^{-k}\right) \right) \xi \left( t\right) 
\label{Eq1.40}
\end{equation}
when functions in $L^{2}\left( \mathbb{T}\right) $ are viewed as $2\pi $%
-periodic functions on $\mathbb{R}$, and taking the limit $n\rightarrow
\infty $ and using (\ref{Eq1.39}) we see 
\begin{equation}
\lim_{n\rightarrow \infty }\left( U_{N}^{n}S_{0}^{n}\xi \right) \left(
t\right) =\left( 2\pi \right) ^{\frac{1}{2}}\hat{\varphi}\left( t\right) \xi
\left( t\right) .  \label{Eq1.41}
\end{equation}
Thus, when the representation $\left\{ S_{i}\right\} $ of $\mathcal{O}_{N}$
on $L^{2}\left( \mathbb{T}\right) $ is given, the formulae (\ref{Eq1.41}), (%
\ref{Eq1.36}), and (\ref{Eq1.37}) in succession give a prescription for
recovering the multiresolution wavelet theory from the representation.
Similarly 
\begin{equation}
\lim_{n\rightarrow \infty }\left( U_{N}^{n}S_{0}^{n-1}S_{k}^{{}}\xi \right)
\left( t\right) =\left( 2\pi \right) ^{\frac{1}{2}}\hat{\psi}_{k}\left(
t\right) \xi \left( t\right) .  \label{Eq1.42}
\end{equation}

The formulae (\ref{Eq1.37}), (\ref{Eq1.41}), and (\ref{Eq1.42}) were derived
under the assumption that the representation of $\mathcal{O}_{N}$ comes from
a wavelet. More fundamentally, if a representation of $\mathcal{O}_{N}$ is
given, Proposition \ref{Pro1.1} gives a necessary and sufficient condition
that it defines functions $m_{i}:\mathbb{T}\rightarrow \mathbb{C}$ with the
unitarity property (\ref{Eq1.11}). If we further assume that $m_{0}\left(
0\right) =\sqrt{N}$ and $m_{0}$ is Lipschitz continuous at $0$, then the
product expansion (\ref{Eq1.39}) converges and defines the function $\hat{%
\varphi}$. But this is still not sufficient for $\varphi $ to be the father
function of a wavelet, as shown by the example between (6.2.4) and (6.2.5)
in \cite{Dau92}. If 
\begin{equation}
m_{0}\left( t\right) =\sum_{k\in \mathbb{Z}}a_{k}e^{-ikt}  \label{Eq1.41bis}
\end{equation}
is the Fourier expansion of $m_{0}$ with $z=e^{-it}$, put 
\begin{equation}
m_{0}^{(k)}\left( z\right) =m_{0}\left( z\right) m_{0}\left( z^{N}\right)
\cdots m_{0}\left( z^{N^{k-1}}\right) .  \label{Eq1.42bis}
\end{equation}
Assume now also that $m_{0}$ is infinitely differentiable. It is then easy
to show that $\varphi $ is a father function for a wavelet, i.e., (\ref
{Eq8.1})--(\ref{Eq8.3}) are valid, if and only if (10.1)
alone is valid, i.e.,
\begin{equation}
\left\{ \varphi \left( \,\cdot \,-k\right) \right\} _{k\in \mathbb{Z}}\text{%
\quad is an orthonormal set.}  \label{Eq1.43}
\end{equation}
Furthermore, this is again equivalent to any of the following properties (%
\ref{Eq1.44})--(\ref{Eq1.47}). 
%TCIMACRO{
%\TeXButton{Custom Aligned Equation}{\begin{align}
%& \left\| \varphi \right\| _{L^{2}\left( \mathbb{R}\right) }=1.
%\label{Eq1.44} \\
%& \begin{minipage}[t]{\displayboxwidth}\raggedright The probability measures
%$\smash{\left| m_{0}^{(k)}\left( z\right) \right| ^{2}\frac{\left| dz\right| }
%{2\pi }}$ converge weakly to Dirac's delta measure on $1\in
%\mathbb{T}$.\end{minipage}  \label{Eq1.45} \\
%& \begin{minipage}[t]{\displayboxwidth}\raggedright There is a compact set
%$K$ of reals, congruent to $\left[ -\pi ,\pi \right] $ modulo $2\pi $, such
%that $K$ contains $0$ in its interior and $\hat{\varphi }\left( t\right) \ne
%0$ for $t\in K$.\end{minipage}  \label{Eq1.46}
%\end{align}} }
%BeginExpansion
\begin{align}
& \left\| \varphi \right\| _{L^{2}\left( \mathbb{R}\right) }=1.
\label{Eq1.44} \\
& \begin{minipage}[t]{\displayboxwidth}\raggedright The probability measures
$\smash{\left| m_{0}^{(k)}\left( z\right) \right| ^{2}\frac{\left| dz\right| }
{2\pi }}$ converge weakly to Dirac's delta measure on $1\in
\mathbb{T}$.\end{minipage}  \label{Eq1.45} \\
& \begin{minipage}[t]{\displayboxwidth}\raggedright There is a compact set
$K$ of reals, congruent to $\left[ -\pi ,\pi \right] $ modulo $2\pi $, such
that $K$ contains $0$ in its interior and $\hat{\varphi }\left( t\right) \ne
0$ for $t\in K$.\end{minipage}  \label{Eq1.46}
\end{align}%
%EndExpansion
The last condition (\ref{Eq1.46}) is due to A.~Cohen \cite{Coh90}, and the
equivalence of the other two conditions is due to Meyer and Paiva \cite
{MePa93}. The latter paper contains an excellent discussion and also shows
that these conditions are equivalent to $\varphi $ being the unique fixed
point of the map 
\begin{equation}
\psi \rightarrow N^{\frac{1}{2}}\sum_{k=-\infty }^{\infty }a_{k}\psi \left(
N\,\cdot \,+k\right)  \label{Eq1.47}
\end{equation}
among a regular class of functions satisfying $\hat{\psi}\left( 2\pi
k\right) =0$ for $k\in \mathbb{Z}\setminus \left\{ 0\right\} $ and $\hat{\psi%
}\left( 0\right) =\left( 2\pi \right) ^{-\frac{1}{2}}$, this fixed point
being an attractor.

Note finally that the condition (\ref{Eq1.45}) can be translated into the
following necessary and sufficient condition that an $S_{0}$ of the form (%
\ref{Eq1.16}) comes from the father function of a wavelet: 
\begin{equation}
\lim_{n\rightarrow \infty }\ip{S_{0}^{n}\openone }{M_{f}^{}S_{0}^{n}\openone
}=f\left( 0\right)  \label{Eq1.48}
\end{equation}
for all $f\in C\left( \mathbb{T}\right) =C\left( \mathbb{R}\diagup 2\pi %
\mathbb{Z}\right) $, where $M_{f}$ denotes the operator of multiplication by 
$f$ on $L^{2}\left( \mathbb{T}\right) $. See Section \ref{Fath} for details.

\section{\label{Coho}%
%TCIMACRO{\TeXButton{frogkern}{\protect\frogkern } }
%BeginExpansion
\protect\frogkern %
%EndExpansion
Cohomology of the map $z\mapsto z^{N}$ with values%
%TCIMACRO{\TeXButton{frogleg}{\protect\frogleg } }
%BeginExpansion
\protect\frogleg %
%EndExpansion
in a topological group $G$}

The terminology we introduce in this section will be used in two
connections. In Section \ref{WoldSm}, with $G=\mathbb{T}$, it will be used
in the characterization of the unitary part of the Wold decomposition of $%
S_{0}$. In Section \ref{Clas}, also with $G=\mathbb{T}$, it will be used in
the characterization of unitary equivalence of two diagonal representations.
In Sections \ref{WoldSc} and \ref{WoldT} somewhat similar terminology, but
with a more general notion of cohomology which is less direct to formulate
in the abstract framework, will be used in the case $G=U\left( N\right) $;
see for example (\ref{EqX.15}). $G$ is a topological group throughout this
section.

Let $\Omega $ be a measure space, $\sigma :\Omega \rightarrow \Omega $ an
endomorphism, and $\mu $ a probability measure on $\Omega $ such that $\mu
\left( \sigma ^{-1}\left( Y\right) \right) =\mu \left( Y\right) $ for all
measurable $Y\subset \Omega $. Extending the terminology in \cite{CFS82}
from the case of automorphisms to the case of endomorphisms, we may define a 
\emph{cocycle} for $\sigma $ with values in $G$ as a map $c:\Omega \times %
\mathbb{N}\rightarrow G$ such that 
\begin{equation}
c\left( x,m+n\right) =c\left( x,m\right) c\left( \sigma ^{m}\left( x\right)
,n\right) .  \label{Eq2.1}
\end{equation}
But then it follows by induction that 
\begin{equation}
c\left( x,m\right) =\begin{cases} 1 &\text{if $m=0$} \\ \prod
_{k=0}^{m-1}\sigma \left( \sigma ^{k}\left( x\right) ,1\right) &\text{if
$m>0$} \end{cases}  \label{Eq2.2}
\end{equation}
so we may and will simply consider a cocycle to be a measurable map $c\left(
\,\cdot \,\right) =c\left( \,\cdot \,,1\right) $ from $\Omega $ into $G$.
Any such map defines a proper cocycle through the formula above.

We say that two cocycles $c_{1}$, $c_{2}$ are \emph{cohomological} if there
is a function $\Delta :\Omega \rightarrow G$ such that 
\begin{equation}
c_{1}\left( x\right) =\Delta \left( \sigma \left( x\right) \right)
^{-1}c_{2}\left( x\right) \Delta \left( x\right)  \label{Eq2.3}
\end{equation}
or 
\begin{equation}
c_{1}\left( x,m\right) =\Delta \left( \sigma ^{m}\left( x\right) \right)
^{-1}c_{2}\left( x,m\right) \Delta \left( x\right) .  \label{Eq2.4}
\end{equation}
In the case that $G$ is abelian, this is the same as saying that $c_{1}$ and 
$c_{2}$ \emph{cobound,} i.e., that $c_{1}\left( x\right) c_{2}\left(
x\right) ^{-1}$ is a coboundary. We also say that $c_{1}$ and $c_{2}$ are
cohomologous. We say in general that a cocycle $c$ is a \emph{coboundary} if
there is another cocycle $\Delta $ such that 
\begin{equation}
c\left( x\right) =\Delta \left( x\right) \Delta \left( \sigma \left(
x\right) \right) ^{-1}  \label{Eq2.5}
\end{equation}
or 
\begin{equation}
c\left( x,m\right) =\Delta \left( x\right) \Delta \left( \sigma ^{m}\left(
x\right) \right) ^{-1}.  \label{Eq2.6}
\end{equation}
In the case that $G$ is abelian, we see that the relation of cohomology is
an equivalence relation.

The question of which cocycles are coboundaries is in general a difficult
one. Recall for example from \cite[Theorem 6.1]{Jor95} that if $c:\mathbb{T}%
\rightarrow \mathbb{T}$ is a Hardy function (i.e., $c\in H^{\infty }\left( %
\mathbb{T}\right) $), and $\sigma \left( z\right) =z^{2}$ for $z\in %
\mathbb{T}$, then the equation 
\begin{equation}
c\left( z\right) \Delta \left( z^{2}\right) =\Delta \left( z\right) 
\label{Eq2.7}
\end{equation}
has a nonzero solution $\Delta \in L^{2}\left( \mathbb{T}\right) $ if and
only if $c$ is a monomial, $c\left( z\right) =z^{n}$, and then $\Delta
\left( z\right) =dz^{-n}$ for a constant $d$. Another criterion which is
more indirect is in \cite[Corollary 3]{Wal96}. The version of this corollary
which is interesting for us is the following: if $c$ is a measurable cocycle
for $z\mapsto z^{N}$ with values in $\mathbb{T}$, then the following
conditions (\ref{Eq2.8}) and (\ref{Eq2.9}) are equivalent. 
%TCIMACRO{
%\TeXButton{Custom Aligned Equation}{\begin{align}
%& \begin{minipage}[t]{\displayboxwidth}\raggedright There is an $f \in
%L^{\infty }\left( \mathbb{T}\right) $ such that the sequence has a nonzero
%w*-limit point as $m\rightarrow \infty $. \end{minipage}  \label{Eq2.8}
%\\
%& \begin{minipage}[t]{\displayboxwidth}\raggedright The cocycle $c$ is a
%coboundary. \end{minipage}  \label{Eq2.9}
%\end{align}} }
%BeginExpansion
\begin{align}
& \begin{minipage}[t]{\displayboxwidth}\raggedright There is an $f \in
L^{\infty }\left( \mathbb{T}\right) $ such that the sequence has a nonzero
w*-limit point as $m\rightarrow \infty $. \end{minipage}  \label{Eq2.8}
\\
& \begin{minipage}[t]{\displayboxwidth}\raggedright The cocycle $c$ is a
coboundary. \end{minipage}  \label{Eq2.9}
\end{align}%
%EndExpansion
Furthermore, if these conditions are fulfilled, the cocycle $\Delta $ having 
$c$ as coboundary is unique up to a phase factor, and 
\begin{equation*}
\lim_{m\rightarrow \infty }\frac{1}{m}\sum_{k=0}^{m-1}f\left(
z^{N^{k}}\right) \bar{c}\left( z^{N^{k-1}}\right) \cdots \bar{c}\left(
z\right) =\bar{\Delta}\left( z\right) \int_{\mathbb{T}}f\left( \eta \right)
\Delta \left( \eta \right) \frac{\left| d\eta \right| }{2\pi }
\end{equation*}
in $L^{2}\left( \mathbb{T}\right) $, and also pointwise for almost all $z$,
for all $f\in L^{\infty }\left( \mathbb{T}\right) $.

The main input in the proof is of course Birkhoff's ergodic theorem, which
immediately gives the implication from (\ref{Eq2.9}) to the conclusion.

\section{\label{WoldSm}%
%TCIMACRO{\TeXButton{frogkern}{\protect\frogkern } }
%BeginExpansion
\protect\frogkern %
%EndExpansion
The Wold decomposition of isometries $S_{m}$ of the form $\left( S_{m}%
%TCIMACRO{\TeXButton{xi}{\protect\xi } }
%BeginExpansion
\protect\xi %
%EndExpansion
\right) \left( z\right) =m\left( z\right) 
%TCIMACRO{\TeXButton{xi}{\protect\xi } }
%BeginExpansion
\protect\xi %
%EndExpansion
\left( z^{N}\right) $}

Equip the circle $\mathbb{T}$ with Haar measure $\frac{\left| dz\right| }{%
2\pi }$, and let $N\in \left\{ 2,3,\dots \right\} $. Formula (\ref{Eq1.3})
now takes the form 
\begin{align}
\int_{\mathbb{T}}f\left( z\right) \frac{\left| dz\right| }{2\pi }& =\int_{%
\mathbb{T}}\frac{1}{N}\sum_{%
\substack{ w \\ \makebox[0pt]{\hss
{$\scriptstyle w^{N}=z $}\hss }}}f\left( w\right) \frac{\left| dz\right| }{%
2\pi }  \label{Eq3.0} \\
& =\int_{\mathbb{T}}f\left( z^{N}\right) \frac{\left| dz\right| }{2\pi }. 
\notag
\end{align}
Let $m:\mathbb{T}\rightarrow \mathbb{C}$ be a measurable function. Define $%
S_{m}:L^{2}\left( \mathbb{T}\right) \rightarrow L^{2}\left( \mathbb{T}%
\right) $ by 
\begin{equation}
\left( S_{m}\xi \right) \left( z\right) =m\left( z\right) \xi \left(
z^{N}\right) .  \label{Eq3.1}
\end{equation}
We have already computed in (\ref{Eq1.17}) that 
\begin{equation}
\left( S_{m}^{*}\xi \right) \left( z\right) =\frac{1}{N}\sum_{\substack{ w
\\ \makebox[0pt]{\hss {$\scriptstyle w^{N}=z $}\hss }}}\bar{m}\left(
w\right) \xi \left( w\right)  \label{Eq3.2}
\end{equation}
and hence 
%TCIMACRO{
%\TeXButton{parameters}{\newlength{\widthsumwithsubstack}
%\newlength{\widthsumalone}
%\newlength{\widthsumleftsubstack}
%\settowidth{\widthsumwithsubstack}{$\displaystyle 
%\sum _{\substack{ w \\ w^{N}=z }}$}
%\settowidth{\widthsumalone}{$\displaystyle \sum $}
%\setlength{\widthsumleftsubstack}{0.5\widthsumwithsubstack}
%\addtolength{\widthsumleftsubstack}{0.5\widthsumalone}}}
%BeginExpansion
\newlength{\widthsumwithsubstack}
\newlength{\widthsumalone}
\newlength{\widthsumleftsubstack}
\settowidth{\widthsumwithsubstack}{$\displaystyle 
\sum _{\substack{ w \\ w^{N}=z }}$}
\settowidth{\widthsumalone}{$\displaystyle \sum $}
\setlength{\widthsumleftsubstack}{0.5\widthsumwithsubstack}
\addtolength{\widthsumleftsubstack}{0.5\widthsumalone}%
%EndExpansion
\begin{equation}
\left( S_{m}^{*}S_{m}^{{}}\xi \right) \left( z\right) =\frac{1}{N}\left( 
%TCIMACRO{
%\TeXButton{Sum}{\makebox[\widthsumleftsubstack][l]{$\displaystyle 
%\sum_{\substack{ w \\ w^{N}=z  }}$\hss }} }
%BeginExpansion
\makebox[\widthsumleftsubstack][l]{$\displaystyle 
\sum_{\substack{ w \\ w^{N}=z  }}$\hss }%
%EndExpansion
\left| m\left( w\right) \right| ^{2}\right) \xi \left( z\right) .
\label{Eq3.3}
\end{equation}
It follows immediately from this spectral representation of $%
S_{m}^{*}S_{m}^{{}}$ that $S_{m}$ is bounded if and only if $m\in L^{\infty
}\left( \mathbb{T}\right) $ and then 
\begin{equation}
\left\| S_{m}\right\| ^{2}=\frac{1}{N}\mathop{\mathrm{ess\,sup}}_{z\in %
\mathbb{T}}\left( 
%TCIMACRO{
%\TeXButton{Sum}{\makebox[\widthsumleftsubstack][l]{$\displaystyle 
%\sum_{\substack{ w \\ w^{N}=z  }}$\hss }} }
%BeginExpansion
\makebox[\widthsumleftsubstack][l]{$\displaystyle 
\sum_{\substack{ w \\ w^{N}=z  }}$\hss }%
%EndExpansion
\left| m\left( w\right) \right| ^{2}\right) \le \left\| m\right\| _{\infty
}^{2}  \label{Eq3.4}
\end{equation}

Furthermore, we see directly from the spectral representation (\ref{Eq3.3})
that $S_{m}$ is an isometry if and only if 
\begin{equation}
\frac{1}{N}\sum_{%
\substack{ w \\ \makebox[0pt]{\hss {$\scriptstyle w^{N}=z
$}\hss }}} \left| m\left( w\right) \right| ^{2}=1  \label{Eq3.5}
\end{equation}
for almost all $z$.

In general, if $S$ is an isometry, define a decreasing sequence of
projections by 
\begin{equation}
E_{k}=S^{k}S^{*\,k}  \label{Eq3.6}
\end{equation}
and let 
\begin{equation}
P_{U}=\mathop{\text{\textup{s-lim}}}_{k\rightarrow \infty }E_{k}.
\label{Eq3.7}
\end{equation}
Then $SP_{U}=P_{U}S$, $P_{U}S$ is a unitary operator on $P_{U}\mathcal{H}$,
and $\left( 1-P_{U}\right) S$ is a shift on $\left( 1-P_{U}\right) \mathcal{H%
}$, i.e., 
\begin{equation}
\bigcap_{n}S^{n}\left( 1-P_{U}\right) \mathcal{H}=\left\{ 0\right\} .
\label{Eq3.8}
\end{equation}
(Note that the two-sided shift is not a shift with this terminology.) The
decomposition 
\begin{equation}
S=SP_{U}\oplus S\left( 1-P_{U}\right)   \label{Eq3.9}
\end{equation}
is the so-called Wold decomposition of $S$ into a unitary operator and a
shift. (For more details on the general Wold decomposition, and some of its
applications, the reader is referred to \cite{SzFo70}, which also serves as
an excellent background reference for the operator theory used in the
present paper.) For $S_{m}$ given by (\ref{Eq3.1}), a calculation now shows
that 
\begin{equation}
\left( E_{k}\xi \right) \left( z\right) =m^{(k)}\left( z\right) \frac{1}{%
N^{k}}\sum_{%
\substack{ w \\ \makebox[0pt]{\hss {$\scriptstyle
w^{N^{k}}=z^{N^{k}} $}\hss }}}\bar{m}^{(k)}\left( w\right) \xi \left(
w\right)   \label{Eq3.10}
\end{equation}
where 
\begin{equation}
m^{(k)}\left( z\right) =\prod_{j=0}^{k-1}m\left( z^{N^{j}}\right) .
\label{Eq3.11}
\end{equation}

Our main result on the Wold decomposition of $S_{m}$ is the following:

\begin{theorem}
\label{Thm3.1}The projection $P_{U}$ corresponding to the unitary part of
the Wold decomposition of the isometry $S_{m}$ is one- or zero-dimensional.
Furthermore, $P_{U}$ is one-dimensional if and only if both conditions 
\textup{(\ref{Eq3.12})} and \textup{(\ref{Eq3.13})} are fulfilled. 
%TCIMACRO{
%\TeXButton{Custom Aligned Equation}{\begin{align}
%& \begin{minipage}[t]{\displayboxwidth}\raggedright $\left| m\left(
%z\right) \right| =1$ for almost all $z\in \mathbb{T}$. \end{minipage}
%\label{Eq3.12} \\
%& \begin{minipage}[t]{\displayboxwidth}\raggedright There exists a
%measurable function $\xi :\mathbb{T}\rightarrow \mathbb{T}$ and a $\lambda
%\in \mathbb{T}$ such that $m\left( z\right) \xi \left( z^{N}\right) =\lambda
%\xi \left( z\right) $ for almost all $z\in \mathbb{T}$. \end{minipage}  
%\label{Eq3.13}
%\end{align}} }
%BeginExpansion
\begin{align}
& \begin{minipage}[t]{\displayboxwidth}\raggedright $\left| m\left(
z\right) \right| =1$ for almost all $z\in \mathbb{T}$. \end{minipage}
\label{Eq3.12} \\
& \begin{minipage}[t]{\displayboxwidth}\raggedright There exists a
measurable function $\xi :\mathbb{T}\rightarrow \mathbb{T}$ and a $\lambda
\in \mathbb{T}$ such that $m\left( z\right) \xi \left( z^{N}\right) =\lambda
\xi \left( z\right) $ for almost all $z\in \mathbb{T}$. \end{minipage}  
\label{Eq3.13}
\end{align}%
%EndExpansion
In this case the range of the projection $P_{U}$ is $\mathbb{C}\xi \subset
L^{2}\left( \mathbb{T}\right) $.

In short, $S_{m}$ is a shift if and only if there exists no phase factor $%
\lambda $ such that $\bar{\lambda}m_{0}$ is a coboundary for the $z\mapsto
z^{N}$ action with values in $\mathbb{T}$.
\end{theorem}

In order to prove Theorem \ref{Thm3.1}, it will be useful to work with the
root mean operator $R=R_{m}$ defined on measurable functions $\xi :\mathbb{T}%
\rightarrow \mathbb{C}$ as follows: 
\begin{equation}
\left( R\xi \right) \left( z\right) =\frac{1}{N}\sum_{%
\substack{ w \\ \makebox[0pt]{\hss {$\scriptstyle w^{N}=z
$}\hss }}}\left| m\left( w\right) \right| ^{2}\xi \left( w\right) .
\label{Eq3.14}
\end{equation}
It follows immediately from (\ref{Eq3.5}) that $R$ is bounded as an operator
from $L^{p}$ $\left( \mathbb{T}\right) $ into $L^{p}\left( \mathbb{T}\right) 
$ for $1\le p\le \infty $, and also $R$ preserves positive functions and,
from (\ref{Eq3.5}), 
\begin{equation}
R\openone=\openone.  \label{Eq3.15}
\end{equation}
Thus 
\begin{equation}
\left\| R\right\| _{\infty \rightarrow \infty }=1  \label{Eq3.16}
\end{equation}
and a computation like the one after (\ref{Eq1.8})--(\ref{Eq1.9}) shows 
\begin{equation}
\left( R^{*}\xi \right) \left( z\right) =\left| m\left( z\right) \right|
^{2}\xi \left( z^{N}\right) .  \label{Eq3.17}
\end{equation}
If $f\in L^{\infty }\left( \mathbb{T}\right) $, again let $M_{f}:L^{2}\left( %
\mathbb{T}\right) \rightarrow L^{2}\left( \mathbb{T}\right) $ denote the
operation of multiplication by $f$, 
\begin{equation}
\left( M_{f}\xi \right) \left( z\right) =f\left( z\right) \xi \left( z\right)
\label{Eq3.18}
\end{equation}
for $\xi \in L^{2}\left( \mathbb{T}\right) $.

We will need the formula 
\begin{align}
E_{k}M_{f}E_{k}& =M_{\left( R^{k}f\right) \left( z^{N^{k}}\right) }E_{k}
\label{Eq3.19} \\
& =M_{\left( R^{k}f\right) \circ \sigma ^{k}}E_{k}  \notag
\end{align}
which follows from (\ref{Eq3.10}) and (\ref{Eq3.14}) by the following
computation: 
%TCIMACRO{
%\TeXButton{parameters}{\newlength{\wnkwidth}
%\settowidth{\wnkwidth}{$\scriptstyle w^{N^{k}}$}}}
%BeginExpansion
\newlength{\wnkwidth}
\settowidth{\wnkwidth}{$\scriptstyle w^{N^{k}}$}%
%EndExpansion
\begin{align*}
\left( E_{k}M_{f}E_{k}\xi \right) \left( z\right) & =m^{\left( k\right)
}\left( z\right) \frac{1}{N^{k}}\sum_{%
\substack{ w \\ \makebox[0pt]{\hss
{$\scriptstyle w^{N^{k}}=z^{N^{k}}$}\hss } }}\bar{m}^{\left( k\right)
}\left( w\right) f\left( w\right) \left( E_{k}\xi \right) \left( w\right) \\
& =m^{(k)}\left( z\right) \frac{1}{\left( N^{k}\right) ^{2}}\sum_{%
\substack{ w,v \\ \makebox[0pt]{\hss {$\scriptstyle
w^{N^{k}}=z^{N^{k}}$}\hss } \\ \makebox[0pt]{\hss {$\scriptstyle 
\makebox[\wnkwidth]{}=v^{N^{k}}$}\hss }}%
}\bar{m}^{(k)}\left( w\right) f\left( w\right) m^{\left( k\right) }\left(
w\right) \bar{m}^{(k)}\left( v\right) \xi \left( v\right) \\
& =\left( \frac{1}{N^{k}}\sum_{%
\substack{ w \\ \makebox[0pt]{\hss
{$\scriptstyle w^{N^{k}}=z^{N^{k}} $}\hss }}}\left| m^{\left( k\right)
}\left( w\right) \right| ^{2}f\left( w\right) \right) \left( \frac{1}{N^{k}}%
m^{(k)}\left( z\right) \sum_{%
\substack{ v \\ \makebox[0pt]{\hss
{$\scriptstyle v^{N^{k}}=z^{N^{k}} $}\hss }}}\bar{m}^{\left( k\right)
}\left( v\right) \xi \left( v\right) \right) \\
& =\left( R^{k}f\right) \left( z^{N^{k}}\right) \left( E_{k}\xi \right)
\left( z\right) .
\end{align*}

\begin{lemma}
\label{Lem3.2}Assume that $P_{U}\neq 0$ and pick $\xi \in P_{U}L^{2}\left( %
\mathbb{T}\right) $ such that $\left\| \xi \right\| _{2}=1$. It follows that 
\begin{equation}
\left| \xi \left( z\right) \right| ^{2}=\lim_{k\rightarrow \infty
}\prod_{n=0}^{k}\left| m\left( z^{N^{n}}\right) \right| ^{2}  \label{Eq3.20}
\end{equation}
and 
\begin{equation}
\left| \xi \left( z\right) \right| =1=\left| m\left( z\right) \right| 
\label{Eq3.21}
\end{equation}
for almost all $z\in \mathbb{T}$.
\end{lemma}

%TCIMACRO{\TeXButton{Begin Proof}{\begin{proof}}}
%BeginExpansion
\begin{proof}%
%EndExpansion
As $P_{U}\xi =\xi $ and $P_{U}\le E_{k}$, we have $E_{k}\xi =\xi $ for all $%
k\in \mathbb{N}$, and using (\ref{Eq3.19}) on an arbitrary $f\in L^{\infty
}\left( \mathbb{T}\right) $ we have 
\begin{align*}
\int_{\mathbb{T}}f\left( z\right) \left| \xi \left( x\right) \right| ^{2}%
\frac{\left| dz\right| }{2\pi }& =\ip{\xi }{\smash{M_f\xi }} \\
& =\ip{E_k\xi }{\smash{M_fE_k\xi }} \\
& =\ip{\xi }{\smash{E_kM_fE_k\xi }} \\
& =\ip{\xi }{\smash{M_{\left( R^kf\right) \left( z^{N^k}\right) }\xi }} \\
& =\int_{\mathbb{T}}\left| \xi \left( z\right) \right| ^{2}\left(
R^{k}f\right) \left( z^{N^{k}}\right) \frac{\left| dz\right| }{2\pi } \\
& =\int_{\mathbb{T}}\frac{1}{N^{k}}\sum_{%
\substack{ w \\ \makebox[0pt]{\hss
{$\scriptstyle w^{N^{k}}=z $}\hss }}}\left| \xi \left( w\right) \right|
^{2}\left( R^{k}f\right) \left( z\right) \frac{\left| dz\right| }{2\pi }.
\end{align*}
Now use 
\begin{equation}
\left( R^{k}f\right) \left( z\right) =\frac{1}{N^{k}}\sum_{\substack{ v \\
\makebox[0pt]{\hss {$\scriptstyle v^{N^{k}}=z $}\hss }}}\left| m^{\left(
k\right) }\left( v\right) \right| ^{2}f\left( v\right)  \label{Eq3.22}
\end{equation}
to compute further 
\begin{align*}
\int_{\mathbb{T}}f\left( z\right) \left| \xi \left( z\right) \right| ^{2}%
\frac{\left| dz\right| }{2\pi }& =\int_{\mathbb{T}}\frac{1}{N^{2k}}\sum_{
\substack{ w,v \\ \makebox[0pt]{\hss {$\scriptstyle w^{N^{k}}=z$}\hss } \\
\makebox[0pt]{\hss {$\scriptstyle v^{N^{k}}=z$}\hss }}}\left| \xi \left(
w\right) \right| ^{2}\left| m^{(k)}\left( v\right) \right| ^{2}f\left(
v\right) \frac{\left| dz\right| }{2\pi } \\
& =\int_{\mathbb{T}}\left| m^{(k)}\left( z\right) \right| ^{2}\frac{1}{N^{k}}%
\sum_{%
\substack{ w \\ \makebox[0pt]{\hss {$\scriptstyle
w^{N^{k}}=z^{N^{k}} $}\hss }}}\left| \xi \left( w\right) \right| ^{2}f\left(
z\right) \frac{\left| dz\right| }{2\pi },
\end{align*}
where the last equality follows from the general formula 
\begin{equation}
\int_{\mathbb{T}}g\left( z\right) \frac{1}{N^{k}}\sum_{\substack{ v \\
\makebox[0pt]{\hss {$\scriptstyle v^{N^{k}}=z $}\hss }}}h\left( v\right) 
\frac{\left| dz\right| }{2\pi }=\int_{\mathbb{T}}g\left( z^{N^{k}}\right)
h\left( z\right) \frac{\left| dz\right| }{2\pi }.  \label{Eq3.23}
\end{equation}
We conclude that 
\begin{equation}
\int_{\mathbb{T}}f\left( z\right) \left| \xi \left( z\right) \right| ^{2}%
\frac{\left| dz\right| }{2\pi }=\int_{\mathbb{T}}f\left( z\right) \left|
m^{(k)}\left( z\right) \right| ^{2}\frac{1}{N^{k}}\sum_{%
\substack{ w \\ \makebox[0pt]{\hss {$\scriptstyle w^{N^{k}}=z^{N^{k}} $}\hss
}}}\left| \xi \left( w\right) \right| ^{2}\frac{\left| dz\right| }{2\pi }.
\label{Eq3.24}
\end{equation}
As this equality is valid for all $f\in L^{\infty }\left( \mathbb{T}\right) $%
, we conclude that 
\begin{equation}
\left| \xi \left( z\right) \right| ^{2}=\left| m^{(k)}\left( z\right)
\right| ^{2}\frac{1}{N^{k}}\sum_{%
\substack{ w \\ \makebox[0pt]{\hss {$\scriptstyle w^{N^{k}}=z^{N^{k}}
$}\hss }}}\left| \xi \left( w\right) \right| ^{2}  \label{Eq3.25}
\end{equation}
for almost all $z\in \mathbb{T}$, $k=1,2,\dots $.

Now let $R_{1}$ be the root mean operator on the measurable functions on $%
\mathbb{T}$ defined by putting $m=1$ in (\ref{Eq3.14}), i.e., 
\begin{equation}
\left( R_{1}\xi \right) \left( z\right) =\frac{1}{N}\sum_{\substack{ w \\
\makebox[0pt]{\hss {$\scriptstyle w^{N}=z $}\hss }}}\xi \left( w\right) .
\label{Eq3.26}
\end{equation}
If $\varphi \in L^{1}\left( \mathbb{T}\right) $, it follows by approximating 
$\varphi $ by functions in $C\left( \mathbb{T}\right) $ that 
\begin{equation}
\lim_{k\rightarrow \infty }\left\| R_{1}^{k}\left( \varphi \right) -\left(
\int_{\mathbb{T}}\varphi \left( z\right) \frac{\left| dz\right| }{2\pi }%
\right) \openone\right\| _{1}=0,  \label{Eq3.27}
\end{equation}
and hence the sequence $R_{1}^{k}\left( \varphi \right) $ converges to the
constant function $\int_{\mathbb{T}}\varphi \left( z\right) \frac{\left|
dz\right| }{2\pi }$ in measure, i.e., 
\begin{equation*}
\lim_{k\rightarrow \infty }\mu \left\{ z\in \mathbb{T}\biggm|\left|
R_{1}^{k}\left( \varphi \right) \left( z\right) -\int_{\mathbb{T}}\varphi
\left( \eta \right) \frac{\left| d\eta \right| }{2\pi }\right| >\varepsilon
\right\} =0
\end{equation*}
for all $\varepsilon >0$.

But repeating the proof of Birkhoff's mean ergodic theorem \cite{CFS82,Wal82}%
, one can show the stronger conclusion that $R_{1}^{k}\left( \varphi \right) 
$ converges almost everywhere to a function which is invariant under all $N$%
-adic rotations, and therefore under all rotations, i.e., $R_{1}^{k}\left(
\varphi \right) $ converges almost everywhere to the constant $\int_{%
\mathbb{T}}\varphi \left( z\right) \frac{\left| dz\right| }{2\pi }$. But (%
\ref{Eq3.25}) says that 
\begin{equation}
\left| \xi \left( z\right) \right| ^{2}=\left| m^{(k)}\left( z\right)
\right| ^{2}R_{1}^{k}\left( \left| \xi \right| ^{2}\right) \left(
z^{N^{k}}\right) ,  \label{Eq3.28}
\end{equation}
and $R_{1}^{k}\left( \left| \xi \right| ^{2}\right) \left( z^{N^{k}}\right)
\rightarrow \left\| \xi \right\| _{2}^{2}=1$ for almost all $z$ by the
remarks above, and hence we have proved (\ref{Eq3.20}): 
\begin{equation*}
\lim_{k\rightarrow \infty }\left| m^{(k)}\left( z\right) \right| ^{2}=\left|
\xi \left( z\right) \right| ^{2}
\end{equation*}
for almost all $z.$ In particular the limit to the left exists for almost
all $z$. Put 
\begin{equation}
m_{\infty }\left( z\right) =\lim_{k\rightarrow \infty }\left| m^{\left(
k\right) }\left( z\right) \right| .  \label{Eq3.29}
\end{equation}
One consequence of (\ref{Eq3.20}) is that, if $\xi \in P_{U}L^{2}\left( %
\mathbb{T}\right) $, then 
\begin{equation}
\left| \xi \left( z\right) \right| =\left\| \xi \right\| _{2}m_{\infty
}\left( z\right)  \label{Eq3.30}
\end{equation}
for almost all $z$; and this immediately implies that the space $%
P_{U}L^{2}\left( \mathbb{T}\right) $ is one-dimensional, establishing the
first statement of Theorem \ref{Thm3.1}.

Now, from the relation 
\begin{equation}
m^{\left( k+1\right) }\left( z\right) =m\left( z\right) m^{\left( k\right)
}\left( z^{N}\right)  \label{Eq3.31}
\end{equation}
and (\ref{Eq3.29}) we deduce 
\begin{equation}
m_{\infty }\left( z\right) =\left| m\left( z\right) \right| m_{\infty
}\left( z^{N}\right) .  \label{Eq3.32}
\end{equation}
But, using this and (\ref{Eq3.5}), we further deduce that 
\begin{align*}
\sum_{\substack{ w \\ \makebox[0pt]{\hss {$\scriptstyle w^{N}=z $}\hss }}%
}m_{\infty }\left( w\right) ^{2}& =\sum_{%
\substack{ w \\ \makebox[0pt]{\hss
{$\scriptstyle w^{N}=z $}\hss }}}\left| m\left( w\right) \right|
^{2}m_{\infty }\left( z\right) ^{2} \\
& =Nm_{\infty }\left( z\right) ^{2},
\end{align*}
so 
\begin{align*}
m_{\infty }\left( z\right) ^{2}& =\frac{1}{N}\sum_{%
\substack{ w \\ \makebox[0pt]{\hss {$\scriptstyle w^{N}=z
$}\hss }}}m_{\infty }\left( w\right) ^{2} \\
& =R_{1}^{{}}\left( m_{\infty }^{2}\right) \left( z\right) .
\end{align*}
Iterating this, we obtain 
\begin{equation*}
m_{\infty }\left( z\right) ^{2}=R_{1}^{k}\left( m_{\infty }^{2}\right)
\left( z\right)
\end{equation*}
for $k=1,2,3$; and, letting $k\rightarrow \infty $, 
\begin{equation*}
m_{\infty }\left( z\right) ^{2}=\int_{\mathbb{T}}m_{\infty }\left( w\right)
^{2}\frac{\left| dw\right| }{2\pi }.
\end{equation*}
Thus $m_{\infty }\left( z\right) $ is a positive constant, and re-inserting
this in (\ref{Eq3.32}) gives 
\begin{equation*}
\left| m\left( z\right) \right| =1
\end{equation*}
for almost all $z$. Thus from (\ref{Eq3.29}), 
\begin{equation*}
m_{\infty }\left( z\right) =1
\end{equation*}
for almost all $z$, and then from (\ref{Eq3.20}), 
\begin{equation*}
\left| \xi \left( z\right) \right| =1
\end{equation*}
for almost all $z$. This ends the proof of (\ref{Eq3.21}) and thus of Lemma 
\ref{Lem3.2}.%
%TCIMACRO{\TeXButton{End Proof}{\end{proof}}}
%BeginExpansion
\end{proof}%
%EndExpansion

%TCIMACRO{
%\TeXButton{Begin Proof of Theorem 3.1}{\begin{proof}[Proof of Theorem 
%\textup{\ref{Thm3.1}}]}}
%BeginExpansion
\begin{proof}[Proof of Theorem \textup{\ref{Thm3.1}}]%
%EndExpansion
We already commented in connection with (\ref{Eq3.30}) that if $P_{U}\ne 0$,
then $P_{U}$ is one-dimensional, and if $P_{U}\ne 0$, then (\ref{Eq3.12})
follows from (\ref{Eq3.21}). But if $P_{U}L^{2}\left( \mathbb{T}\right) =%
\mathbb{C}\xi $ with $\left\| \xi \right\| _{2}=1$, it follows from
unitarity of $S_{m}P_{U}=P_{U}S_{m}$ that $\xi $ must be an eigenvector of $%
S_{m}$ with eigenvalue $\lambda $ of modulus one, $S_{m}\xi $ =$\lambda \xi $%
, or 
\begin{equation*}
m\left( z\right) \xi \left( z^{N}\right) =\lambda \xi \left( z\right) ,
\end{equation*}
which is (\ref{Eq3.13}).

Conversely, if (\ref{Eq3.12}) and (\ref{Eq3.13}) are fulfilled, it is
obvious that $P_{U}\ne 0$, since $\xi \in P_{U}L^{2}\left( \mathbb{T}\right) 
$ (it suffices instead of (\ref{Eq3.12}) and (\ref{Eq3.13}) merely to assume
that $S_{m}$ has an eigenvector with eigenvalues of modulus $1$).%
%TCIMACRO{\TeXButton{End Proof}{\end{proof}}}
%BeginExpansion
\end{proof}%
%EndExpansion

\section{\label{WoldSc}%
%TCIMACRO{\TeXButton{frogkern}{\protect\frogkern } }
%BeginExpansion
\protect\frogkern %
%EndExpansion
The Wold decomposition of operators $S_{C}$ on $L^{2}\left( \mathbb{T}; 
\mathbb{C}^{n}\right) $ of the form $\left( S_{C}%
%TCIMACRO{\TeXButton{xi}{\protect\xi } }
%BeginExpansion
\protect\xi %
%EndExpansion
\right) \left( z\right) =C\left( z\right) 
%TCIMACRO{\TeXButton{xi}{\protect\xi } }
%BeginExpansion
\protect\xi %
%EndExpansion
\left( z^{N}\right) $}

In this section we will consider a situation which is more general in some
respects, and more special in other respects, than in Section \ref{WoldSm}.
Let $L^{2}\left( \mathbb{T};\mathbb{C}^{n}\right) \cong L^{2}\left( %
\mathbb{T}\right) \otimes \mathbb{C}^{n}$ be the Hilbert space of $L^{2}$%
-functions on $\mathbb{T}$ with values in the Hilbert space $\mathbb{C}^{n}$%
. Let $C:\mathbb{T}\rightarrow M_{n}=\mathcal{B}\left( \mathbb{C}^{n}\right) 
$ be a measurable bounded function, and define an operator $S_{C}\in 
\mathcal{B}\left( L^{2}\left( \mathbb{T};\mathbb{C}^{n}\right) \right) $ by 
\begin{equation}
\left( S_{C}\xi \right) \left( z\right) =C\left( z\right) \xi \left(
z^{N}\right) .  \label{EqX.1}
\end{equation}
One verifies as in Section \ref{WoldSm} that 
\begin{equation}
\left( S_{C}^{*}\xi \right) \left( z\right) =\frac{1}{N}\sum_{%
\substack{ w \\ \makebox[0pt]{\hss {$\scriptstyle w^{N}=z
$}\hss }}}C\left( w\right) ^{*}\xi \left( w\right)  \label{EqX.2}
\end{equation}
and hence 
\begin{equation}
\left( S_{C}^{*}S_{C}^{{}}\xi \right) \left( z\right) =\frac{1}{N}\sum_{%
\substack{ w \\ \makebox[0pt]{\hss {$\scriptstyle w^{N}=z
$}\hss }}}C\left( w\right) ^{*}C\left( w\right) \xi \left( z\right) .
\label{EqX.3}
\end{equation}
Thus $S_{C}$ is an isometry if and only if 
\begin{equation}
\frac{1}{N}\sum_{%
\substack{ w \\ \makebox[0pt]{\hss {$\scriptstyle w^{N}=z
$}\hss }}}C\left( w\right) ^{*}C\left( w\right) =\openone_{n}  \label{EqX.4}
\end{equation}
for almost all $z\in \mathbb{T}$. So far everything generalizes Section \ref
{WoldSm}, but in order to prove an analogue of Theorem \ref{Thm3.1} we
assume a condition which is a bit stronger, namely that each $C\left(
z\right) $ is unitary, 
\begin{equation}
C\left( z\right) ^{*}C\left( z\right) =\openone_{n}  \label{EqX.5}
\end{equation}
for almost every $z\in \mathbb{T}$. Define as before 
\begin{equation}
E_{k}=S_{C}^{k}S_{C}^{*\,k}  \label{EqX.6}
\end{equation}
and let 
\begin{equation}
P_{U}=\mathop{\text{\textup{s-lim}}}_{k\rightarrow \infty }E_{k}
\label{EqX.9}
\end{equation}
be the projection onto the subspace corresponding to the unitary part of the
Wold decomposition of $S_{C}$. Again one verifies 
\begin{equation}
\left( E_{k}\xi \right) \left( z\right) =C^{(k)}\left( z\right) \frac{1}{%
N^{k}}\sum_{%
\substack{ w \\ \makebox[0pt]{\hss {$\scriptstyle w^{N^{k}}=z^{N^{k}}
$}\hss }}}C^{(k)}\left( w\right) ^{*}\xi \left( w\right)  \label{EqX.10}
\end{equation}
where 
\begin{equation}
C^{(k)}\left( z\right) =C\left( z\right) C\left( z^{N}\right) \cdots C\left(
z^{N^{k-1}}\right) .  \label{EqX.11}
\end{equation}

The analogue of Theorem \ref{Thm3.1} is now

\begin{theorem}
\label{ThmX.1}Assume that $S_{C}$ is defined by \textup{(\ref{EqX.1})} and
assume that the unitarity condition \textup{(\ref{EqX.5})} is satisfied.
Then the projection $P_{U}$ corresponding to the unitary part of the Wold
decomposition is at most $n$-dimensional. If $\dim P_{U}=m\le n$, the range
of $P_{U}$ can be characterized as follows: there is a projection $P_{0}\in
M_{n}$ of dimension $m$, and a measurable function $\Delta :\mathbb{T}%
\rightarrow M_{n}$, such that 
\begin{equation}
\Delta \left( z\right) ^{*}\Delta \left( z\right) =P_{0}  \label{EqX.12}
\end{equation}
for all $z\in \mathbb{T}$, i.e., $\Delta \left( z\right) $ is a partial
isometry with initial projection $P_{0}$, and a function $\xi \in
L^{2}\left( \mathbb{T};\mathbb{C}^{n}\right) $ is in the range of $P_{U}$ if
and only if there is a vector $v\in P_{0}\mathbb{C}^{n}$ such that 
\begin{equation}
\xi \left( z\right) =\Delta \left( z\right) v  \label{EqX.13}
\end{equation}
for almost all $z$. Furthermore, there is a partial unitary $U_{0}\in M_{n}$
with 
\begin{equation}
U_{0}^{{}}U_{0}^{*}=U_{0}^{*}U_{0}^{{}}=P_{0}^{{}}  \label{EqX.14}
\end{equation}
such that 
\begin{equation}
C\left( z\right) \Delta \left( z^{N}\right) =\Delta \left( z\right) U_{0}.
\label{EqX.15}
\end{equation}
Here $P_{0}$ is the unique maximal projection with the property that there
exist $\Delta \left( \,\cdot \,\right) $ and $U_{0}$ satisfying \textup{(\ref
{EqX.12}), (\ref{EqX.14}),} and \textup{(\ref{EqX.15});} and then $U_{0}$ is
uniquely determined, and $\Delta \left( \,\cdot \,\right) $ is uniquely
determined up to a phase factor.
\end{theorem}

%TCIMACRO{\TeXButton{Begin Proof}{\begin{proof}}}
%BeginExpansion
\begin{proof}%
%EndExpansion
Let $\xi \in P_{U}\left( \mathcal{H}\right) $. For all $k$ we have 
\begin{equation}
E_{k}\xi =\xi ,  \label{EqX.16}
\end{equation}
and it follows from (\ref{EqX.10}) that 
\begin{equation}
C^{(k)}\left( z\right) ^{*}\xi \left( z\right) =\frac{1}{N^{k}}\sum_{%
\substack{ w \\ \makebox[0pt]{\hss {$\scriptstyle w^{N^{k}}=z^{N^{k}}
$}\hss }}}C^{(k)}\left( w\right) ^{*}\xi \left( w\right)  \label{EqX.17}
\end{equation}
for almost all $z\in \mathbb{T}$. But replacing the $z$ to the left with any 
$\eta \in \mathbb{T}$ with $\eta ^{N^{k}}=z^{N^{k}}$, we see that any of the
vectors $C^{(k)}\left( \eta \right) ^{*}\xi \left( \eta \right) $ is a
convex combination of all the vectors of this form with equal weight, and it
follows that 
\begin{equation}
C^{(k)}\left( z\right) ^{*}\xi \left( z\right) =C^{\left( k\right) }\left(
w\right) ^{*}\xi \left( w\right)  \label{EqX.18}
\end{equation}
whenever $z^{N^{k}}=w^{N^{k}}$. At this point we use the unitarity of $%
C^{(k)}\left( w\right) ^{*}$ to deduce 
\begin{equation}
\left\| \xi \left( z\right) \right\| =\left\| \xi \left( w\right) \right\|
\label{EqX.19}
\end{equation}
whenever $z^{N^{k}}=w^{N^{k}}$, and letting $k\rightarrow \infty $ and using
Luzin's theorem (for any $\varepsilon >0$ there is a closed subset $F\subset %
\mathbb{T}$ such that $\mu \left( \mathbb{T}-F\right) <\varepsilon $ and $%
z\rightarrow \left\| \xi \left( z\right) \right\| $ is continuous on $F$,
where $\mu $ is Haar measure) we deduce that $\left\| \xi \left( z\right)
\right\| $ is equal to a constant for almost all $z$. Now, if $\xi ,\eta \in
P_{U}\left( \mathcal{H}\right) $, then all linear combinations of $\xi $ and 
$\eta $ are in $P_{U}\left( \mathcal{H}\right) $, and it follows from the
polarization identity 
\begin{equation}
\ip{\xi \left( z\right) }{\eta \left( z\right) }=\frac{1}{4}%
\sum_{k=0}^{3}i^{k}\left\| \xi \left( z\right) +i^{k}\eta \left( z\right)
\right\|  \label{EqX.20}
\end{equation}
that $z\rightarrow \ip{\xi \left( z\right) }{\eta \left( z\right) }$ is
equal to a constant almost everywhere. It follows that $P_{U}\left( \mathcal{%
H}\right) $ can at most be $n$-dimensional, and modifying the
representatives $\xi \left( z\right) $ on a set of measure zero, we may
assume that 
\begin{equation}
z\rightarrow \ip{\xi \left( z\right) }{\eta \left( z\right) }  \label{EqX.21}
\end{equation}
is a constant for any two $\xi ,\eta \in P_{U}\left( \mathcal{H}\right) $.
But then, if $P_{0}$ is the projection onto the set of $\xi \left( 1\right) $%
, $\xi \in P_{U}\left( \mathcal{H}\right) $, we may for each $v\in P_{0}%
\mathbb{C}^{n}$ find a $\xi $ with $\xi \left( 1\right) =v$, and define 
\begin{equation}
\Delta \left( z\right) v=\Delta \left( z\right) \xi \left( 1\right) =\xi
\left( z\right) .  \label{EqX.22}
\end{equation}
Because of (\ref{EqX.21}), each $\Delta \left( z\right) $ is a partial
isometry with initial projection $P_{0}$, and the statements around (\ref
{EqX.12}) and (\ref{EqX.13}) in the theorem are proved. Furthermore, we have
defined a unitary operator $V:P_{0}\left( \mathbb{C}^{n}\right) \rightarrow
P_{U}\left( \mathcal{H}\right) $ by 
\begin{equation*}
\left( Vv\right) \left( z\right) =\Delta \left( z\right) v.
\end{equation*}
But as $P_{U}S_{C}P_{U}$ is unitary on $P_{U}\mathcal{H}$, we have that 
\begin{equation*}
U_{0}=V^{*}P_{U}S_{C}P_{U}V=V^{*}S_{C}V
\end{equation*}
is unitary on $P_{0}\mathbb{C}^{n}$. But as 
\begin{equation*}
VU_{0}v=S_{C}Vv
\end{equation*}
for $v\in P_{0}\mathbb{C}^{n}$ we have 
\begin{equation*}
\Delta \left( z\right) U_{0}v=C\left( z\right) \Delta \left( z^{N}\right) v
\end{equation*}
and (\ref{EqX.15}) follows. Conversely, it is easy to check from (\ref
{EqX.15}) that $\xi \left( z\right) =\Delta \left( z\right) v$ is in the
range of each $E_{k}$. This ends the proof of Theorem \ref{ThmX.1}.%
%TCIMACRO{\TeXButton{End Proof}{\end{proof}}}
%BeginExpansion
\end{proof}%
%EndExpansion

\section{\label{WoldT}%
%TCIMACRO{\TeXButton{frogkern}{\protect\frogkern } }
%BeginExpansion
\protect\frogkern %
%EndExpansion
The Wold decomposition of operators $T$ on $L^{2}\left( \mathbb{T}\right) $
of the form $\left( T%
%TCIMACRO{\TeXButton{xi}{\protect\xi } }
%BeginExpansion
\protect\xi %
%EndExpansion
\right) \left( z\right) =\frac{1}{\protect\sqrt{N}}\sum%
\limits_{k=0}^{N-1}m_{k}\left( z\right) 
%TCIMACRO{\TeXButton{xi}{\protect\xi } }
%BeginExpansion
\protect\xi %
%EndExpansion
\left( 
%TCIMACRO{\TeXButton{rho}{\protect\rho } }
%BeginExpansion
\protect\rho %
%EndExpansion
^{k}z^{N}\right) $}

Here $\rho =\rho _{N}=e^{\frac{2\pi i}{N}}$, and the coefficient functions $%
m_{k}\in L^{\infty }\left( \mathbb{T}\right) $ satisfy the unitarity
condition (\ref{Eq1.11}), i.e., 
\begin{equation}
C\left( z\right) :=\frac{1}{\sqrt{N}} 
%TCIMACRO{
%\TeXButton{pmatrix}{\begin{pmatrix}
% m_{0}(z) & m_{1}(z)
% & \dots  & m_{N-1}(z) \\ 
% m_{0}(\rho z) & m_{1}(\rho z)
% & \dots  & m_{N-1}(\rho z) \\ 
% \vdots  & \vdots  & \ddots  & \vdots  \\ 
% m_{0}(\rho ^{N-1}z) & m_{1}(\rho
% ^{N-1}z) & \dots  & m_{N-1}(\rho ^{N-1}z)
% \end{pmatrix}} }
%BeginExpansion
\begin{pmatrix}
 m_{0}(z) & m_{1}(z)
 & \dots  & m_{N-1}(z) \\ 
 m_{0}(\rho z) & m_{1}(\rho z)
 & \dots  & m_{N-1}(\rho z) \\ 
 \vdots  & \vdots  & \ddots  & \vdots  \\ 
 m_{0}(\rho ^{N-1}z) & m_{1}(\rho
 ^{N-1}z) & \dots  & m_{N-1}(\rho ^{N-1}z)
 \end{pmatrix}%
%EndExpansion
\label{EqY.1}
\end{equation}
is unitary for almost every $z\in \mathbb{C}$. This condition implies that
the operator $T$ defined by 
\begin{equation}
\left( T\xi \right) \left( z\right) =\frac{1}{\sqrt{N}}\sum_{k=0}^{N-1}m_{k}%
\left( z\right) \xi \left( \rho ^{k}z^{N}\right)  \label{EqY.2}
\end{equation}
is an isometry, since $T$ has the form 
\begin{equation}
T=\frac{1}{\sqrt{N}}\sum_{k=0}^{N-1}S_{k}U^{k},  \label{EqY.3}
\end{equation}
where $S_{0},\dots ,S_{N-1}$ is the representation of $\mathcal{O}_{N}$
given by (\ref{Eq1.16}), and $U$ is the unitary operator on $L^{2}\left( %
\mathbb{T}\right) $ defined by 
\begin{equation*}
\left( U\xi \right) \left( z\right) =\xi \left( \rho z\right) .
\end{equation*}
Let $S_{C}$ be the isometry on $L^{2}\left( \mathbb{T};\mathbb{C}^{N}\right) 
$ defined by (\ref{EqX.1}), 
\begin{equation}
\left( S_{C}\xi \right) \left( z\right) =C\left( z\right) \xi \left(
z^{N}\right) .  \label{EqY.4}
\end{equation}
We now verify that $S_{C}$ is a dilation of $T$. Define an isometric
embedding $V:L^{2}\left( \mathbb{T}\right) \rightarrow L^{2}\left( \mathbb{T}%
;\mathbb{C}^{N}\right) $ by 
\begin{equation}
\left( V\xi \right) \left( z\right) =\frac{1}{\sqrt{N}} 
%TCIMACRO{
%\TeXButton{pmatrix}{\begin{pmatrix}
% \xi (z) \\ \xi (\rho z) \\ \vdots \\ \xi (\rho ^{N-1}z)
% \end{pmatrix}} }
%BeginExpansion
\begin{pmatrix}
 \xi (z) \\ \xi (\rho z) \\ \vdots \\ \xi (\rho ^{N-1}z)
 \end{pmatrix}%
%EndExpansion
.  \label{EqY.5}
\end{equation}
The dilation property is then given by 
\begin{equation}
S_{C}V=VT,  \label{EqY.6}
\end{equation}
which is easily verified from (\ref{EqY.1}), (\ref{EqY.2}), (\ref{EqY.4}),
and (\ref{EqY.5}).

If $S$ is a general isometry, let $\mathcal{H}_{U}\left( S\right) $ denote
the subspace of the Hilbert space corresponding to the unitary part of the
Wold decomposition of $S$, i.e., 
\begin{equation}
\mathcal{H}_{U}\left( S\right) =\bigcap_{k}S^{k}S^{*\,k}\mathcal{H}.
\label{EqY.7}
\end{equation}
It turns out that the unitary subspaces of $T$ and of its dilation $S_{C}$
are the same!

\begin{proposition}
\label{ProY.1}With the assumptions and notation above, 
\begin{equation}
\mathcal{H}_{U}\left( S_{C}\right) =V\left( \mathcal{H}_{U}\left( T\right)
\right) .  \label{EqY.8}
\end{equation}
\end{proposition}

%TCIMACRO{\TeXButton{Begin Proof}{\begin{proof}}}
%BeginExpansion
\begin{proof}%
%EndExpansion
Since $S_{C}$ is a dilation of $T$ in the sense of (\ref{EqY.6}), it is
clear that 
\begin{equation}
V\left( \mathcal{H}_{U}\left( T\right) \right) \subset \mathcal{H}_{U}\left(
S\right),  \label{EqY.9}
\end{equation}
and to prove the reverse inclusion it suffices to show that any $\eta \in 
\mathcal{H}_{U}\left( S\right) $ is in the range of $V$, i.e., that there is
a $\xi \in L^{2}\left( \mathbb{T}\right) $ such that 
\begin{equation}
\eta \left( z\right) = 
%TCIMACRO{
%\TeXButton{pmatrix}{\begin{pmatrix}
% \xi (z) \\ \xi (\rho z) \\ \xi (\rho ^{2}z) \\ \vdots \\ \xi (\rho ^{N-1}z)
% \end{pmatrix}} }
%BeginExpansion
\begin{pmatrix}
 \xi (z) \\ \xi (\rho z) \\ \xi (\rho ^{2}z) \\ \vdots \\ \xi (\rho ^{N-1}z)
 \end{pmatrix}%
%EndExpansion
.  \label{EqY.10}
\end{equation}
But by Theorem \ref{ThmX.1}, $\eta $ has the form 
\begin{equation*}
\eta \left( z\right) =\Delta \left( z\right) v
\end{equation*}
for a suitable $v\in \mathbb{C}^{N}$, and by linearity we may assume that $v$
is an eigenvector of the partial unitary matrix $U_{0}$, i.e., 
\begin{equation*}
U_{0}v=\lambda v,
\end{equation*}
where $\lambda \in \mathbb{T}$. We then obtain from (\ref{EqX.15}) that 
\begin{equation*}
C\left( z\right) \eta \left( z^{N}\right) =\lambda \eta \left( z\right) .
\end{equation*}
If 
\begin{equation*}
\eta \left( z\right) = 
%TCIMACRO{
%\TeXButton{pmatrix}{\begin{pmatrix}
% \xi _{0}(z) \\ \vdots \\ \xi _{N-1}(z)
% \end{pmatrix}} }
%BeginExpansion
\begin{pmatrix}
 \xi _{0}(z) \\ \vdots \\ \xi _{N-1}(z)
 \end{pmatrix}%
%EndExpansion
,
\end{equation*}
we thus obtain from (\ref{EqY.1}) that 
\begin{multline*}
%TCIMACRO{
%\TeXButton{pmatrix}{\begin{pmatrix}
% \xi _{0}(z) \\ \xi _{1}(z) \\ \vdots \\ \xi _{N-1}(z)
% \end{pmatrix}} }
%BeginExpansion
\begin{pmatrix}
 \xi _{0}(z) \\ \xi _{1}(z) \\ \vdots \\ \xi _{N-1}(z)
 \end{pmatrix}%
%EndExpansion
\\
=\frac{\bar{\lambda}}{\sqrt{N}} 
%TCIMACRO{
%\TeXButton{pmatrix}{\begin{pmatrix}
% m_{0}(z) & m_{1}(z)
% & \dots  & m_{N-1}(z) \\ 
% m_{0}(\rho z) & m_{1}(\rho z)
% & \dots  & m_{N-1}(\rho z) \\ 
% \vdots  & \vdots  & \ddots  & \vdots  \\ 
% m_{0}(\rho ^{N-1}z) & m_{1}(\rho
% ^{N-1}z) & \dots  & m_{N-1}(\rho ^{N-1}z)
% \end{pmatrix}} }
%BeginExpansion
\begin{pmatrix}
 m_{0}(z) & m_{1}(z)
 & \dots  & m_{N-1}(z) \\ 
 m_{0}(\rho z) & m_{1}(\rho z)
 & \dots  & m_{N-1}(\rho z) \\ 
 \vdots  & \vdots  & \ddots  & \vdots  \\ 
 m_{0}(\rho ^{N-1}z) & m_{1}(\rho
 ^{N-1}z) & \dots  & m_{N-1}(\rho ^{N-1}z)
 \end{pmatrix}%
%EndExpansion
%TCIMACRO{
%\TeXButton{pmatrix}{\begin{pmatrix}
% \xi _{0}(z^N) \\ \xi _{1}(z^N) \\ \vdots \\ \xi _{N-1}(z^N)
% \end{pmatrix}} }
%BeginExpansion
\begin{pmatrix}
 \xi _{0}(z^N) \\ \xi _{1}(z^N) \\ \vdots \\ \xi _{N-1}(z^N)
 \end{pmatrix}%
%EndExpansion
,
\end{multline*}
and hence, using $\rho ^{kN}=1$, 
\begin{align*}
\xi _{k}\left( z\right) & =\bar{\lambda}\frac{1}{\sqrt{N}}%
\sum_{j=0}^{N-1}m_{j}\left( \rho ^{k}z\right) \xi _{j}\left( z^{N}\right) \\
& =\xi _{0}\left( \rho ^{k}z\right) .
\end{align*}
Thus $\eta $ has the special form (\ref{EqY.10}), and we have established
the reverse inclusion of (\ref{EqY.9}), and thus Proposition \ref{ProY.1}.%
%TCIMACRO{\TeXButton{End Proof}{\end{proof}}}
%BeginExpansion
\end{proof}%
%EndExpansion

Let us summarize the results of this section and the previous one.

\begin{corollary}
\label{CorY.2}The subspace corresponding to the unitary part of the Wold
decomposition of the operator $T$ defined by \textup{(\ref{EqY.2})} has
dimension $n\le N$. Furthermore, there exists a projection $P_{0}\in M_{N}$
of dimension $n$, and measurable functions $d_{0}\left( z\right) ,\dots
,d_{N-1}\left( z\right) $ from $\mathbb{T}$ into $\mathbb{C}$ such that 
\begin{equation}
\sum_{k\in \mathbb{Z}_{N}}\bar{d}_{i}\left( \rho ^{k}z\right) d_{j}\left(
\rho ^{k}z\right) =\left( P_{0}\right) _{ij}  \label{EqY.11}
\end{equation}
for $i,j\in \mathbb{Z}_{N}$ and almost all $z$, such that $\xi \in \mathcal{H%
}_{U}\left( T\right) $ if and only if there are scalars $v_{0},\dots ,v_{N-1}
$ with 
\begin{equation}
\xi \left( z\right) =\sum_{k}d_{k}\left( z\right) v_{k}.  \label{EqY.12}
\end{equation}
If 
\begin{equation}
\Delta \left( z\right) =%
%TCIMACRO{
%\TeXButton{pmatrix}{\begin{pmatrix}
% d_0(z) & d_1(z) & \dots & d_{N-1}(z) \\
% d_0(\rho z) & d_1(\rho z) & \dots & d_{N-1}(\rho z) \\
% \vdots & \vdots & \ddots & \vdots \\
% d_0(\rho ^{N-1}z) & d_1(\rho ^{N-1}z) & \dots & d_{N-1}(\rho ^{N-1}z) 
% \end{pmatrix}} }
%BeginExpansion
\begin{pmatrix}
 d_0(z) & d_1(z) & \dots & d_{N-1}(z) \\
 d_0(\rho z) & d_1(\rho z) & \dots & d_{N-1}(\rho z) \\
 \vdots & \vdots & \ddots & \vdots \\
 d_0(\rho ^{N-1}z) & d_1(\rho ^{N-1}z) & \dots & d_{N-1}(\rho ^{N-1}z) 
 \end{pmatrix}%
%EndExpansion
\label{EqY.13}
\end{equation}
then there exists a partial unitary $U_{0}\in M_{N}$ with 
\begin{equation}
U_{0}^{{}}U_{0}^{*}=U_{0}^{*}U_{0}^{{}}=P_{0}^{{}}  \label{EqY.14}
\end{equation}
such that 
\begin{equation}
C\left( z\right) \Delta \left( z^{N}\right) =\Delta \left( z\right) U_{0}
\label{EqY.15}
\end{equation}
for almost all $z\in \mathbb{T}$. In particular, if $v=\left( 
\begin{smallmatrix}
v_{0}\\
\smash[t]{\vdots}\rule{0pt}{10.5pt} \\
v_{N-1}
\end{smallmatrix}\right) $ is taken to be an eigenvector for $U_{0}$ with
eigenvalue $\lambda \in \mathbb{T}$, then $\xi \left( z\right)
=\sum_{k}d_{k}\left( z\right) v_{k}$ is a Haar vector in the sense 
\begin{equation}
\frac{1}{\sqrt{N}}\sum_{k=0}^{N-1}m_{k}\left( z\right) \xi \left( \rho
^{k}z^{N}\right) =\bar{\lambda}\xi \left( z\right) .  \label{EqY.16}
\end{equation}
\end{corollary}

%TCIMACRO{\TeXButton{Begin Proof}{\begin{proof}}}
%BeginExpansion
\begin{proof}%
%EndExpansion
Let $U_{0}$, $P_{0}$, $\Delta \left( z\right) $ be the objects constructed
in Theorem \ref{ThmX.1} from the matrix $C\left( z\right) $. If we write $%
\Delta \left( z\right) $ as a column vector of row vectors, 
\begin{equation}
\Delta \left( z\right) = 
%TCIMACRO{
%\TeXButton{pmatrix}{\begin{pmatrix}
% \Delta _0(z) \\ \vdots \\ \Delta _{N-1}(z)
% \end{pmatrix}} }
%BeginExpansion
\begin{pmatrix}
 \Delta _0(z) \\ \vdots \\ \Delta _{N-1}(z)
 \end{pmatrix}%
%EndExpansion
,  \label{EqY.17}
\end{equation}
it follows from 
\begin{equation}
C\left( z\right) = 
%TCIMACRO{
%\TeXButton{pmatrix}{\begin{pmatrix}
% C_0(z) \\ C_0(\rho z) \\ \vdots \\ C_0(\rho ^{N-1}z)
% \end{pmatrix}} }
%BeginExpansion
\begin{pmatrix}
 C_0(z) \\ C_0(\rho z) \\ \vdots \\ C_0(\rho ^{N-1}z)
 \end{pmatrix}%
%EndExpansion
,  \label{EqY.18}
\end{equation}
where $C_{0}\left( z\right) $ is the first row of $C\left( z\right) $, and (%
\ref{EqX.15}), that 
\begin{equation}
%TCIMACRO{
%\TeXButton{pmatrix}{\begin{pmatrix}
% \Delta _0(z) \\ \vdots \\ \Delta _{N-1}(z)
% \end{pmatrix}} }
%BeginExpansion
\begin{pmatrix}
 \Delta _0(z) \\ \vdots \\ \Delta _{N-1}(z)
 \end{pmatrix}%
%EndExpansion
= 
%TCIMACRO{
%\TeXButton{pmatrix}{\begin{pmatrix}
% C_0(z) \\ C_0(\rho z) \\ \vdots \\ C_0(\rho ^{N-1}z)
% \end{pmatrix}} }
%BeginExpansion
\begin{pmatrix}
 C_0(z) \\ C_0(\rho z) \\ \vdots \\ C_0(\rho ^{N-1}z)
 \end{pmatrix}%
%EndExpansion
%TCIMACRO{
%\TeXButton{pmatrix}{\begin{pmatrix}
% \Delta _0(z^N) \\ \vdots \\ \Delta _{N-1}(z^N)
% \end{pmatrix}} }
%BeginExpansion
\begin{pmatrix}
 \Delta _0(z^N) \\ \vdots \\ \Delta _{N-1}(z^N)
 \end{pmatrix}%
%EndExpansion
U_{0}^{*},  \label{EqY.19}
\end{equation}
and hence 
\begin{equation}
\Delta _{k}\left( z\right) =\Delta _{0}\left( \rho ^{k}z\right)
\label{EqY.20}
\end{equation}
for $k\in \mathbb{Z}_{N}$. It follows that $\Delta \left( z\right) $ has the
form in (\ref{EqY.13}), and then (\ref{EqY.11}), (\ref{EqY.12}), (\ref
{EqY.14}), and (\ref{EqY.15}) are transcriptions of (\ref{EqX.12}), (\ref
{EqX.13}), (\ref{EqX.14}), and (\ref{EqX.15}), respectively.%
%TCIMACRO{\TeXButton{End Proof}{\end{proof}}}
%BeginExpansion
\end{proof}%
%EndExpansion

\section{\label{CharReal}%
%TCIMACRO{\TeXButton{frogkern}{\protect\frogkern } }
%BeginExpansion
\protect\frogkern %
%EndExpansion
Characterizations of Cuntz algebra representations with $S_{0}$ a shift and
realizations of these representations on a Hardy space}

Recall that an isometry $S$ on a Hilbert space $\mathcal{H}$ is called (by
us) a shift iff $\bigcap_{n=1}^{\infty }S^{n}\mathcal{H}=\left\{ 0\right\} $%
. Then putting $\mathcal{K}=\mathcal{H}\ominus S\mathcal{H}$, letting $\xi
_{0,j}$ be an orthonormal basis for $\mathcal{K}$ and putting $\xi
_{i,j}=S^{i}\xi _{0,j}$, $i\in \mathbb{N}\cup \left\{ 0\right\} $, $\left\{
\xi _{i,j}\right\} $ is an orthonormal basis for $\mathcal{H}$. Let $%
H_{+}^{2}\left( \mathbb{T}\right) $ be the Hardy subspace of $L^{2}\left( %
\mathbb{T}\right) $, i.e., the closed linear span of the orthonormal set of
functions $z\mapsto z^{n}$, $n=1,2,3,\dots $, and define a unitary operator 
\begin{equation}
V:\mathcal{H}\rightarrow H_{+}^{2}\left( \mathbb{T}\right) \otimes \mathcal{K%
}=\mathcal{H}_{+}\left( \mathcal{K}\right)   \label{Eq4.1}
\end{equation}
by 
\begin{equation}
V\xi _{i,j}=z^{i+1}\otimes \xi _{0,j}.  \label{Eq4.2}
\end{equation}
Viewing the elements in $\mathcal{H}_{+}\left( \mathcal{K}\right) $ as
functions from $\mathbb{T}$ into $\mathcal{K}$, 
\begin{equation*}
\xi \in \mathcal{H}_{+}\left( \mathcal{K}\right) \Longleftrightarrow \xi
\left( z\right) =\sum_{n=1}^{\infty }\xi _{n}z^{n},
\end{equation*}
where $\xi _{n}\in \mathcal{K}$ and $\left\| \xi \right\|
^{2}=\sum_{n=1}^{\infty }\left\| \xi _{n}\right\| ^{2}$, $%
S_{0}^{+}=VS_{0}^{{}}V_{{}}^{*}$ is nothing but the operator of
multiplication by $z$: 
\begin{equation}
\left( S_{0}^{+}\xi \right) \left( z\right) =z\xi \left( z\right) .
\label{Eq4.3}
\end{equation}

We will now generalize this description of a shift to a representation $%
s_{i}\rightarrow S_{i}$ of the Cuntz algebra $\mathcal{O}_{N}$ on $\mathcal{H%
}$ such that $S_{0}$ is a shift, when $N=2,3,\dots $.

\begin{lemma}
\label{Lem4.1}There is a $1$\textup{--}$1$ correspondence between
representations $s_{i}\rightarrow S_{i}$ of $\mathcal{O}_{N}$ on $\mathcal{H}
$ such that $S_{0}$ is a shift, and representations of $\mathcal{O}_{\infty }
$ on $\mathcal{H}$ such that the sum of the ranges of the isometries is $%
\openone$. If the representatives of the generators of $\mathcal{O}_{\infty }
$ are denoted by $T_{k}^{\left( \infty ,j\right) }$, where $j=1,\dots ,N-1$, 
$k=1,2,3,\dots $, so that 
\begin{equation}
T_{k_{1}}^{\left( \infty ,j_{1}\right) *}\,T_{k_{2}}^{\left( \infty
,j_{2}\right) }=\delta _{j_{1}j_{2}}\delta _{k_{1}k_{2}}\openone
\label{Eq4.4}
\end{equation}
and 
\begin{equation}
\sum_{j=1\vphantom{k=1}}^{N-1}\sum_{k=1\vphantom{j=1}}^{\infty
}T_{k}^{\left( \infty ,j\right) }\,T_{k}^{\left( \infty ,j\right) *}=\openone%
,  \label{Eq4.5}
\end{equation}
then the $1$\textup{--}$1$ correspondence is given by 
\begin{equation}
T_{k}^{\left( \infty ,j\right) }=S_{0}^{k-1}S_{j}^{{}},\quad j=1,\dots
,N-1,\quad k=1,2,\dots ,  \label{Eq4.6}
\end{equation}
and by 
\begin{equation}
S_{0}^{{}}=\sum_{k=1\vphantom{j=1}}^{\infty }\sum_{j=1\vphantom{k=1}%
}^{N-1}T_{k+1}^{\left( \infty ,j\right) }\,T_{k}^{\left( \infty ,j\right) *}
\label{Eq4.7}
\end{equation}
and 
\begin{equation}
S_{j}^{{}}=T_{1}^{\left( \infty ,j\right) },\quad j=1,\dots ,N-1,
\label{Eq4.8}
\end{equation}
where all infinite sums converge in the strong operator topology.
\end{lemma}

%TCIMACRO{\TeXButton{Begin Proof}{\begin{proof}}}
%BeginExpansion
\begin{proof}%
%EndExpansion
If $S_{0},\dots ,S_{N-1}$ is a representation of $\mathcal{O}_{N}$ on $%
\mathcal{H}$ with $S_{0}$ a shift, define $T_{k}^{\left( \infty ,j\right) }$
by (\ref{Eq4.6}). One uses the Cuntz relations 
\begin{equation}
S_{i}^{*}S_{j}^{{}}=\delta _{ij}^{{}}\openone  \label{Eq4.9}
\end{equation}
to verify (\ref{Eq4.4}). The other Cuntz relation, 
\begin{equation}
\sum_{j}S_{j}^{{}}S_{j}^{*}=\openone,  \label{Eq4.10}
\end{equation}
implies 
\begin{align*}
\sum_{k=1\vphantom{j=1}}^{\infty }\sum_{j=1\vphantom{k=1}}^{N-1}T_{k}^{%
\left( \infty ,j\right) }\,T_{k}^{\left( \infty ,j\right) *}& =\sum_{k=1%
\vphantom{j=1}}^{\infty }\sum_{j=1\vphantom{k=1}%
}^{N-1}S_{0}^{k-1}S_{j}^{{}}S_{j}^{*}S_{0}^{*\,k-1} \\
& =\sum_{k=1}^{\infty }S_{0}^{k-1}\left( \openone-S_{0}^{{}}S_{0}^{*}\right)
S_{0}^{*\,k-1} \\
& =\openone-\lim_{k\rightarrow \infty }S_{0}^{k}S_{0}^{*\,k},
\end{align*}
but the last limit is zero since $S_{0}$ is a shift, and (\ref{Eq4.5})
follows. Furthermore, (\ref{Eq4.8}) is immediate from (\ref{Eq4.6}), while 
\begin{align*}
\sum_{k=1\vphantom{j=1}}^{\infty }\sum_{j=1\vphantom{k=1}}^{N-1}T_{k+1}^{%
\left( \infty ,j\right) }\,T_{k}^{\left( \infty ,j\right) *}& =\sum_{k=1%
\vphantom{j=1}}^{\infty }\sum_{j=1\vphantom{k=1}%
}^{N-1}S_{0}^{k}S_{j}^{{}}S_{j}^{*}S_{0}^{*\,k-1} \\
& =\sum_{k=1}^{\infty }S_{0}^{k}\left( \openone-S_{0}^{{}}S_{0}^{*}\right)
S_{0}^{*\,k-1} \\
& =S_{0}\left( \openone-\lim_{k\rightarrow \infty
}S_{0}^{k}S_{0}^{*\,k}\right) \\
& =S_{0},
\end{align*}
so (\ref{Eq4.7}) is verified.

Conversely, if $T_{k}^{\left( \infty ,j\right) }$ are given satisfying (\ref
{Eq4.4}) and (\ref{Eq4.5}), one verifies that $S_{0}$ and $S_{j}$, $%
j=1,\dots ,N-1,$ satisfy the Cuntz relations (\ref{Eq4.9}) and (\ref{Eq4.10}%
), that (\ref{Eq4.6}) is valid, and that $\lim_{k\rightarrow \infty
}S_{0}^{k}S_{0}^{*\,k}=0$, i.e., $S_{0}$ is a shift.%
%TCIMACRO{\TeXButton{End Proof}{\end{proof}}}
%BeginExpansion
\end{proof}%
%EndExpansion

We will now use this Lemma to construct the announced Hardy-space structure
on $\mathcal{H}$.

\begin{theorem}
\label{Thm4.2}Let $s_{i}\rightarrow S_{i}$ be a representation of $\mathcal{O%
}_{N}$ on a Hilbert space $\mathcal{K}$ such that $S_{0}$ is a shift. Then
there exists a unitary operator 
\begin{equation}
V:\mathcal{H}_{+}\left( \bigoplus_{j=1}^{N-1}\mathcal{K}\right)
\longrightarrow \mathcal{K}  \label{Eq4.11}
\end{equation}
such that if $S_{j}^{+}=V_{{}}^{*}S_{j}^{{}}V$, then 
\begin{equation}
S_{0}^{+}=M_{z}^{{}}=\text{multiplication by }z  \label{Eq4.12}
\end{equation}
and 
\begin{equation}
S_{j}^{+}\psi =z\left( \left( \bigoplus_{i=1}^{j-1}0\right) \oplus V\psi
\oplus \left( \bigoplus_{i=j+1}^{N-1}0\right) \right)   \label{Eq4.13}
\end{equation}
for $j=1,\dots ,N-1$.
\end{theorem}

%TCIMACRO{\TeXButton{Begin Proof}{\begin{proof}}}
%BeginExpansion
\begin{proof}%
%EndExpansion
Define $T_{k}^{\left( \infty ,j\right) }$, $j=1,\dots ,N-1$, $k=1,2,\dots $,
by (\ref{Eq4.6}), and define $V:\mathcal{H}_{+}\left( \bigoplus_{j=1}^{N-1}%
\mathcal{K}\right) \longrightarrow \mathcal{K}$ by 
\begin{equation}
V\left( \sum_{k=1}^{\infty }\left( \bigoplus_{j=1}^{N-1}\psi _{k}^{\left(
j\right) }\right) z_{{}}^{k}\right) =\sum_{k=1\vphantom{j=1}}^{\infty
}\sum_{j=1\vphantom{k=1}}^{N-1}T_{k}^{\left( \infty ,j\right) }\psi
_{k}^{\left( j\right) }.  \label{Eq4.14}
\end{equation}
It follows from (\ref{Eq4.4}) and (\ref{Eq4.5}) in Lemma \ref{Lem4.1} that $%
V $ is indeed unitary and 
\begin{equation}
V_{{}}^{*}\psi =\sum_{k=1}^{\infty }\left(
\bigoplus_{j=1}^{N-1}T_{k}^{\left( \infty ,j\right) *}\psi \right)
z_{{}}^{k}.  \label{Eq4.15}
\end{equation}
Thus, if $\psi \left( z\right) =\sum_{k=1}^{\infty }\left(
\bigoplus_{j=1}^{N-1}\psi _{k}^{\left( j\right) }\right) z_{{}}^{k}$ is in $%
\mathcal{H}_{+}\left( \bigoplus_{j=1}^{N-1}\mathcal{K}\right) $, one checks 
\begin{align*}
\left( S_{0}^{+}\psi \right) \left( z\right) & =\left(
V_{{}}^{*}S_{0}^{{}}V\psi \right) \left( z\right) \\
& =\left( V_{{}}^{*}\left( S_{0}^{{}}\left( \sum_{k=1\vphantom{j=1}}^{\infty
}\sum_{j=1\vphantom{k=1}}^{N-1}T_{k}^{\left( \infty ,j\right) }\psi
_{k}^{\left( j\right) }\right) \right) \right) \left( z\right) \\
& =\left( V_{{}}^{*}\left( \sum_{k=1\vphantom{j=1}}^{\infty }\sum_{j=1%
\vphantom{k=1}}^{N-1}T_{k+1}^{\left( \infty ,j\right) }\psi _{k}^{\left(
j\right) }\right) \right) \left( z\right) \\
& =\sum_{k=2}^{\infty }\left( \bigoplus_{j=1}^{N-1}\psi _{k-1}^{\left(
j\right) }\right) z_{{}}^{k}=z\psi \left( z\right) ,
\end{align*}
where we used 
\begin{equation*}
S_{0}^{{}}T_{k}^{\left( \infty ,j\right) }=T_{k+1}^{\left( \infty ,j\right)
},
\end{equation*}
which follows from (\ref{Eq4.6}). Thus (\ref{Eq4.12}) is valid. For $%
j=1,2,\dots ,N-1$ we use (\ref{Eq4.8}) and (\ref{Eq4.15}) to calculate 
\begin{align*}
\left( S_{j}^{+}\psi \right) \left( z\right) & =\left(
V_{{}}^{*}S_{j}^{{}}V\psi \right) \left( z\right) =\left( V_{{}}^{*}\left(
S_{j}^{{}}\left( \sum_{k=1\vphantom{i=1}}^{\infty }\sum_{i=1\vphantom{k=1}%
}^{N-1}T_{k}^{\left( \infty ,i\right) }\psi _{k}^{\left( i\right) }\right)
\right) \right) \left( z\right) \\
& =z\left( \left( \bigoplus_{i=1}^{j-1}0\right) \oplus V\psi \oplus \left(
\bigoplus_{i=j+1}^{N-1}0\right) \right) ,
\end{align*}
which is (\ref{Eq4.13}).%
%TCIMACRO{\TeXButton{End Proof}{\end{proof}}}
%BeginExpansion
\end{proof}%
%EndExpansion

\section{\label{CharPi}%
%TCIMACRO{\TeXButton{frogkern}{\protect\frogkern } }
%BeginExpansion
\protect\frogkern %
%EndExpansion
Characterization of representations $%
%TCIMACRO{\TeXButton{pi}{\protect\pi } }
%BeginExpansion
\protect\pi %
%EndExpansion
$ of $\mathcal{O}_{N}$ with $%
%TCIMACRO{\TeXButton{pi}{\protect\pi } }
%BeginExpansion
\protect\pi %
%EndExpansion
\left( \mathcal{D}_{N}\right) \subset M_{L^{\infty }\left( \mathbb{T}\right)
}$}

Recall that $\limfunc{UHF}\nolimits_{N}$ is the $C^{*}$-subalgebra of $%
\mathcal{O}_{N}$ which is the closed linear span of $s_{I}^{{}}s_{J}^{*}$
with $\left| I\right| =\left| J\right| $, and $\mathcal{D}_{N}$ is the
canonical maximal abelian subalgebra of $\limfunc{UHF}\nolimits_{N}$, i.e., $%
\mathcal{D}_{N}$ is the closed linear span of operators $s_{I}^{{}}s_{I}^{*}$%
. Thus, if $\limfunc{UHF}\nolimits_{N}\cong \bigotimes_{n=1}^{\infty }M_{N}$%
, then $\mathcal{D}_{N}\cong \bigotimes_{n=1}^{\infty }\mathbb{C}^{N}$.

\begin{theorem}
\label{Thm5.1}Consider a representation $\pi $ of $\mathcal{O}_{N}$ on $%
L^{2}\left( \mathbb{T}\right) $ of the form \textup{(\ref{Eq1.16}),} 
\begin{equation*}
\left( S_{i}\xi \right) \left( z\right) =m_{i}\left( z\right) \xi \left(
z^{N}\right) ,
\end{equation*}
where the functions $m_{i}$ satisfy the appropriate form of the unitarity
condition \textup{(\ref{Eq1.11}).} Let $M_{L^{\infty }\left( \mathbb{T}%
\right) }$ be the image of $L^{\infty }\left( \mathbb{T}\right) $ acting as
multiplication operators on $L^{2}\left( \mathbb{T}\right) $. The following
conditions are equivalent: 
%TCIMACRO{
%\TeXButton{Custom Aligned Equation}{\begin{align}
%& \begin{minipage}[t]{\displayboxwidth}\raggedright $\pi \left( 
%\mathcal{D}_{N}\right) ^{\prime \prime }\subset M_{L^{\infty
%}\left( \mathbb{T}\right) }$;\end{minipage}  \label{Eq5.1} \\
%& \begin{minipage}[t]{\displayboxwidth}\raggedright $\pi \left( 
%\mathcal{D}_{N}\right) ^{\prime \prime }=M_{L^{\infty }\left( 
%\mathbb{T}\right) }$;\end{minipage}  \label{Eq5.2} \\
%& \begin{minipage}[t]{\displayboxwidth}\raggedright $m_{i}\left( 
%z\right) =\sqrt{N}\chi _{A_{i}}\left( z\right) u\left( z\right) 
%$,\end{minipage}
%\label{Eq5.3}
%\end{align}} }
%BeginExpansion
\begin{align}
& \begin{minipage}[t]{\displayboxwidth}\raggedright $\pi \left( 
\mathcal{D}_{N}\right) ^{\prime \prime }\subset M_{L^{\infty
}\left( \mathbb{T}\right) }$;\end{minipage}  \label{Eq5.1} \\
& \begin{minipage}[t]{\displayboxwidth}\raggedright $\pi \left( 
\mathcal{D}_{N}\right) ^{\prime \prime }=M_{L^{\infty }\left( 
\mathbb{T}\right) }$;\end{minipage}  \label{Eq5.2} \\
& \begin{minipage}[t]{\displayboxwidth}\raggedright $m_{i}\left( 
z\right) =\sqrt{N}\chi _{A_{i}}\left( z\right) u\left( z\right) 
$,\end{minipage}
\label{Eq5.3}
\end{align}%
%EndExpansion
where $u$ is a measurable function $\mathbb{T}\rightarrow \mathbb{T}$, and $%
A_{0},\dots ,A_{N-1}$ are $N$ measurable subsets of $\mathbb{T}$ with the
property that if $\rho =e^{\frac{2\pi i}{N}}$, then, for almost all $z\in %
\mathbb{T}$, the $N$ equidistant points $z,\rho z,\rho ^{2}z,\dots ,\rho
^{N-1}z$ lie with one in each of the $N$ sets $A_{0},\dots ,A_{N-1}$ \textup{%
(}i.e., $A_{0},\dots ,A_{N-1}$ form a partition of $\mathbb{T}$ up to null
sets, and, for each $k$, the $N$ sets $A_{k},\rho A_{k},\dots ,\rho
^{N-1}A_{k}$ form a partition of $\mathbb{T}$\textup{).} Any $m_{i}$ of this
form does indeed define a representation of $\mathcal{O}_{N}$.
\end{theorem}

%TCIMACRO{\TeXButton{Begin Proof}{\begin{proof}}}
%BeginExpansion
\begin{proof}%
%EndExpansion
(\ref{Eq5.2})$\Rightarrow $(\ref{Eq5.1}) is trivial, and we prove (\ref
{Eq5.1})$\Rightarrow $(\ref{Eq5.3})$\Rightarrow $(\ref{Eq5.2}).

\emph{Ad} (\ref{Eq5.1})$\Rightarrow $(\ref{Eq5.3}): Assume (\ref{Eq5.1}).
From (\ref{Eq1.16})--(\ref{Eq1.17}), we have 
\begin{equation}
\left( S_{i}^{{}}S_{i}^{*}\xi \right) \left( z\right) =\frac{1}{N}%
m_{i}^{{}}\left( z\right) \sum_{%
\substack{ w \\ \makebox[0pt]{\hss {$\scriptstyle w^{N}=z^{N}
$}\hss }}}\bar{m}_{i}^{{}}\left( w\right) \xi \left( w\right) .
\label{Eq5.4}
\end{equation}
But in order that $S_{i}^{{}}S_{i}^{*}$ be a multiplication operator, we
must thus require that 
\begin{equation*}
m_{i}\left( z\right) \bar{m}_{i}\left( w\right) =0
\end{equation*}
almost everywhere, whenever $w^{N}=z^{N}$ and $w\ne z$. If $A_{i}$ is the
support of $m_{i}$, it follows that the $N$ sets $A_{i}$, $\rho A_{i},\dots
,\rho ^{N-1}A_{i}$ are disjoint (up to null sets). Thus, for given $z\in %
\mathbb{T}$ and $i\in \mathbb{Z}_{N}$, there is at most one $j$ such that $%
\rho ^{\,j}z\in A_{i}$. But, by the unitarity condition (\ref{Eq1.11}), for
given $i$ and $z$, 
\begin{equation*}
\sum_{j}\left| m_{i}\left( \rho ^{\,j}z\right) \right| ^{2}=N,
\end{equation*}
and hence there must exist a $j$ with $\rho ^{\,j}z\in A_{i}$. We have thus
proved that the points $z,\rho z,\rho ^{2}z,\dots ,\rho ^{N-1}z$ lie one
each in the sets $A_{0},A_{1},\dots ,A_{N-1}$, so these sets form a
partition with the stated properties. But then necessarily the functions $%
m_{i}$ must have the form 
\begin{equation*}
m_{i}\left( z\right) =\sqrt{N}\chi _{A_{i}}\left( z\right) u\left( z\right) ,
\end{equation*}
where $u:\mathbb{T}\rightarrow \mathbb{C}$ is a measurable function. But, by
unitarity again, 
\begin{equation*}
1=\frac{1}{N}\sum_{i}\left| m_{i}\left( z\right) \right| ^{2}=\left| u\left(
z\right) \right| ^{2},
\end{equation*}
so that $u$ actually maps into the circle $\mathbb{T}$. Thus (\ref{Eq5.3})
is valid, and we just comment at this point that if $m_{i}\left( z\right) $
is given by (\ref{Eq5.3}), the matrix (\ref{Eq1.11}) is a permutation matrix
for any given $x$, and thus the unitarity condition is fulfilled from the
conditions in (\ref{Eq5.3}) and the last statement of the theorem follows.

\emph{Ad} (\ref{Eq5.3})$\Rightarrow $(\ref{Eq5.2}): Assume (\ref{Eq5.3}).
Using (\ref{Eq1.16})--(\ref{Eq1.17}), and putting $\rho =\rho _{n}=e^{\frac{%
2\pi i}{N^{n}}}$, one has for $I=\left( i_{1},\dots ,i_{n}\right) $: 
\begin{multline}
\left( S_{I}^{{}}S_{I}^{*}\xi \right) \left( z\right)  \label{Eq5.5} \\
\begin{aligned} =\; & m_{i_{1}}\left( z\right) m_{i_{2}}\left( z^{N}\right)
\cdots m_{i_{n}}\left( z^{N^{n-1}}\right) \frac{1}{N^{n}} \\ & \cdot
\sum_{k=0}^{N^{n}-1}\bar{m}_{i_{n}}\left( \rho ^{kN^{n-1}}z^{N^{n-1}}\right)
\bar{m}_{i_{n-1}}\left( \rho ^{kN^{n-2}}z^{N^{n-2}}\right) \cdots
\bar{m}_{i_{1}}\left( \rho ^{k}z\right) \xi \left( \rho ^{k}z\right) \\ =\;
& \chi _{A_{i_{1}}}\left( z\right) \chi _{A_{i_{2}}}\left( z^{N}\right)
\cdots \chi _{A_{i_{n}}}\left( z^{N^{n-1}}\right) \xi \left( z\right) .
\end{aligned}
\end{multline}
But now defining a coding map $\sigma :\mathbb{T}\rightarrow
\prod_{k=1}^{\infty }\mathbb{Z}_{n}$ by $\sigma \left( z\right) =\left(
i_{1},i_{2},i_{3},\dots \right) $, if $z^{N^{n-1}}\in A_{i_{n}}$, it follows
from the properties of $A_{i}$ that $\sigma $ is a measure-preserving map
from $\mathbb{T}$ into $\prod_{k=1}^{\infty }\mathbb{Z}_{n}$ when both
groups are equipped with normalized Haar measure. Replacing $\mathbb{T}$ by $%
\prod_{k=1}^{\infty }\mathbb{Z}_{n}$ by means of this map, the relation (\ref
{Eq5.5}) takes the form 
\begin{equation}
\left( S_{I}^{{}}S_{I}^{*}\xi \right) \left( j_{1},j_{2},\dots \right)
=\delta _{i_{1}j_{1}}\delta _{i_{2}j_{2}}\cdots \delta _{i_{n}j_{n}}\xi
\left( j_{1},j_{2},\dots \right) ,  \label{Eq5.6}
\end{equation}
and it is clear from this relation that the von Neumann algebra generated by 
$\pi \left( \mathcal{D}_{N}\right) $ is exactly $L^{\infty }\left( \bigcross%
_{n=1}^{\infty }\mathbb{Z}_{N}\right) $. This establishes (\ref{Eq5.3})$%
\Rightarrow $(\ref{Eq5.2}), and Theorem \ref{Thm5.1} is proved.%
%TCIMACRO{\TeXButton{End Proof}{\end{proof}}}
%BeginExpansion
\end{proof}%
%EndExpansion

From the last part of the above proof, we also have the following

\begin{corollary}
\label{Cor5.2}If $\pi $ is a representation of $\mathcal{O}_{N}$ on a
Hilbert space $\mathcal{H}$ satisfying the equivalent conditions \textup{(%
\ref{Eq5.1})--(\ref{Eq5.3})} in Theorem \textup{\ref{Thm5.1},} then there is
a unitary operator $U:\mathcal{H}\rightarrow L^{2}\left( \bigcross%
_{n=1}^{\infty }\mathbb{Z}_{n}\right) $ and a measurable function $v:%
\bigcross_{n=1}^{\infty }\mathbb{Z}_{n}\rightarrow \mathbb{T}$ such that 
\begin{equation}
\left( US_{i}U^{*}\xi \right) \left( j_{1},j_{2},\dots \right) =\sqrt{N}%
\delta _{ij_{1}}v\left( j_{1},j_{2},\dots \right) \xi \left(
j_{2},j_{3},\dots \right)   \label{Eq5.7}
\end{equation}
for all $\left( j_{1},j_{2},\dots \right) \in \bigcross_{n=1}^{\infty }%
\mathbb{Z}_{n}$ and all $\xi \in L^{2}\left( \bigcross_{n=1}^{\infty }%
\mathbb{Z}_{n}\right) $.
\end{corollary}

Let us finally remark that the representations described in Corollary \ref
{Cor5.2} are all irreducible, even in restriction to $\limfunc{UHF}%
\nolimits_{N}$, by \cite[Proposition 7.1]{BrJo96a}.

\section{\label{Clas}%
%TCIMACRO{\TeXButton{frogkern}{\protect\frogkern } }
%BeginExpansion
\protect\frogkern %
%EndExpansion
Classification of representations with $%
%TCIMACRO{\TeXButton{pi}{\protect\pi } }
%BeginExpansion
\protect\pi %
%EndExpansion
\left( \mathcal{D}_{N}\right) ^{\prime \prime }\subset M_{L^{\infty }\left( 
\mathbb{T}\right) }$ up to unitary equivalence}

We will now make a further study of the representations of $\mathcal{O}_{N}$
introduced in Section \ref{CharPi}. By Corollary \ref{Cor5.2}, these act on $%
L^{2}\left( \prod_{1}^{\infty }\mathbb{Z}_{N}\right) $ and are labeled by
measurable functions 
\begin{equation}
u:\prod_{1}^{\infty }\mathbb{Z}_{N}\rightarrow \mathbb{T}.  \label{Eq6.1}
\end{equation}
In the case that $u$ depends only on a finite number of the coordinates in $%
\prod_{1}^{\infty }\mathbb{Z}_{N}$, these representations were also
considered in Section 7 in \cite{BrJo96a}. The formulae (\ref{Eq1.12})--(\ref
{Eq1.13}) now take the form 
\begin{align}
\left( S_{i}\xi \right) \left( x_{1},x_{2},\dots \right) & =\sqrt{N}\delta
_{x_{1}i}u\left( x_{1},x_{2},\dots \right) \xi \left( x_{2},x_{3},\dots
\right) ,  \label{Eq6.2} \\
\left( S_{i}^{*}\xi \right) \left( x_{1},x_{2},\dots \right) & =\frac{1}{%
\sqrt{N}}\bar{u}\left( i,x_{1},x_{2},\dots \right) \xi \left(
i,x_{1},x_{2},\dots \right) .  \label{Eq6.3}
\end{align}
We define $\pi ^{u}$ as the representation defined by $u$. By emulating the
proof of irreducibility of $\pi ^{u}$ from \cite{BrJo96a} we can now
establish the following

\begin{proposition}
\label{Pro6.1}Let $T$ be a bounded operator on $L^{2}\left(
\prod_{1}^{\infty }\mathbb{Z}_{N}\right) $, and let $u,u^{\prime
}:\prod_{1}^{\infty }\mathbb{Z}_{N}\rightarrow \mathbb{T}$ be measurable
functions. Then the following conditions are equivalent. 
%TCIMACRO{
%\TeXButton{Custom Aligned Equation}{\begin{align}
%& \begin{minipage}[t]{\displayboxwidth}\raggedright $T\pi ^{u}\left( x\right)
%=\pi ^{u^{\prime }}\!\left( x\right) T$ for all $x\in
%\mathcal{O}_{N}$.\end{minipage}  \label{Eq6.4} \\
%& \begin{minipage}[t]{\displayboxwidth}\raggedright $T =M_{f}$ where $f\in
%L^{\infty }\left( \prod_{1}^{\infty }\mathbb{Z}_{N}\right) $ is a function
%satisfying\end{minipage}  \label{Eq6.5} \\
%& \begin{minipage}[t]{\displayboxwidth} 
%\makebox[\displayboxwidth]{$\displaystyle f\left( x_{1},x_{2},\dots \right) 
%u\left( x_{1},x_{2},\dots \right) =u^{\prime }\left( x_{1},x_{2},\dots \right) 
%f\left( x_{2},x_{3},\dots \right) $}\end{minipage}  \notag \\
%& \begin{minipage}[t]{\displayboxwidth}\raggedright for all 
%$\left( x_{1},x_{2},\dots \right) \in \prod_{1}^{\infty
%}\mathbb{Z}_{N}$.\end{minipage}  \notag
%\end{align}}}
%BeginExpansion
\begin{align}
& \begin{minipage}[t]{\displayboxwidth}\raggedright $T\pi ^{u}\left( x\right)
=\pi ^{u^{\prime }}\!\left( x\right) T$ for all $x\in
\mathcal{O}_{N}$.\end{minipage}  \label{Eq6.4} \\
& \begin{minipage}[t]{\displayboxwidth}\raggedright $T =M_{f}$ where $f\in
L^{\infty }\left( \prod_{1}^{\infty }\mathbb{Z}_{N}\right) $ is a function
satisfying\end{minipage}  \label{Eq6.5} \\
& \begin{minipage}[t]{\displayboxwidth} 
\makebox[\displayboxwidth]{$\displaystyle f\left( x_{1},x_{2},\dots \right) 
u\left( x_{1},x_{2},\dots \right) =u^{\prime }\left( x_{1},x_{2},\dots \right) 
f\left( x_{2},x_{3},\dots \right) $}\end{minipage}  \notag \\
& \begin{minipage}[t]{\displayboxwidth}\raggedright for all 
$\left( x_{1},x_{2},\dots \right) \in \prod_{1}^{\infty
}\mathbb{Z}_{N}$.\end{minipage}  \notag
\end{align}%
%EndExpansion
\end{proposition}

\begin{remark}
\label{Rem6.2}In particular, if $u=u^{\prime }$, (\ref{Eq6.5}) entails $%
f=f\circ \sigma $; and hence $f$ is constant by ergodicity, and this
confirms the irreducibility of $\pi ^{u}$.
\end{remark}

%TCIMACRO{\TeXButton{Begin Proof}{\begin{proof}}}
%BeginExpansion
\begin{proof}%
%EndExpansion
\emph{Ad} (\ref{Eq6.4})$\Rightarrow $(\ref{Eq6.5}): By (\ref{Eq5.6}), we
have 
\begin{equation}
\left( \pi ^{u}\left( s_{I}^{{}}s_{I}^{*}\right) \xi \right) \left(
x_{1},x_{2},\dots \right) =\delta _{i_{1}x_{1}}^{{}}\delta
_{i_{2}x_{2}}^{{}}\cdots \delta _{i_{n}x_{n}}^{{}}\xi \left(
x_{1},x_{2},\dots \right) ,  \label{Eq6.6}
\end{equation}
where the right side is independent of $u$, and hence the intertwining
operator $T$ must commute with $M_{L_{{}}^{\infty }\left( \prod_{1}^{\infty }%
\mathbb{Z}_{N}\right) }$, and as the latter algebra is maximal abelian there
must be an $f\in L_{{}}^{\infty }\left( \prod_{1}^{\infty }\mathbb{Z}%
_{N}\right) $ such that 
\begin{equation*}
T=M_{f}.
\end{equation*}
But then 
\begin{align*}
\left( T\pi ^{u}\left( s_{i}\right) \xi \right) \left( x_{1},x_{2},\dots
\right) & =f\left( x_{1},x_{2},\dots \right) \sqrt{N}\delta _{x_{1}i}u\left(
x_{1},x_{2},\dots \right) \xi \left( x_{2},x_{3},\dots \right) \\
%TCIMACRO{\TeXButton{and}{\intertext{and}} }
%BeginExpansion
\intertext{and}%
%EndExpansion
\left( \pi ^{u^{\prime }}\left( s_{i}\right) T\xi \right) \left(
x_{1},x_{2},\dots \right) & =\sqrt{N}\delta _{x_{1}i}u^{\prime }\left(
x_{1},x_{2},\dots \right) f\left( x_{2},x_{3},\dots \right) \xi \left(
x_{2},x_{3},\dots \right) .
\end{align*}
The intertwining (\ref{Eq6.4}) for $x=s_{i}$ implies that 
\begin{equation*}
f\left( x_{1},x_{2},\dots \right) u\left( x_{1},x_{2},\dots \right)
=u^{\prime }\left( x_{1},x_{2},\dots \right) f\left( x_{2},x_{3},\dots
\right) .
\end{equation*}
This ends the proof of (\ref{Eq6.4})$\Rightarrow $(\ref{Eq6.5}). For the
converse implication, note that the relation in (\ref{Eq6.5}) and the
computation above imply 
\begin{equation*}
M_{f}\pi ^{u}\left( s_{i}\right) =\pi ^{u^{\prime }}\!\left( s_{i}\right)
M_{f}.
\end{equation*}
But as $M_{f}$ is normal, it follows from Fuglede's theorem (or a direct
computation) that 
\begin{equation*}
M_{f}^{{}}\pi _{{}}^{u}\left( s_{i}^{*}\right) =\pi _{{}}^{u^{\prime
}}\left( s_{i}^{*}\right) M_{f}^{{}},
\end{equation*}
and hence (\ref{Eq6.4}) is valid.%
%TCIMACRO{\TeXButton{End Proof}{\end{proof}}}
%BeginExpansion
\end{proof}%
%EndExpansion

If we view $u,u^{\prime }$ as functions on $\mathbb{T}$, the result in
Proposition \ref{Pro6.1} can be stated in terms of the cohomology theory of
Section \ref{Coho}:

\begin{corollary}
\label{Cor6.2}Let $u,u^{\prime }:\mathbb{T}\rightarrow \mathbb{T}$ be
measurable functions and $\pi ^{u}$, $\pi ^{u^{\prime }}\!$ be the
associated irreducible representations of $\mathcal{O}_{N}$. Then the
following conditions are equivalent. 
%TCIMACRO{
%\TeXButton{Custom Aligned Equation}{\begin{align}
%& \begin{minipage}[t]{\displayboxwidth}$\pi ^{u}$ and $\pi ^{u^{\prime }}\!$
%are unitarily equivalent.\end{minipage}  \label{Eq6.7} \\
%& \begin{minipage}[t]{\displayboxwidth}\raggedright The cocycles $u,u^{\prime
%}$ cobound, i.e., there exists a measurable function $\Delta
%:\mathbb{T}\rightarrow \mathbb{T}$ such that\end{minipage}  \label{Eq6.8} \\
%& \begin{minipage}[t]{\displayboxwidth} 
%\makebox[\displayboxwidth]{$\displaystyle \Delta \left( z\right) 
%u\left( z\right) =u^{\prime }\left( z\right) \Delta \left( z^{N}\right) 
%$}\end{minipage}  \notag \\
%& \begin{minipage}[t]{\displayboxwidth}\raggedright for almost all $z\in
%\mathbb{T}$.\end{minipage}  \notag
%\end{align}}}
%BeginExpansion
\begin{align}
& \begin{minipage}[t]{\displayboxwidth}$\pi ^{u}$ and $\pi ^{u^{\prime }}\!$
are unitarily equivalent.\end{minipage}  \label{Eq6.7} \\
& \begin{minipage}[t]{\displayboxwidth}\raggedright The cocycles $u,u^{\prime
}$ cobound, i.e., there exists a measurable function $\Delta
:\mathbb{T}\rightarrow \mathbb{T}$ such that\end{minipage}  \label{Eq6.8} \\
& \begin{minipage}[t]{\displayboxwidth} 
\makebox[\displayboxwidth]{$\displaystyle \Delta \left( z\right) 
u\left( z\right) =u^{\prime }\left( z\right) \Delta \left( z^{N}\right) 
$}\end{minipage}  \notag \\
& \begin{minipage}[t]{\displayboxwidth}\raggedright for almost all $z\in
\mathbb{T}$.\end{minipage}  \notag
\end{align}%
%EndExpansion
\end{corollary}

%TCIMACRO{\TeXButton{Begin Proof}{\begin{proof}}}
%BeginExpansion
\begin{proof}%
%EndExpansion
By Proposition \ref{Pro6.1}, $\pi ^{u}$ and $\pi ^{u^{\prime }}\!$ are
unitarily equivalent if and only if there is a nonzero function $f:\mathbb{T}%
\rightarrow \mathbb{C}$ with 
\begin{equation*}
f\left( z\right) u\left( z\right) =u^{\prime }\left( z\right) f\left(
z^{N}\right) .
\end{equation*}
But as $\left| u\left( z\right) \right| =\left| u^{\prime }\left( z\right)
\right| =1$ we obtain 
\begin{equation*}
\left| f\left( z\right) \right| =\left| f\left( z^{N}\right) \right| ,
\end{equation*}
and by ergodicity of $z\mapsto z^{N}$, $z\mapsto \left| f\left( z\right)
\right| $ is equal to a constant almost everywhere. Let $\Delta \left(
z\right) =f\left( z\right) /\left| f\left( z\right) \right| $. Then $\Delta $
satisfies (\ref{Eq6.8}). Conversely, if $\Delta $ satisfies (\ref{Eq6.8}),
then $M_{\Delta }$ is an intertwiner between the two representations.%
%TCIMACRO{\TeXButton{End Proof}{\end{proof}}}
%BeginExpansion
\end{proof}%
%EndExpansion

In conclusion, the unitary equivalence classes of the representations $\pi
^{u}$ are 
%TCIMACRO{\TeXButton{labeled}{\hbox{labeled }} }
%BeginExpansion
\hbox{labeled }%
%EndExpansion
by the cohomology classes of the cocycles $u$, which are discussed in 
%TCIMACRO{\TeXButton{Section}{\hbox{Section }} }
%BeginExpansion
\hbox{Section }%
%EndExpansion
\ref{Coho}.

\section{\label{Comp}%
%TCIMACRO{\TeXButton{frogkern}{\protect\frogkern } }
%BeginExpansion
\protect\frogkern %
%EndExpansion
Computation of the Hardy-space realizations for the examples coming from
wavelets}

In Section \ref{CharReal} we defined a certain Hardy-space realization of
representations of $\mathcal{O}_{N}$ in the slightly special case that $%
S_{0} $ is a shift. The construction depended on some seemingly arbitrary
choices. However, in the case that $N=2$ and the representation of $\mathcal{%
O}_{2}$ comes from a wavelet in $L^{2}\left( \mathbb{R}\right) $ as
described in \cite{Jor95}, it turns out that these choices are quite
canonical. We will describe this in the present section.

Let us first give a short rundown of the multiresolution analysis of
wavelets from \cite{Dau92}, \cite{MePa93}. The starting point is a function $%
\varphi \in L^{2}\left( \mathbb{R}\right) $, called the \emph{scaling
function} or \emph{father function} with the properties (\ref{Eq7.1}), (\ref
{Eq7.3}), (\ref{Eq7.4a}), and (\ref{Eq7.4b}) below. 
\begin{equation}
\begin{minipage}[t]{\displayboxwidth}\raggedright The set $\left\{ \varphi
\left( \,\cdot\, -k\right) \mid k\in \mathbb{Z}\right\} $ is an orthonormal
set in $L^{2}\left( \mathbb{R}\right) $.\end{minipage}  \label{Eq7.1}
\end{equation}
If we define $\mathcal{V}_{0}$ as the closed linear span of the functions $%
\varphi \left( \,\cdot \,-k\right) $, it is then clear that $\mathcal{V}_{0}$
is invariant under $\mathbb{Z}$-translation. Define \emph{scaling} on $%
L^{2}\left( \mathbb{R}\right) $ as the unitary operator $U$: 
\begin{equation}
\left( U\xi \right) \left( x\right) =2^{-\frac{1}{2}}\xi \left( x/ 2\right) .
\label{Eq7.2}
\end{equation}
Then we assume 
\begin{equation}
U\varphi \in \mathcal{V}_{0}.  \label{Eq7.3}
\end{equation}
We define the multiresolution associated to $\varphi $ as the sequence of
subspaces $\mathcal{V}_{n}=U^{n}\mathcal{V}_{0}$, and this sequence is
decreasing by (\ref{Eq7.3}). The final assumptions on $\varphi $ are 
\begin{subequations}
\label{Eq7.4}
\begin{align}
\bigwedge_{n}\mathcal{V}_{n}& =\left\{ 0\right\} ,  \label{Eq7.4a} \\
\bigvee_{n}\mathcal{V}_{n}& =L^{2}\left( \mathbb{R}\right) .  \label{Eq7.4b}
\end{align}
\end{subequations}
From $\varphi $ one now constructs a \emph{wavelet} or \emph{mother function}
$\psi $ as follows: first use (\ref{Eq7.3}) and expand $U\varphi $ in the
orthonormal basis $\varphi \left( \,\cdot \,-k\right) $: 
\begin{equation}
U\varphi =\sum_{k}a_{k}\varphi \left( \,\cdot \,-k\right) .  \label{Eq7.5}
\end{equation}
Equivalently, using the Fourier transform 
\begin{equation}
\hat{\varphi}\left( t\right) =\frac{1}{\sqrt{2\pi }}\int_{-\infty }^{\infty
}dx\,e^{-ixt}\varphi \left( x\right) ,  \label{Eq7.6}
\end{equation}
the relation (\ref{Eq7.5}) takes the form 
\begin{equation}
\sqrt{2}\hat{\varphi}\left( 2t\right) =m_{0}\left( t\right) \hat{\varphi}%
\left( t\right) ,  \label{Eq7.7}
\end{equation}
where 
\begin{equation}
m_{0}\left( t\right) =\sum_{k}a_{k}e^{-ikt}.  \label{Eq7.8}
\end{equation}
Thus $m_{0}$ is a function of $z=e^{-it}$, and as such $m_{0}\in L^{2}\left( %
\mathbb{T}\right) $. The orthonormality of $\left\{ \varphi \left( \,\cdot
\,-k\right) \right\} $ entails 
\begin{equation}
\left| m_{0}\left( z\right) \right| ^{2}+\left| m_{0}\left( -z\right)
\right| ^{2}=2.  \label{Eq7.9}
\end{equation}

If $W_{0}=\mathcal{V}_{0}^{\perp }\cap \mathcal{V}_{-1}^{{}}$, then for a $%
\xi \in L^{2}\left( \mathbb{R}\right) $ it can be shown that $\xi \in W_{0}$
if and only if $\xi $ has the form 
\begin{equation}
\hat{\xi}\left( 2t\right) =z\bar{m}_{0}\left( -z\right) f\left( z^{2}\right) 
\hat{\varphi}\left( t\right)  \label{Eq7.10}
\end{equation}
for some function $f$, where $z=e^{-it}$. Now, \emph{define} $\psi $ as the
particular function obtained from $\varphi $ in this manner with $f=\frac{1}{%
\sqrt{2}}$, i.e., 
\begin{equation}
\sqrt{2}\hat{\psi}\left( 2t\right) =z\bar{m}_{0}\left( -z\right) \hat{\varphi%
}\left( t\right) =m_{1}\left( z\right) \hat{\varphi}\left( t\right) .
\label{Eq7.11}
\end{equation}
Then the functions $\psi _{n,k}$ defined by 
\begin{equation}
\psi _{n,k}\left( x\right) =2^{-\frac{n}{2}}\psi \left( 2^{-n}x-k\right)
\label{Eq7.12}
\end{equation}
form an orthonormal basis for $L^{2}\left( \mathbb{R}\right) $. In fact, it
follows from the reasoning in \cite{Dau92} that this does not depend on the
specific choice of $f$ above, and any choice of $f$ such that $\left|
f\left( z\right) \right| =\frac{1}{\sqrt{2}}$ almost everywhere will do. For
the specific choice of $f$ we have the explicit expression 
\begin{equation}
\psi \left( x\right) =\sqrt{2}\sum_{k}\left( -1\right) ^{k}\bar{a}%
_{1-k}\varphi \left( 2x-k\right)  \label{Eq7.13}
\end{equation}
as an orthogonal decomposition. For us it is more important to note that the
functions $m_{0},m_{1}$ satisfy the unitarity condition, i.e., the matrix 
\begin{equation}
2^{-\frac{1}{2}} 
%TCIMACRO{
%\TeXButton{pmatrix}{\begin{pmatrix}
%                    m_0(z) & m_0(-z) \\
%                    m_1(z) & m_1(-z)
%                    \end{pmatrix}} }
%BeginExpansion
\begin{pmatrix}
                    m_0(z) & m_0(-z) \\
                    m_1(z) & m_1(-z)
                    \end{pmatrix}%
%EndExpansion
\label{Eq7.14}
\end{equation}
is unitary for all $z\in \mathbb{T}$. This is indeed the case for any $m_{1}$
of the form 
\begin{equation}
m_{1}\left( z\right) =z\bar{m}_{0}\left( -z\right) f\left( z^{2}\right)
\label{Eq7.15}
\end{equation}
where $\left| f\left( z\right) \right| =1$ for all $z$, and, conversely, for 
$m_{0}$ given with (\ref{Eq7.9}), unitarity of (\ref{Eq7.14}) implies (\ref
{Eq7.15}).

Conversely, if $m_{0}$ satisfies (\ref{Eq7.9}) and $m_{1}\left( z\right) =z%
\bar{m}_{0}\left( -z\right) $, iteration of (\ref{Eq7.7}) and (\ref{Eq7.11})
give formal product expansions of $\hat{\varphi}$ and $\hat{\psi}$.
Moreover, it can be shown \cite[Theorem 6.3.6]{Dau92}, that if $m_{0}$ is a
trigonometric polynomial that satisfies $\left| m_{0}\left( t\right) \right|
^{2}+\left| m_{0}\left( t+\pi \right) \right| ^{2}=2$ and $m_{0}\left(
0\right) =\sqrt{2}$, and there exists no nontrivial finite subset $F\subset %
\mathbb{T}$ with $F^{2}\subset F$ such that $m_{0}|_{-F}=0$, then $\varphi
,\psi $ defined by 
\begin{align*}
\hat{\varphi}\left( t\right) & =\left( 2\pi \right) ^{-\frac{1}{2}%
}\prod_{k=1}^{\infty }\left( \frac{m_{0}\left( t\cdot 2^{-k}\right) }{\sqrt{2%
}}\right) , \\
\hat{\psi}\left( t\right) & =e^{-\frac{it}{2}}\bar{m}_{0}\left( \frac{t}{2}%
+\pi \right) \hat{\varphi}\left( \frac{t}{2}\right) 
\end{align*}
are compactly supported functions in $L^{2}\left( \mathbb{R}\right) $ which
are the father and mother functions of a wavelet, and in particular 
\begin{align*}
\varphi \left( x\right) & =\sum_{k}a_{k}\varphi \left( 2x-k\right) , \\
\psi \left( x\right) & =\sqrt{2}\sum_{k}\left( -1\right) ^{k}\bar{a}%
_{-k+1}\varphi \left( 2x-k\right) ,
\end{align*}
where $a_{k}$ are the Fourier coefficients of $m_{0}$: 
\begin{equation*}
m_{0}\left( t\right) =\sum_{k}a_{k}e^{-ikt}.
\end{equation*}
Compare this with the different conditions (\ref{Eq1.44})--(\ref{Eq1.47}).
The case when such non-trivial subsets $F$ of $\mathbb{T}$, as specified
above (and called \emph{cycles\/\/}), do exist, is discussed in Remark \ref
{RemAA.2} below.

Let us now compute the Hardy-space realization from Theorem \ref{Thm4.2} of
the representation of $\mathcal{O}_{2}$ on $L^{2}\left( \mathbb{T}\right) $
defined by (\ref{Eq7.14}). We must for the moment assume that $S_{0}$ is a
shift, i.e., by Theorem \ref{Thm3.1} we must assume that there does \emph{not%
} exist a measurable function $\xi :\mathbb{T}\rightarrow \mathbb{T}$ and a $%
\lambda \in \mathbb{T}$ such that $m\left( z\right) \xi \left( z^{2}\right)
=\lambda \xi \left( z\right) $ for almost all $z\in \mathbb{T}$. We will
show later, in Lemma \ref{Lem7.3}, that this condition is automatically
fulfilled in this situation. The unitary $V:\mathcal{H}_{+}\left( \mathcal{K}%
\right) \rightarrow \mathcal{K}$, where $\mathcal{K}=L^{2}\left( \mathbb{T}%
\right) $, given in general by (\ref{Eq4.14}), is now defined by 
\begin{equation}
V\left( \sum_{k=1}^{\infty }\psi _{k}^{{}}z_{{}}^{k}\right)
=\sum_{k=1}^{\infty }T_{k}^{\infty }\psi _{k}^{{}}  \label{Eq7.15bis}
\end{equation}
for $\psi _{k}\in L^{2}\left( \mathbb{T}\right) $ with $\sum_{k}\left\| \psi
_{k}\right\| ^{2}<\infty $. By (\ref{Eq4.6}), 
\begin{equation}
T_{k}^{\infty }=S_{0}^{k-1}S_{1}^{{}},\quad k=1,2,\dots ,  \label{Eq7.16}
\end{equation}
and hence 
\begin{align}
\left( T_{k}^{\infty }\xi \right) \left( z\right) & =m_{0}^{{}}\left(
z\right) m_{0}^{{}}\left( z^{2}\right) \cdots m_{0}^{{}}\left(
z^{2^{k-2}}\right) m_{1}^{{}}\left( z^{2^{k-1}}\right) \xi \left(
z^{2^{k}}\right)   \label{Eq7.17} \\
& =m_{0}^{\left( k-1\right) }\left( z\right) m_{1}^{{}}\left(
z^{2^{k-1}}\right) \xi \left( z^{2^{k}}\right) .  \notag
\end{align}
(As a general reference to the use of Hardy spaces in operator theory, we
give \cite[Chapter V]{SzFo70}.)

We will connect the Hardy-space description with the wavelet formalism by
means of a unitary 
\begin{equation}
\mathcal{F}_{\varphi }:\mathcal{V}_{0}\rightarrow L^{2}\left( \mathbb{T}%
\right) =\mathcal{K},  \label{Eq7.18}
\end{equation}
an isometric operator 
\begin{equation}
M_{\hat{\varphi}}:\mathcal{K}\rightarrow L^{2}\left( 
%TCIMACRO{
%\TeXButton{RR hat}{\smash{\hat{\mathbb{R}}}\vphantom{\mathbb{R}}} }
%BeginExpansion
\smash{\hat{\mathbb{R}}}\vphantom{\mathbb{R}}%
%EndExpansion
\right) ,  \label{Eq7.19}
\end{equation}
and another unitary operator 
\begin{equation}
J:L^{2}\left( 
%TCIMACRO{
%\TeXButton{RR hat}{\smash{\hat{\mathbb{R}}}\vphantom{\mathbb{R}}} }
%BeginExpansion
\smash{\hat{\mathbb{R}}}\vphantom{\mathbb{R}}%
%EndExpansion
\right) \rightarrow \mathcal{K}\otimes L^{2}\left( \mathbb{T}\right) .
\label{Eq7.20}
\end{equation}
Let us define these. $\mathcal{F}_{\varphi }$ is defined by the requirement
that it maps $\varphi \left( \,\cdot \,-k\right) $ into $e^{-ikt}$, and as $%
\left\{ \varphi \left( \,\cdot \,-k\right) \right\} $ and $\left\{
e^{-ikt}\right\} $ are orthonormal bases for $\mathcal{V}_{0}$ and $%
L^{2}\left( \mathbb{T}\right) $ respectively, $\mathcal{F}_{\varphi }$ is
unitary. $M_{\hat{\varphi}}$ is defined by 
\begin{equation}
\left( M_{\hat{\varphi}}\xi \right) \left( t\right) =\hat{\varphi}\left(
t\right) \xi \left( e^{-it}\right) .  \label{Eq7.21}
\end{equation}
We then have 
\begin{align*}
M_{\hat{\varphi}}\mathcal{F}_{\varphi }\left( \varphi \left( \,\cdot
\,-k\right) \right) \left( t\right) & =\hat{\varphi}\left( t\right) e^{-ikt}
\\
& =\mathcal{F}\left( \varphi \left( \,\cdot \,-k\right) \right) ,
\end{align*}
where $\mathcal{F}$ denotes Fourier transform as defined by (\ref{Eq7.6}).
As $\left\{ \varphi \left( \,\cdot \,-k\right) \right\} $ is an orthonormal
basis for $\mathcal{V}_{0}$ this establishes that the diagram 
\begin{equation}
\begin{array}{ccc}
\mathcal{V}_{0} & 
%TCIMACRO{
%\TeXButton{F phi right}{\overset{\mathcal{F}_\varphi }{\longrightarrow }} }
%BeginExpansion
\overset{\mathcal{F}_\varphi }{\longrightarrow }%
%EndExpansion
& \mathcal{K}=L^{2}\left( \mathbb{T}\right) \\ 
%TCIMACRO{\TeXButton{hook down}{\hookdownarrow } }
%BeginExpansion
\hookdownarrow %
%EndExpansion
&  & 
%TCIMACRO{
%\TeXButton{M phi hat down}{\makebox[0pt]{\hss \rule{0.3pt}{6pt}\hss 
%}\setbox\frogdown=\hbox to 0pt{\hss $\displaystyle \downarrow $\hss 
%}\lower\ht\frogdown\box\frogdown\raisebox{-0.5ex
%}{\makebox[0pt][l]{$\scriptstyle \mkern6mu M_{\hat{\varphi }}$\hss }}} }
%BeginExpansion
\makebox[0pt]{\hss \rule{0.3pt}{6pt}\hss 
}\setbox\frogdown=\hbox to 0pt{\hss $\displaystyle \downarrow $\hss 
}\lower\ht\frogdown\box\frogdown\raisebox{-0.5ex
}{\makebox[0pt][l]{$\scriptstyle \mkern6mu M_{\hat{\varphi }}$\hss }}%
%EndExpansion
\vphantom{\rule[-20pt]{0pt}{36pt}} \\ 
L^{2}\left( \mathbb{R}\right) & 
%TCIMACRO{
%\TeXButton{F right}{\smash[t]{\overset{\mathcal{F}}{\longrightarrow }}} }
%BeginExpansion
\smash[t]{\overset{\mathcal{F}}{\longrightarrow }}%
%EndExpansion
& L^{2}\left( 
%TCIMACRO{
%\TeXButton{RR hat}{\smash{\hat{\mathbb{R}}}\vphantom{\mathbb{R}}} }
%BeginExpansion
\smash{\hat{\mathbb{R}}}\vphantom{\mathbb{R}}%
%EndExpansion
\right)
\end{array}
\label{Eq7.22}
\end{equation}
is commutative, and as $\mathcal{F}_{\varphi }$ and $\mathcal{F}$ are
unitaries, it follows that $M_{\hat{\varphi}}$ is an isometry.

If $\psi _{n,k}$ is the orthonormal basis for $L^{2}\left( \mathbb{R}\right) 
$ given by (\ref{Eq7.12}), then the Fourier transforms 
\begin{equation}
\hat{\psi}_{n,k}\left( t\right) =2^{\frac{n}{2}}e^{-i2^{n}\!kt}\hat{\psi}%
\left( 2^{n}t\right)  \label{Eq7.23}
\end{equation}
form an orthonormal basis for $L^{2}\left( 
%TCIMACRO{
%\TeXButton{RR hat}{\smash{\hat{\mathbb{R}}}\vphantom{\mathbb{R}}}}
%BeginExpansion
\smash{\hat{\mathbb{R}}}\vphantom{\mathbb{R}}%
%EndExpansion
\right) $, and we define $J$ by the requirement 
\begin{equation}
\left( J\hat{\psi}_{n,k}\right) \left( e^{-it},z\right) =e^{-ikt}z^{n}.
\label{Eq7.24}
\end{equation}
$J$ maps the orthonormal basis $\hat{\psi}_{n,k}$ for $L^{2}\left( 
%TCIMACRO{
%\TeXButton{RR hat}{\smash{\hat{\mathbb{R}}}\vphantom{\mathbb{R}}}}
%BeginExpansion
\smash{\hat{\mathbb{R}}}\vphantom{\mathbb{R}}%
%EndExpansion
\right) $ into an orthonormal basis for $\mathcal{K}\otimes L^{2}\left( %
\mathbb{T}\right) =L^{2}\left( \mathbb{T}\right) \otimes L^{2}\left( %
\mathbb{T}\right) $. If $w$ is a $2\pi $-periodic function we see from (\ref
{Eq7.24}) that 
\begin{equation}
J\left( w\left( \,\cdot \,\right) \hat{\psi}\left( \,\cdot \,\right) \right)
\left( e^{-it},z\right) =w\left( e^{-it}\right) ,  \label{Eq7.25}
\end{equation}
and in particular 
\begin{equation}
J\left( \hat{\psi}\right) =1.  \label{Eq7.26}
\end{equation}
More generally, from (\ref{Eq7.24}), 
\begin{equation}
J\left( 2^{\frac{n}{2}}w\left( 2^{n}\,\cdot \,\right) \hat{\psi}\left(
2^{n}\,\cdot \,\right) \right) =w\left( e^{-it}\right) z^{n}.  \label{Eq7.27}
\end{equation}
It is also interesting to note that if $U$ is the scaling map given by (\ref
{Eq7.3}) a simple computation shows 
\begin{equation}
J\mathcal{F}U\mathcal{F}^{*}J^{*}=M_{z},  \label{Eq7.28}
\end{equation}
i.e., $U$ transforms into the operator of multiplication by $z$. This is
because 
\begin{equation}
U\psi _{n,k}=\psi _{n+1,k}.  \label{Eq7.29}
\end{equation}
Let us now connect this to the Hardy-space representation. We have 
\begin{equation}
\mathcal{H}_{+}^{{}}\left( \mathcal{K}\right) =\mathcal{K}\otimes
H_{+}^{2}\left( \mathbb{T}\right) ,  \label{Eq7.30}
\end{equation}
where $H_{+}^{2}\left( \mathbb{T}\right) $ consists of all vectors in $%
L^{2}\left( \mathbb{T}\right) $ with a Fourier expansion of the form $%
\sum_{k=1}^{\infty }a_{k}z^{k}$.

\begin{theorem}
\label{Thm7.1}With the preceding notation and assumptions, the operator $%
S_{0}:L^{2}\left( \mathbb{T}\right) \rightarrow L^{2}\left( \mathbb{T}%
\right) $ defined by $\left( S_{0}\xi \right) \left( z\right) =m_{0}\left(
z\right) \xi \left( z^{2}\right) $ is a shift, and the following diagram
commutes: 
\begin{equation}
\begin{array}{ccccc}
\mathcal{V}_{0} & 
%TCIMACRO{
%\TeXButton{F phi rightleft}{\mkern18mu \smash[b]{\raisebox{0.5ex
%}{\makebox[0pt]{\hss $\displaystyle \overset{\mathcal{F}_\varphi 
%}{\longrightarrow }$\hss }}\raisebox{-0.5ex}{\makebox[0pt]{\hss 
%$\displaystyle \underset{\mathcal{F}_{\varphi }^{-1}}{\longleftarrow 
%}$\hss }}}\mkern18mu } }
%BeginExpansion
\mkern18mu \smash[b]{\raisebox{0.5ex
}{\makebox[0pt]{\hss $\displaystyle \overset{\mathcal{F}_\varphi 
}{\longrightarrow }$\hss }}\raisebox{-0.5ex}{\makebox[0pt]{\hss 
$\displaystyle \underset{\mathcal{F}_{\varphi }^{-1}}{\longleftarrow 
}$\hss }}}\mkern18mu %
%EndExpansion
& \mathcal{K}=L^{2}\left( \mathbb{T}\right)  & 
%TCIMACRO{
%\TeXButton{V rightleft}{\mkern18mu \smash[b]{\raisebox{0.5ex
%}{\makebox[0pt]{\hss $\displaystyle \overset{V^*}{\longrightarrow 
%}$\hss }}\raisebox{-0.5ex}{\makebox[0pt]{\hss $\displaystyle 
%\underset{V}{\longleftarrow }$\hss }}}\mkern18mu } }
%BeginExpansion
\mkern18mu \smash[b]{\raisebox{0.5ex
}{\makebox[0pt]{\hss $\displaystyle \overset{V^*}{\longrightarrow 
}$\hss }}\raisebox{-0.5ex}{\makebox[0pt]{\hss $\displaystyle 
\underset{V}{\longleftarrow }$\hss }}}\mkern18mu %
%EndExpansion
& \mathcal{H}_{+}^{{}}\left( \mathcal{K}\right) =\mathcal{K}\otimes
H_{+}^{2}\left( \mathbb{T}\right)  \\ 
%TCIMACRO{\TeXButton{hook down}{\hookdownarrow } }
%BeginExpansion
\hookdownarrow %
%EndExpansion
&  & 
%TCIMACRO{
%\TeXButton{M phi hat down}{\makebox[0pt]{\hss \rule{0.3pt}{6pt}\hss 
%}\setbox\frogdown=\hbox to 0pt{\hss $\displaystyle \downarrow $\hss 
%}\lower\ht\frogdown\box\frogdown\raisebox{-0.5ex
%}{\makebox[0pt][l]{$\scriptstyle \mkern6mu M_{\hat{\varphi }}$\hss }}}}
%BeginExpansion
\makebox[0pt]{\hss \rule{0.3pt}{6pt}\hss 
}\setbox\frogdown=\hbox to 0pt{\hss $\displaystyle \downarrow $\hss 
}\lower\ht\frogdown\box\frogdown\raisebox{-0.5ex
}{\makebox[0pt][l]{$\scriptstyle \mkern6mu M_{\hat{\varphi }}$\hss }}%
%EndExpansion
\vphantom{\rule[-20pt]{0pt}{36pt}} &  & 
%TCIMACRO{\TeXButton{hook down}{\hookdownarrow } }
%BeginExpansion
\hookdownarrow %
%EndExpansion
\\ 
L^{2}\left( \mathbb{R}\right)  & 
%TCIMACRO{
%\TeXButton{F rightleft}{\mkern18mu \smash[t]{\raisebox{0.5ex
%}{\makebox[0pt]{\hss $\displaystyle \overset{\mathcal{F}}{\longrightarrow 
%}$\hss }}\raisebox{-0.5ex}{\makebox[0pt]{\hss $\displaystyle 
%\underset{\mathcal{F}^{-1}}{\longleftarrow }$\hss }}}\mkern18mu } }
%BeginExpansion
\mkern18mu \smash[t]{\raisebox{0.5ex
}{\makebox[0pt]{\hss $\displaystyle \overset{\mathcal{F}}{\longrightarrow 
}$\hss }}\raisebox{-0.5ex}{\makebox[0pt]{\hss $\displaystyle 
\underset{\mathcal{F}^{-1}}{\longleftarrow }$\hss }}}\mkern18mu %
%EndExpansion
& L^{2}\left( 
%TCIMACRO{
%\TeXButton{RR hat}{\smash{\hat{\mathbb{R}}}\vphantom{\mathbb{R}}}}
%BeginExpansion
\smash{\hat{\mathbb{R}}}\vphantom{\mathbb{R}}%
%EndExpansion
\right)  & 
%TCIMACRO{
%\TeXButton{J rightleft}{\mkern18mu \smash[t]{\raisebox{0.5ex
%}{\makebox[0pt]{\hss $\displaystyle \overset{J}{\longrightarrow }$\hss 
%}}\raisebox{-0.5ex}{\makebox[0pt]{\hss $\displaystyle 
%\underset{J^{-1}}{\longleftarrow }$\hss }}}\mkern18mu } }
%BeginExpansion
\mkern18mu \smash[t]{\raisebox{0.5ex
}{\makebox[0pt]{\hss $\displaystyle \overset{J}{\longrightarrow }$\hss 
}}\raisebox{-0.5ex}{\makebox[0pt]{\hss $\displaystyle 
\underset{J^{-1}}{\longleftarrow }$\hss }}}\mkern18mu %
%EndExpansion
& \mathcal{K}\otimes L^{2}\left( \mathbb{T}\right) 
\end{array}
\label{Eq7.31}
\end{equation}
\end{theorem}

%TCIMACRO{\TeXButton{Begin Proof}{\begin{proof}}}
%BeginExpansion
\begin{proof}%
%EndExpansion
We will establish that $S_{0}$ is a shift in Lemma \ref{Lem7.3}, and hence
the map $V$ is well-defined by Theorem \ref{Thm4.2}. We have already
established commutativity of the left triangle in (\ref{Eq7.22}), so the
right triangle remains, i.e., if $\psi _{k}\in L^{2}\left( \mathbb{T}\right)
=\mathcal{K}$ with $\sum_{k=1}^{\infty }\left\| \psi _{k}\right\|
^{2}<\infty $, we must show 
\begin{equation}
JM_{\hat{\varphi}}V\left( \sideset{}{^{\smash{\oplus}}}{\sum}%
\limits_{k=1}^{\infty }\psi _{k}\right) \left( t,z\right)
=\sum_{k=1}^{\infty }\psi _{k}\left( t\right) z^{k},  \label{Eq7.32}
\end{equation}
where $\sideset{}{^{\smash{\oplus}}}{\sum}\limits_{k=1}^{\infty }\psi
_{k}=\sum\limits_{k=1}^{\infty }\psi _{k}z^{k}$. But by (\ref{Eq7.15}) and (%
\ref{Eq7.17}), 
\begin{equation}
V\left( \sideset{}{^{\smash{\oplus}}}{\sum}\limits_{k=1}^{\infty }\psi
_{k}\right) \left( t\right) =\sum_{k=1}^{\infty }m_{0}\left( t\right)
m_{0}\left( 2t\right) \cdots m_{0}\left( 2^{k-2}t\right) m_{1}\left(
2^{k-1}t\right) \psi _{k}\left( 2^{k}t\right) ,  \label{Eq7.33}
\end{equation}
so by (\ref{Eq7.21}) and the iterated versions of (\ref{Eq7.7}) and (\ref
{Eq7.11}), 
\begin{multline}
\left( M_{\hat{\varphi}}V\right) \left( \sideset{}{^{\smash{\oplus}}}{\sum}%
\limits_{k=1}^{\infty }\psi _{k}\right) \left( t\right)  \label{Eq7.34} \\
\begin{aligned}{\relax} & =\sum_{k=1}^{\infty }\hat{\varphi}\left( t\right)
m_{0}\left( t\right) m_{0}\left( 2t\right) \cdots m_{0}\left(
2^{k-2}t\right) m_{1}\left( 2^{k-1}t\right) \psi _{k}\left( 2^{k}t\right) \\
& =\sum_{k=1}^{\infty }2^{\frac{k}{2}}\hat{\psi}\left( 2^{k}t\right) \psi
_{k}\left( 2^{k}t\right) , \end{aligned}
\end{multline}
where $\psi $ is the mother function. But now apply (\ref{Eq7.27}) to deduce
(\ref{Eq7.32}). This proves Theorem \ref{Thm7.1}.%
%TCIMACRO{\TeXButton{End Proof}{\end{proof}}}
%BeginExpansion
\end{proof}%
%EndExpansion

\begin{corollary}
\label{Cor7.2}The operator $S_{0}$ on $L^{2}\left( \mathbb{T}\right) $, 
\begin{equation*}
\left( S_{0}\xi \right) \left( z\right) =m_{0}\left( z\right) \xi \left(
z^{2}\right) 
\end{equation*}
is a compression of the scaling operator $U$ on $L^{2}\left( \mathbb{R}%
\right) $, 
\begin{equation*}
\left( U\xi \right) \left( x\right) =2^{-\frac{1}{2}}\xi \left( x/2\right) ,
\end{equation*}
in the sense that 
\begin{equation}
S_{0}^{{}}=M_{\hat{\varphi}}^{*}\mathcal{F}U\mathcal{F}_{{}}^{-1}M_{\hat{%
\varphi}}^{{}}.  \label{Eq7.35}
\end{equation}
\end{corollary}

%TCIMACRO{\TeXButton{Begin Proof}{\begin{proof}}}
%BeginExpansion
\begin{proof}%
%EndExpansion
This follows from (\ref{Eq4.12}) and (\ref{Eq7.28}). Both operators act as
multiplication by $z$ on the respective spaces 
\begin{equation*}
\mathcal{K}\otimes H_{+}^{2}\left( \mathbb{T}\right) \subset \mathcal{K}%
\otimes L_{{}}^{2}\left( \mathbb{T}\right) . 
%TCIMACRO{
%\TeXButton{qed}{\settowidth{\qedskip}{$\displaystyle 
%\mathcal{K}\otimes H_{+}^{2}\left( \mathbb{T}\right) \subset 
%\mathcal{K}\otimes L_{{}}^{2}\left( \mathbb{T}\right) 
%.$}\addtolength{\qedskip}{-\textwidth}\makebox[0pt
%][l]{\makebox[-0.5\qedskip][r]{\qed}\hss}}}
%BeginExpansion
\settowidth{\qedskip}{$\displaystyle 
\mathcal{K}\otimes H_{+}^{2}\left( \mathbb{T}\right) \subset 
\mathcal{K}\otimes L_{{}}^{2}\left( \mathbb{T}\right) 
.$}\addtolength{\qedskip}{-\textwidth}\makebox[0pt
][l]{\makebox[-0.5\qedskip][r]{\qed}\hss}%
%EndExpansion
\end{equation*}
\renewcommand{\qed}{\relax}%
%TCIMACRO{\TeXButton{End Proof}{\end{proof}}}
%BeginExpansion
\end{proof}%
%EndExpansion

Referring to the diagram (\ref{Eq7.31}) in Theorem \ref{Thm7.1}, we will use
the term $z$\emph{-transform} for the map from any vertex into the lower
right-hand vertex $\mathcal{K}\otimes L^{2}\left( \mathbb{T}\right) $, where 
$z$ is the variable in $\mathbb{T}$. For example, it follows from (\ref
{Eq7.24}) that the $z$-transform of $\psi _{n,k}$ is $e^{-ikt}z^{n}$, and in
particular the $z$-transform of the mother wavelet $\psi $ itself is $1$.
Let us compute the $z$-transform $F\left( e^{-it},z\right) $ of the father
wavelet $\varphi $. Since $\varphi \in \mathcal{V}_{0}$, it follows by using
the three possible paths from (\ref{Eq7.31}) that $F$ has the form 
\begin{align}
\sum_{k=1}^{\infty }z^{n}w_{k}\left( e^{-it}\right) & =V^{*}\mathcal{F}%
_{\varphi }\left( \varphi \right) \left( e^{-it},z\right)  \label{Eq7.36} \\
& =JM_{\hat{\varphi}}\mathcal{F}_{\varphi }\left( \varphi \right) \left(
e^{-it},z\right)  \notag \\
& =J\mathcal{F}\left( \varphi \right) \left( e^{-it},z\right) ,  \notag
\end{align}
and the first of these identities leads to the following expression for $%
F\left( e^{-it},z\right) $: 
\begin{align}
F\left( e^{-it},z\right) & =V^{*}\mathcal{F}_{\varphi }\left( \varphi
\right) \left( e^{-it},z\right)  \label{Eq7.37} \\
& =\left( V^{*}1\right) \left( e^{-it},z\right)  \notag \\
& =\sum_{n=1}^{\infty }\frac{z^{n}}{2^{n}}\sum_{%
\substack{w \\ \makebox[0pt]{\hss
{$\scriptstyle w^{2^n}=e^{-it}$}\hss }}}\bar{m}_{1}^{{}}\left(
w^{2^{n-1}}\right) \bar{m}_{0}^{(n-1)}\left( w\right) .  \notag
\end{align}

Let us compute this expansion for the Haar wavelet 
\begin{equation}
\varphi \left( x\right) =\chi _{\left[ 0,1\right] }\left( x\right) .
\label{Eq7.38}
\end{equation}
Then 
\begin{equation}
U\left( \varphi \right) \left( x\right) =2^{-\frac{1}{2}}\varphi \left(
x/2\right) =2^{-\frac{1}{2}}\varphi \left( x\right) +2^{-\frac{1}{2}}\varphi
\left( x-1\right) ,  \label{Eq7.39}
\end{equation}
so in (\ref{Eq7.5}) $a_{0}=a_{1}=2^{-\frac{1}{2}}$ and all other
coefficients are zero. Thus from (\ref{Eq7.8}) and (\ref{Eq7.11}), 
\begin{align}
m_{0}\left( t\right) & =2^{\frac{1}{2}}e^{-i\frac{t}{2}}\cos \left( t/
2\right) ,  \label{Eq7.40} \\
m_{1}\left( t\right) & =e^{-it}\bar{m}_{0}\left( t+\pi \right)  \notag \\
& =-2^{\frac{1}{2}}e^{-i\frac{t}{2}}\sin \left( t/2\right) .  \notag
\end{align}
Hence, from (\ref{Eq7.37}), 
\begin{equation}
F\left( e^{-it},z\right) =\sum_{n=1}^{\infty }\left( \frac{z}{\sqrt{2}}%
\right) ^{n}=\frac{z}{\sqrt{2}}\left( 1-\frac{z}{\sqrt{2}}\right) ^{-1}.
\label{Eq7.41}
\end{equation}
Note that in this case, from (\ref{Eq7.13}), 
\begin{align}
\psi \left( x\right) & =\varphi \left( 2x\right) -\varphi \left( 2x-1\right)
\label{Eq7.42} \\
& =\chi _{\left[ 0,\frac{1}{2}\right] }\left( x\right) -\chi _{\left[ \frac{1%
}{2},1\right] }\left( x\right) ,  \notag
\end{align}
and the expression (\ref{Eq7.41}) for the $z$-transform corresponds to the
expansion 
\begin{align}
\varphi \left( x\right) & =\sum_{n=1}^{\infty }2^{-\frac{n}{2}}\psi
_{n,0}\left( x\right)  \label{Eq7.43} \\
& =\sum_{n=1}^{\infty }2^{-n}\psi \left( 2^{-n}x\right) ,  \notag
\end{align}
which can be verified by hand.

Finally, let us compute (\ref{Eq7.41}) by combining (\ref{Eq7.42}) with the
obvious relation 
\begin{equation*}
\varphi \left( x\right) =\varphi \left( 2x\right) +\varphi \left(
2x-1\right) .
\end{equation*}
Adding these, we see 
\begin{equation*}
2\varphi \left( 2x\right) =\varphi \left( x\right) +\psi \left( x\right) ,
\end{equation*}
i.e., 
\begin{equation*}
2^{\frac{1}{2}}\varphi =U\varphi +U\psi .
\end{equation*}
Now take the $z$-transform and use (\ref{Eq7.28}) to obtain 
\begin{equation*}
2^{\frac{1}{2}}F=zF+z,
\end{equation*}
which gives (\ref{Eq7.41}). In general it seems more difficult to obtain
simple expressions for $F$ using (\ref{Eq7.5}) and (\ref{Eq7.13}), $U\varphi
=\sum_{k}a_{k}\varphi \left( \,\cdot \,-k\right) $, $U\psi =\sum_{k}\left(
-1\right) ^{k}\bar{a}_{1-k}\varphi \left( \,\cdot \,-k\right) $, since the $%
z $-transform does not have any particularly simple property with respect to
translation by $1$ in $L^{2}\left( \mathbb{R}\right) $, and a derivation of $%
F$ along these lines leads back to (\ref{Eq7.37}).

\begin{lemma}
\label{Lem7.3}The function $m_{0}$ defined from the father function $\varphi 
$ by \textup{(\ref{Eq7.5})} and \textup{(\ref{Eq7.8})} has the property that 
$\left| m_{0}\left( z\right) \right| \ne 1$ for a set $z$ of positive
measure. In particular the operator $S_{0}:L^{2}\left( \mathbb{T}\right)
\rightarrow L^{2}\left( \mathbb{T}\right) $ defined by $\left( S_{0}\xi
\right) \left( z\right) =m_{0}\left( z\right) \xi \left( z^{2}\right) $ is a
shift.
\end{lemma}

%TCIMACRO{\TeXButton{Begin Proof}{\begin{proof}}}
%BeginExpansion
\begin{proof}%
%EndExpansion
The last statement follows from the former by Theorem \ref{Thm3.1}, and
hence we only need to show that we cannot have $\left| m_{0}\left( z\right)
\right| =1$ almost everywhere. If \emph{ad absurdum} this is the case, it
follows from (\ref{Eq7.7}) that 
\begin{equation*}
\sqrt{2}\left| \hat{\varphi}\left( 2t\right) \right| =\left| \hat{\varphi}%
\left( t\right) \right| 
\end{equation*}
for almost all $t\in \mathbb{R}$. This means that $\left| \hat{\varphi}%
\right| $ is an eigenvector with eigenvalue $1$ of the unitary operator $%
\hat{U}$ on $L^{2}\left( \mathbb{R}\right) $ defined by 
\begin{equation*}
\left( \hat{U}\xi \right) \left( t\right) =\sqrt{2}\xi \left( 2t\right) .
\end{equation*}
But this unitary is a multiple of the two-sided shift by the following
reasoning: we have a decomposition $L^{2}\left( \mathbb{R}\right)
=L^{2}\left( \mathbb{R}_{+}\right) \oplus L^{2}\left( \mathbb{R}_{-}\right) $
into $\hat{U}$- and $\hat{U}^{*}$-invariant subspaces, and it suffices to
consider $L^{2}\left( \mathbb{R}_{+}\right) $. One checks that $%
V:L^{2}\left( \mathbb{R}_{+},dt\right) \rightarrow L^{2}\left( \mathbb{R}%
,ds\right) $ defined by $\left( V\eta \right) \left( s\right) =\eta \left(
e^{s}\right) e^{\frac{s}{2}}$ for $\eta \in L^{2}\left( \mathbb{R}%
_{+},dt\right) $ is unitary and 
\begin{equation*}
\left( V\hat{U}V^{*}\xi \right) \left( s\right) =\xi \left( s+\ln 2\right) 
\end{equation*}
for $\xi \in L^{2}\left( \mathbb{R},ds\right) $. If one furthermore defines
a map $W:L^{2}\left( \mathbb{R}\right) \rightarrow L^{2}\left( \mathbb{T}%
\times \left[ 0,2\pi \right\rangle \right) $ by 
\begin{equation*}
\left( W\xi \right) \left( z\right) \left( s\right) =\sum_{k=-\infty
}^{\infty }z^{-k}\xi \left( s+k\ln 2\right) ,
\end{equation*}
one checks that 
\begin{equation*}
WV\hat{U}=M_{z}WV,
\end{equation*}
i.e., $\hat{U}$ is unitarily equivalent with multiplication by $z$ on $%
L^{2}\left( \mathbb{T}\times \left[ 0,\ln 2\right\rangle \right) \otimes %
\mathbb{C}^{2}$, where $z$ is the $\mathbb{T}$ variable. But this operator
has absolutely continuous spectrum, so we cannot have $1$ as a discrete
eigenvalue. Thus $\left| m_{0}\left( z\right) \right| \ne 1$ on a set of $z$
of positive measure.%
%TCIMACRO{\TeXButton{End Proof}{\end{proof}}}
%BeginExpansion
\end{proof}%
%EndExpansion

\section{\label{Scal}Wavelets of scale $N$}

The construction in Section \ref{Comp} can be generalized in various
directions. One generalization which is well known in wavelet theory is to
replace the strict orthogonality requirement (\ref{Eq7.1}) on the translates
of $\varphi $ by a weaker requirement like, say, 
\begin{equation*}
\left\| \sum_{n\in \mathbb{Z}}\xi _{n}\varphi \left( \,\cdot \,-n\right)
\right\| _{L^{2}\left( \mathbb{R}\right) }^{2}\le c\left\| \xi \right\|
_{\ell ^{2}}^{2}
\end{equation*}
for all $\xi =\left( \xi _{n}\right) _{n\in \mathbb{Z}}$ in $\ell ^{2}=\ell
^{2}\left( \mathbb{Z}\right) $. This will be considered in Section \ref{Fath}%
, and has interest when going from $\mathcal{O}_{2}$-representations back to
wavelets. But before that we will consider another generalization which is
interesting for us but seems to have been merely postulated in wavelet
theory without proper proofs \cite{GrMa92}, \cite{Mey87}, \cite{MRF96}: the
replacement of scale $2$ by scale $N$, with $N\in \left\{ 3,4,5,\dots
\right\} $. In this case we still start with a father function $\varphi $
with the properties (\ref{Eq7.1}),(\ref{Eq7.3}), and (\ref{Eq7.4}) replaced
with their natural generalizations 
\begin{equation}
\begin{minipage}[t]{\displayboxwidth}\centering $\left\{ \varphi \left(
\,\cdot\, -k\right) \mid k\in \mathbb{Z}\right\} $ is an orthonormal set in
$L^{2}\left( \mathbb{R}\right) $,\end{minipage}  \label{Eq8.1}
\end{equation}
\begin{align}
\mathcal{V}_{0}& =\overline{\func{span}}\left\{ \varphi \left( \,\cdot
\,-k\right) \right\} ,  \notag \\
\left( U_{N}\xi \right) \left( x\right) & =N^{-\frac{1}{2}}\xi \left(
x/N\right) ,  \notag \\
U_{N}\varphi & \in \mathcal{V}_{0},  \label{Eq8.2}
\end{align}
\begin{subequations}
\vspace{-\abovedisplayskip}\label{Eq8.3}
\begin{align}
\bigwedge_{n}U_{N}^{n}\mathcal{V}_{0}^{}& =\left\{ 0\right\} ,  
\label{Eq8.3a} \\
\bigvee_{n}U_{N}^{n}\mathcal{V}_{0}^{}& =L^{2}\left( \mathbb{R}\right) .  
\label{Eq8.3b}
\end{align}
\end{subequations}
Again, define $\left( a_{n}\right) \in \ell _{2}$ by 
\begin{equation}
U_{N}\varphi =\smash{\sum_{k}}\vphantom{\sum}a_{k}\varphi \left( \,\cdot
\,-k\right) ,  \label{Eq8.4}
\end{equation}
i.e., 
\begin{equation}
\sqrt{N}\hat{\varphi}\left( Nt\right) =m_{0}\left( t\right) \hat{\varphi}%
\left( t\right) ,  \label{Eq8.5}
\end{equation}
where 
\begin{equation}
m_{0}\left( t\right) =\sum_{k}a_{k}e^{-ikt}.  \label{Eq8.6}
\end{equation}
Thus $m_{0}$ may be viewed as a function on $\mathbb{T}$. As in 
\cite[(5.1.20)]{Dau92}, the orthonormality of $\varphi \left( \,\cdot
\,-k\right) $ is now, by Fourier transform, equivalent to the condition 
\begin{equation}
\func{PER}\left( \left| \hat{\varphi}\right| ^{2}\right) \left( t\right)
:=\sum_{k}\left| \hat{\varphi}\left( t+2\pi k\right) \right| ^{2}=\left(
2\pi \right) ^{-1}  \label{Eq8.7}
\end{equation}
for almost all $t$. Also, using (\ref{Eq8.5}), one has 
\begin{multline}
\func{PER}\left( \left| \hat{\varphi}\right| ^{2}\right) \left( t\right) 
\label{Eq8.8} \\
\begin{aligned}{\relax} & =\frac{1}{N}\sum_{k}\left| m_{0}\left( 
\frac{t+2\pi k}{N}\right) \hat{\varphi}\left( \frac{t+2\pi k}{N}\right) 
\right| ^{2} \\ & =\frac{1}{N}\sum_{m=0\vphantom{n\in 
\mathbb{Z}}}^{N-1}\,\sum_{n\in \mathbb{Z}}\left| m_{0}\left( 
\frac{t+2\pi m}{N}+2\pi n\right) \right| ^{2}\left| \hat{\varphi}\left( 
\frac{t+2\pi m}{N}+2\pi n\right) \right| ^{2} \\ 
& =\frac{1}{N}\sum_{\substack{w \\ \makebox[0pt]{\hss
{$\scriptstyle w^{N}=e^{-it}$}\hss }}}\left| m_{0}\left( w\right) 
\right| ^{2}\func{PER}\left( \left| \hat{\varphi}\right| ^{2}\right) 
\left( w\right) , \end{aligned}
\end{multline}
and combining this with (\ref{Eq8.7}) we see that orthonormality of $\left\{
\varphi \left( \,\cdot \,-k\right) \right\} $ implies 
\begin{equation}
\sum_{k\in \mathbb{Z}_{N}}\left| m_{0}\left( t+2\pi k/N\right) \right|
^{2}=N.  \label{Eq8.9}
\end{equation}
By the example on pages 177--178 of \cite{Dau92}, the relation (\ref{Eq8.9})
does not conversely imply that $\left\{ \varphi \left( \,\cdot \,-k\right)
\right\} $ is orthonormal. If we define an operator $R:L^{\infty }\left( %
\mathbb{T}\right) \rightarrow L^{\infty }\left( \mathbb{T}\right) $ by 
\begin{equation}
\left( R\xi \right) \left( z\right) =\frac{1}{N}\sum_{%
\substack{w \\ \makebox[0pt]{\hss
{$\scriptstyle w^{N}=z$}\hss }}}\left| m_{0}\left( w\right) \right| ^{2}\xi
\left( w\right) ,  \label{Eq8.10}
\end{equation}
then we see that (\ref{Eq8.9}) together with the assumption that the
eigensubspace of $R$ corresponding to eigenvalue $1$ is one-dimensional,
implies that (\ref{Eq8.7}) holds, i.e., $\left\{ \varphi \left( \,\cdot
\,-k\right) \right\} $ is orthonormal. Conversely, one can use the
ergodicity of $z\mapsto z^{N}$ to show that (\ref{Eq8.7}) implies the
eigensubspace of $R$ corresponding to eigenvalue $1$ is one-dimensional; see
Section \ref{AppA}.

More generally, if $\xi \in \mathcal{V}_{-1}^{}=U_{N}^{-1}\mathcal{V}_{0}^{}$
then $U_{N}\xi $ has a decomposition 
\begin{equation}
U_{N}\xi =\sum_{k}\xi _{k}\varphi \left( \,\cdot \,-k\right) ,
\label{Eq8.11}
\end{equation}
and defining 
\begin{equation}
m_{\xi }\left( t\right) =\sum_{k}\xi _{k}e^{-ikt}  \label{Eq8.12}
\end{equation}
we have 
\begin{equation}
\sqrt{N}\hat{\xi}\left( Nt\right) =m_{\xi }\left( t\right) \hat{\varphi}%
\left( t\right) .  \label{Eq8.13}
\end{equation}
Define the operator $T$ on $L^{2}\left( \mathbb{R}\right) $ by 
\begin{equation}
\left( T\xi \right) \left( x\right) =\xi \left( x-1\right) .  \label{Eq8.14}
\end{equation}
If $\xi ,\eta \in \mathcal{V}_{-1}$ one then uses (\ref{Eq8.13}) and (\ref
{Eq8.7}) to compute, as in \cite[(5.1.21)--(5.1.24)]{Dau92}, that the
vectors $\xi $ and $T^{k}\eta $ are orthogonal for all $k\in \mathbb{Z}$ if
and only if 
\begin{equation}
\sum_{k\in \mathbb{Z}_{N}}m_{\xi }\left( t+2\pi k/N\right) \bar{m}_{\eta
}\left( t+2\pi k/N\right) =0  \label{Eq8.15}
\end{equation}
for almost all $t\in \mathbb{R}$, and also $\left\{ \xi \left( \,\cdot
\,-k\right) \right\} $ is an orthonormal set if and only if 
\begin{equation}
\sum_{k\in \mathbb{Z}_{N}}\left| m_{\xi }\left( t+2\pi k/N\right) \right|
^{2}=N.  \label{Eq8.16}
\end{equation}
For the latter statement, one uses (\ref{Eq8.13}) and the same computation
as in (\ref{Eq8.8}) to show 
\begin{align*}
\func{PER}\left( \left| \smash{\hat{\xi}}\vphantom{\hat{\varphi}}\right|
^{2}\right) \left( z\right) & =R\left( \func{PER}\left( \left| \hat{\varphi}%
\right| ^{2}\right) \right) \\
& =R\left( \left( 2\pi \right) ^{-1}1\right) \\
& =\frac{1}{2\pi N}\sum_{%
\substack{w \\ \makebox[0pt]{\hss
{$\scriptstyle w^{N}=z$}\hss }}}\left| m_{\xi }\left( w\right) \right| ^{2},
\end{align*}
and hence $\func{PER}\left( \left| \smash{\hat{\xi}}\vphantom{\hat{\varphi}}%
\right| ^{2}\right) =\left( 2\pi \right) ^{-1}$ almost everywhere if and
only if (\ref{Eq8.16}) holds. Here $R$ is defined by $m_{\xi }$ rather than $%
m_{0}$ as in (\ref{Eq8.10}). Thus, if $\xi \in \mathcal{V}_{-1}$, the
vectors $T^{k}\xi $ are mutually orthogonal if and only if the function $%
\vec{m}_{\xi }:\mathbb{T}\rightarrow \mathbb{C}^{N}$ defined by 
\begin{equation}
\vec{m}_{\xi }^{{}}\left( z\right) =\left( m_{\xi }^{{}}\left( z\right)
,m_{\xi }^{{}}\left( \rho _{N}^{{}}z\right) ,\dots ,m_{\xi }^{{}}\left( \rho
_{N}^{N-1}z\right) \right) ,  \label{Eq8.17}
\end{equation}
where $\rho _{N}=e^{\frac{2\pi i}{N}}$, takes values in the sphere of radius 
$N^{\frac{1}{2}}\ $for almost all $z$, and $T^{k}\xi $ and $T^{l}\eta $ are
mutually orthogonal if and only if 
\begin{equation}
\ip{\vec{m}_\xi (z)}{\vec{m}_\eta (z)}=0  \label{Eq8.18}
\end{equation}
for almost all $z$.

Now given $m_{\xi }:\mathbb{T}\rightarrow \mathbb{C}$ with the property (\ref
{Eq8.16}), the corresponding $\xi \in L^{2}\left( \mathbb{R}\right) $ can be
defined from the relations (\ref{Eq8.11})--(\ref{Eq8.13}). In this way one
may construct a set of functions $\psi _{1},\dots ,\psi _{N-1}$ in $%
L^{2}\left( \mathbb{R}\right) $ such that if $m_{i}\left( z\right) =m_{\psi
_{i}}\left( z\right) $ for $i=1,\dots ,N-1$ and $m_{0}\left( z\right) $ is
given by (\ref{Eq8.9}), then the vectors 
\begin{equation*}
N^{-\frac{1}{2}}\vec{m}_{0}\left( z\right) ,N^{-\frac{1}{2}}\vec{m}%
_{1}\left( z\right) ,\dots ,N^{-\frac{1}{2}}\vec{m}_{N-1}\left( z\right)
\end{equation*}
form an orthonormal basis of $\mathbb{C}^{N}$ for each $z\in \mathbb{T}$.
This can for example be done by choosing a fixed measurable map $F$ from the
unit sphere $S^{2N-1}$ in $\mathbb{C}^{N}$ into $N$-dimensional orthogonal
frames in $\mathbb{C}^{N}$, i.e., 
\begin{equation}
F\left( \vec{x}\right) =\left( \vec{F}_{0}\left( \vec{x}\right) ,\dots ,\vec{%
F}_{N-1}\left( \vec{x}\right) \right) ,  \label{Eq8.19}
\end{equation}
where the vectors $\vec{F}_{0}\left( \vec{x}\right) ,\dots ,\vec{F}%
_{N-1}\left( \vec{x}\right) $ form an orthonormal basis and we assume $\vec{F%
}_{0}\left( \vec{x}\right) =\vec{x}$. It is no problem finding such
measurable maps, but they can be chosen continuous if and only if $N=2,4,8$;
see Remark \ref{Rem8.2} below. Having chosen $F$, one cannot now just set 
\begin{equation}
\vec{m}_{i}\left( z\right) =N^{\frac{1}{2}}\vec{F}_{i}\left( N^{-\frac{1}{2}}%
\vec{m}_{0}\left( z\right) \right) ,  \label{Eq8.20}
\end{equation}
as this may break the particular symmetry enjoyed by the vector functions of
the form (\ref{Eq8.17}), that is, 
\begin{equation}
\vec{m}\left( \rho z\right) =V\vec{m}\left( z\right) ,  \label{Eq8.21}
\end{equation}
where $V$ is the permutation matrix 
\begin{equation}
V= 
%TCIMACRO{
%\TeXButton{pmatrix}{\begin{pmatrix}
%               0 & 1 & 0 & \dots  & 0 & 0 \\
%               0 & 0 & 1 & \dots  & 0 & 0 \\
%               \vdots  & \vdots  & \vdots  & \ddots  & \vdots  & \vdots  \\
%               0 & 0 & 0 & \dots  & 1 & 0 \\
%               0 & 0 & 0 & \dots  & 0 & 1 \\
%               1 & 0 & 0 & \dots  & 0 & 0 \\
%               \end{pmatrix}} }
%BeginExpansion
\begin{pmatrix}
               0 & 1 & 0 & \dots  & 0 & 0 \\
               0 & 0 & 1 & \dots  & 0 & 0 \\
               \vdots  & \vdots  & \vdots  & \ddots  & \vdots  & \vdots  \\
               0 & 0 & 0 & \dots  & 1 & 0 \\
               0 & 0 & 0 & \dots  & 0 & 1 \\
               1 & 0 & 0 & \dots  & 0 & 0 \\
               \end{pmatrix}%
%EndExpansion
.  \label{Eq8.22}
\end{equation}
So one defines $\vec{m}_{i}\left( z\right) $ by (\ref{Eq8.20}) just when $%
0\le \func{Arg}z<2\pi /N$, and then extend the definition to all $z$ by
requiring (\ref{Eq8.21}). Since $V$ is unitary it follows that the ensuing
functions $N^{-\frac{1}{2}}\vec{m}_{0}\left( z\right) ,\dots ,N^{-\frac{1}{2}%
}\vec{m}_{N-1}\left( z\right) $ form an orthonormal basis for each $z\in %
\mathbb{T}$.

Thus we have established most of the following probably known extension of 
\cite[Theorem 5.1.1]{Dau92}, but we have not found the result explicitly in
the literature, other than as a postulate \cite{GrMa92}, \cite{Mey87}, \cite
{MRF96}.

\begin{theorem}
\label{Thm8.1}Let $\varphi \in L^{2}\left( \mathbb{R}\right) $ be a function
satisfying \textup{(\ref{Eq8.1}), (\ref{Eq8.2}), }and \textup{(\ref{Eq8.3}).}
Then there exist $N-1$ functions $\psi _{1},\dots ,\psi _{N-1}$ such that $%
\left\{ T^{k}\psi _{m}\right\} $, $k\in \mathbb{Z}$, $m=1,\dots ,N-1$ forms
an orthogonal basis for $\mathcal{V}_{0}^{\perp }\cap \mathcal{V}_{-1}^{{}}$%
, and thus $\left\{ U_{N}^{n}T_{{}}^{k}\psi _{m}^{{}}\right\} $, $n,k\in %
\mathbb{Z}$, $m=1,\dots ,N-1$ forms an orthonormal basis for $L^{2}\left( %
\mathbb{R}\right) $. Furthermore, the sequences $\psi _{1},\dots ,\psi _{N-1}
$ of such mother functions are in one-to-one correspondence with the
sequences $m_{1},\dots ,m_{N-1}$ of functions in $L^{2}\left( \mathbb{T}%
\right) $ with the property that $N^{-\frac{1}{2}}\vec{m}_{0}\left( z\right)
,N^{-\frac{1}{2}}\vec{m}_{1}\left( z\right) ,\dots ,N^{-\frac{1}{2}}\vec{m}%
_{N-1}\left( z\right) $ is an orthonormal set for almost all $z\in \mathbb{T}
$. The correspondence is given by 
\begin{equation}
\sqrt{N}\hat{\psi}_{k}\left( Nt\right) =m_{k}\left( e^{-it}\right) \hat{%
\varphi}\left( t\right)   \label{Eq8.23}
\end{equation}
for $k=1,\dots ,N-1$. Furthermore, if $\psi _{1},\dots ,\psi _{M}$ is any
sequence in $\mathcal{V}_{0}^{\perp }\cap \mathcal{V}_{-1}^{{}}$ such that $%
\left\{ T^{k}\psi _{m}\right\} $ forms an orthonormal set, then $M\le N-1$,
and $\left\{ T^{k}\psi _{m}\right\} $ then is an orthonormal basis if and
only if $M=N-1$.
\end{theorem}

%TCIMACRO{\TeXButton{Begin Proof}{\begin{proof}}}
%BeginExpansion
\begin{proof}%
%EndExpansion
The only thing remaining to prove is that $\left\{ T^{k}\psi _{m}\right\} $, 
$k\in \mathbb{Z}$, $m=1,\dots ,N-1$ really forms a basis for $\mathcal{V}%
_{0}^{\perp }\cap \mathcal{V}_{-1}^{{}}$ when $\psi _{m}$ is constructed as
before. But any $\xi \in \mathcal{V}_{-1}$ has the form 
\begin{equation*}
\sqrt{N}\hat{\xi}\left( Nt\right) =m\left( t\right) \hat{\varphi}\left(
t\right) ,
\end{equation*}
where $m\in L^{2}\left( \mathbb{T}\right) $, by the reasoning leading to (%
\ref{Eq8.13}), and $\xi $ being orthogonal to $\mathcal{V}_{0}$ means 
\begin{equation*}
\sum_{k=1}^{N-1}m\left( \rho ^{k}z\right) \bar{m}_{0}\left( \rho
^{k}z\right) =0
\end{equation*}
by the reasoning leading to (\ref{Eq8.15}). But this means that $\vec{m}%
\left( z\right) $ is orthogonal to $\vec{m}_{0}\left( z\right) $ for almost
all $z\in \mathbb{T}$, and it follows that $\vec{m}\left( z\right) $ is a
linear combination of $\vec{m}_{1}\left( z\right) ,\dots ,\vec{m}%
_{N-1}\left( z\right) $ for almost all $z$: 
\begin{equation*}
\vec{m}\left( z\right) =\sum_{k=1}^{N-1}\mu _{k}\left( z\right) \vec{m}%
_{k}\left( z\right) .
\end{equation*}
The symmetries 
\begin{align*}
\vec{m}\left( \rho z\right) & =V\vec{m}\left( z\right) , \\
\vec{m}_{k}\left( \rho z\right) & =V\vec{m}_{k}\left( z\right)
\end{align*}
imply that 
\begin{equation*}
\mu _{k}\left( \rho z\right) =\mu _{k}\left( z\right) ,
\end{equation*}
and hence 
\begin{equation*}
\mu _{k}\left( z\right) =\lambda _{k}\left( z^{N}\right)
\end{equation*}
for a suitable function $\lambda _{k}$ on $\mathbb{T}$. Thus 
\begin{align*}
\sqrt{N}\hat{\xi}\left( Nt\right) & =\sum_{k=1}^{N-1}\lambda _{k}\left(
Nt\right) m_{k}\left( t\right) \hat{\varphi}\left( t\right) \\
& =\sum_{k=1}^{N}\lambda _{k}\left( Nt\right) \sqrt{N}\hat{\psi}_{k}\left(
Nt\right) ,
\end{align*}
i.e., 
\begin{equation*}
\hat{\xi}\left( t\right) =\sum_{k=1}^{N}\lambda _{k}\left( t\right) \hat{\psi%
}_{k}\left( t\right) .
\end{equation*}
But putting 
\begin{equation*}
\lambda _{k}^{{}}\left( t\right) =\sum_{m}a_{m}^{(k)}e_{{}}^{-ikt},
\end{equation*}
this means 
\begin{equation*}
\xi =\sum_{k,m}a_{m}^{(k)}\psi _{k}^{{}}\left( \,\cdot \,-m\right) ,
\end{equation*}
so $\left\{ T^{m}\psi _{k}\right\} $ is a basis.%
%TCIMACRO{\TeXButton{End Proof}{\end{proof}}}
%BeginExpansion
\end{proof}%
%EndExpansion

\begin{remark}
\label{Rem8.2}In Daubechies's approach to wavelets, the selection function $F
$ in (\ref{Eq8.19}) has a particularly simple form like 
\begin{equation}
F_{0}=%
%TCIMACRO{
%\TeXButton{pmatrix}{\begin{pmatrix}
%              1 & 0 \\
%              0 & 1
%              \end{pmatrix}}}
%BeginExpansion
\begin{pmatrix}
              1 & 0 \\
              0 & 1
              \end{pmatrix}%
%EndExpansion
,\quad F_{1}=%
%TCIMACRO{
%\TeXButton{pmatrix}{\begin{pmatrix}
%              0 & 1 \\
%              -1 & 0
%              \end{pmatrix}}}
%BeginExpansion
\begin{pmatrix}
              0 & 1 \\
              -1 & 0
              \end{pmatrix}%
%EndExpansion
.  \label{Eq8.25}
\end{equation}
When $N>2$ one cannot always choose the selection function this simple, and
by a celebrated theorem of Adams, it cannot even be chosen continuous except
in the cases $N=2,4,8$. If $S^{n-1}$ is the unit sphere in $\mathbb{R}^{n}$,
Adams's theorem \cite{Ada62} says that the highest number of pointwise
linearly independent vector fields that may be defined on $S^{n-1}$ is $\rho
\left( n\right) -1$, where the function $\rho \left( n\right) $ is defined
as follows: let $b$ be the multiplicity of $2$ in the prime decomposition of 
$n$; write $b=c+4d$ where $c\in \left\{ 0,1,2,3\right\} $ and $d\in \left\{
0,1,2,\dots \right\} $; and put $\rho \left( n\right) =2^{c}+8d$. One checks
that $\rho \left( 2N\right) -1\ge N$ if and only if $N=2,4,8$ (where $\rho
\left( 2N\right) -1=3,7,8$). In the cases $N=2,4,8$ where the maps $F$ can
be chosen continuous, they may actually be chosen very simple: if $N=2,$ let 
$F_{1}$ be multiplication by $i$ on $\mathbb{R}^{2}=\mathbb{C}$; if $N=4$,
let $F_{1}$, $F_{2}$, $F_{3}$ be multiplication by $i$, $j$, $k$
respectively on the quaternions; and when $N=8$ use the same method with the
Cayley numbers, and then view the resulting (real) orthogonal $N\times N$
matrices as unitary matrices.
\end{remark}

At this stage it is no surprise that the appropriate analogue of Theorem \ref
{Thm7.1} is also true in this more general setting. If $\mathcal{K}%
=L^{2}\left( \mathbb{T}\right) $, define the unitary 
\begin{equation}
V:\mathcal{H}_{+}\left( \bigoplus_{j=1}^{N-1}\mathcal{K}\right) \rightarrow 
\mathcal{K}  \label{Eq8.26}
\end{equation}
as in (\ref{Eq4.14}). Put $\bigoplus_{j=1}^{N-1}\mathcal{K}=\mathbb{C}%
^{N-1}\otimes \mathcal{K}$ in such a way that the $j$'th component of $%
\bigoplus_{j=1}^{N-1}\mathcal{K}$ identifies with $e_{j}\otimes \mathcal{K}$%
, where $\left\{ e_{1},\dots ,e_{N-1}\right\} $ is the standard basis for $%
\mathbb{C}^{N-1}$. In this case we define the unitary map 
\begin{equation}
J:L^{2}\left( 
%TCIMACRO{
%\TeXButton{RR hat}{\smash{\hat{\mathbb{R}}}\vphantom{\mathbb{R}}} }
%BeginExpansion
\smash{\hat{\mathbb{R}}}\vphantom{\mathbb{R}}%
%EndExpansion
\right) \rightarrow \mathbb{C}^{N-1}\otimes \mathcal{K}\otimes L^{2}\left( %
\mathbb{T}\right)  \label{Eq8.27}
\end{equation}
by the requirement 
\begin{equation}
J\left( U_{N}^{n}T_{}^{k}\psi _{m}^{}\right) \left( e^{-it},z\right)
=e_{m}\otimes e^{-ikt}\otimes z^{n}.  \label{Eq8.28}
\end{equation}
The following result is now proved exactly as Theorem \ref{Thm7.1}:

\begin{corollary}
\label{Cor8.3}With the preceding notation and assumptions, the operator $%
S_{0}:L^{2}\left( \mathbb{T}\right) \rightarrow L^{2}\left( \mathbb{T}%
\right) $ defined by $\left( S_{0}\xi \right) \left( z\right) =m_{0}\left(
z\right) \xi \left( z^{N}\right) $ is a shift, and the following diagram
commutes: 
\begin{equation}
\begin{array}{ccccc}
\mathcal{V}_{0} & 
%TCIMACRO{
%\TeXButton{F phi rightleft}{\mkern18mu \smash[b]{\raisebox{0.5ex
%}{\makebox[0pt]{\hss $\displaystyle \overset{\mathcal{F}_\varphi 
%}{\longrightarrow }$\hss }}\raisebox{-0.5ex}{\makebox[0pt]{\hss 
%$\displaystyle \underset{\mathcal{F}_{\varphi }^{-1}}{\longleftarrow 
%}$\hss }}}\mkern18mu } }
%BeginExpansion
\mkern18mu \smash[b]{\raisebox{0.5ex
}{\makebox[0pt]{\hss $\displaystyle \overset{\mathcal{F}_\varphi 
}{\longrightarrow }$\hss }}\raisebox{-0.5ex}{\makebox[0pt]{\hss 
$\displaystyle \underset{\mathcal{F}_{\varphi }^{-1}}{\longleftarrow 
}$\hss }}}\mkern18mu %
%EndExpansion
& \mathcal{K}=L^{2}\left( \mathbb{T}\right)  & 
%TCIMACRO{
%\TeXButton{V rightleft}{\mkern18mu \smash[b]{\raisebox{0.5ex
%}{\makebox[0pt]{\hss $\displaystyle \overset{V^*}{\longrightarrow 
%}$\hss }}\raisebox{-0.5ex}{\makebox[0pt]{\hss $\displaystyle 
%\underset{V}{\longleftarrow }$\hss }}}\mkern18mu } }
%BeginExpansion
\mkern18mu \smash[b]{\raisebox{0.5ex
}{\makebox[0pt]{\hss $\displaystyle \overset{V^*}{\longrightarrow 
}$\hss }}\raisebox{-0.5ex}{\makebox[0pt]{\hss $\displaystyle 
\underset{V}{\longleftarrow }$\hss }}}\mkern18mu %
%EndExpansion
& \mathbb{C}^{N-1}\otimes \mathcal{K}\otimes H_{+}^{2}\left( \mathbb{T}%
\right)  \\ 
%TCIMACRO{\TeXButton{hook down}{\hookdownarrow } }
%BeginExpansion
\hookdownarrow %
%EndExpansion
&  & 
%TCIMACRO{
%\TeXButton{M phi hat down}{\makebox[0pt]{\hss \rule{0.3pt}{6pt}\hss 
%}\setbox\frogdown=\hbox to 0pt{\hss $\displaystyle \downarrow $\hss 
%}\lower\ht\frogdown\box\frogdown\raisebox{-0.5ex}{\makebox[0pt
%][l]{$\scriptstyle \mkern6mu M_{\hat{\varphi }}$\hss }}}}
%BeginExpansion
\makebox[0pt]{\hss \rule{0.3pt}{6pt}\hss 
}\setbox\frogdown=\hbox to 0pt{\hss $\displaystyle \downarrow $\hss 
}\lower\ht\frogdown\box\frogdown\raisebox{-0.5ex}{\makebox[0pt
][l]{$\scriptstyle \mkern6mu M_{\hat{\varphi }}$\hss }}%
%EndExpansion
\vphantom{\rule[-20pt]{0pt}{36pt}} &  & 
%TCIMACRO{\TeXButton{hook down}{\hookdownarrow } }
%BeginExpansion
\hookdownarrow %
%EndExpansion
\\ 
L^{2}\left( \mathbb{R}\right)  & 
%TCIMACRO{
%\TeXButton{F rightleft}{\mkern18mu \smash[t]{\raisebox{0.5ex
%}{\makebox[0pt]{\hss $\displaystyle \overset{\mathcal{F}}{\longrightarrow 
%}$\hss }}\raisebox{-0.5ex}{\makebox[0pt]{\hss $\displaystyle 
%\underset{\mathcal{F}^{-1}}{\longleftarrow }$\hss }}}\mkern18mu } }
%BeginExpansion
\mkern18mu \smash[t]{\raisebox{0.5ex
}{\makebox[0pt]{\hss $\displaystyle \overset{\mathcal{F}}{\longrightarrow 
}$\hss }}\raisebox{-0.5ex}{\makebox[0pt]{\hss $\displaystyle 
\underset{\mathcal{F}^{-1}}{\longleftarrow }$\hss }}}\mkern18mu %
%EndExpansion
& L^{2}\left( 
%TCIMACRO{
%\TeXButton{RR hat}{\smash{\hat{\mathbb{R}}}\vphantom{\mathbb{R}}}}
%BeginExpansion
\smash{\hat{\mathbb{R}}}\vphantom{\mathbb{R}}%
%EndExpansion
\right)  & 
%TCIMACRO{
%\TeXButton{J rightleft}{\mkern18mu \smash[t]{\raisebox{0.5ex
%}{\makebox[0pt]{\hss $\displaystyle \overset{J}{\longrightarrow 
%}$\hss }}\raisebox{-0.5ex}{\makebox[0pt]{\hss $\displaystyle 
%\underset{J^{*}}{\longleftarrow }$\hss }}}\mkern18mu } }
%BeginExpansion
\mkern18mu \smash[t]{\raisebox{0.5ex
}{\makebox[0pt]{\hss $\displaystyle \overset{J}{\longrightarrow 
}$\hss }}\raisebox{-0.5ex}{\makebox[0pt]{\hss $\displaystyle 
\underset{J^{*}}{\longleftarrow }$\hss }}}\mkern18mu %
%EndExpansion
& \mathbb{C}^{N-1}\otimes \mathcal{K}\otimes L^{2}\left( \mathbb{T}\right) 
\end{array}
\label{Eq8.29}
\end{equation}
In particular, the operator $S_{0}$ is a compression of the scaling operator 
$\left( U_{N}\xi \right) \left( x\right) =N^{-\frac{1}{2}}\xi \left(
x/N\right) $ in the sense 
\begin{equation}
S_{0}^{{}}=M_{\hat{\varphi}}^{*}\mathcal{F}U_{N}^{{}}\mathcal{F}_{{}}^{-1}M_{%
\hat{\varphi}}^{{}},  \label{Eq8.30}
\end{equation}
both operators being represented by multiplication by $z$ by the $z$%
-transform.
\end{corollary}

We end this section by proving the formula (\ref{Eq1.37}) in the
Introduction.

\begin{corollary}
\label{Cor8.4}Adopt the preceding notation and assumptions. For any $\xi \in
L^{2}\left( \mathbb{R}\right) $, let 
\begin{equation*}
\xi =\sum_{i=1\vphantom{j,k\in \mathbb{Z}}}^{N-1}\sum_{j,k\in \mathbb{Z}%
}a_{jk}^{(i)}\left( \xi \right) U_{N}^{j}T_{{}}^{k}\psi _{i}^{{}}
\end{equation*}
be the orthonormal decomposition of $\xi $ relative to the wavelet basis in
Theorem \textup{\ref{Thm8.1},} and if $\xi \in \mathcal{V}_{0}$ \textup{(}%
i.e., $a_{jk}^{(i)}\left( \xi \right) =0$ for $j\le 0$\textup{),} let $f\in
L^{2}\left( \mathbb{T}\right) =L^{2}\left( \mathbb{R}\diagup 2\pi \mathbb{Z}%
\right) $ be the unique function such that 
\begin{equation*}
\hat{\xi}\left( t\right) =f\left( t\right) \hat{\varphi}\left( t\right) .
\end{equation*}
Then 
\begin{equation*}
a_{jk}^{(i)}\left( \xi \right) =\left( S_{i}^{*}S_{0}^{*\,j-1}f\right)
^{\sim }\left( k\right) 
\end{equation*}
for $i=1,\dots ,N-1$, $j=1,2,\dots $, $k\in \mathbb{Z}$, where $\left(
\,\cdot \,\right) ^{\sim }$ refers to the Fourier transform \textup{(\ref
{Eq1.38})} on $L^{2}\left( \mathbb{T}\right) $.
\end{corollary}

%TCIMACRO{\TeXButton{Begin Proof}{\begin{proof}}}
%BeginExpansion
\begin{proof}%
%EndExpansion
We have 
\begin{align*}
a_{jk}^{(i)}\left( \xi \right) & =N^{-\frac{j}{2}}\int_{\mathbb{R}}\bar{\psi}%
_{i}\left( N^{-j}x-k\right) \xi \left( x\right) \,dx \\
& =N^{\frac{j}{2}}\int_{\hat{\mathbb{R}}}e_{i}^{+ikN^{j}t}\Bar{\Hat{\psi}}%
_{i}\left( N^{j}t\right) f\left( t\right) \hat{\varphi}\left( t\right) \,dt
\\
& =N^{\frac{j}{2}}\int_{\mathbb{T}}e^{+ikN^{j}\!t}f\left( t\right) \func{PER}%
\left( \Bar{\Hat{\psi}}_{i}\left( N^{j}\,\cdot \,\right) \hat{\varphi}\left(
\,\cdot \,\right) \right) \left( t\right) \,dt.
\end{align*}
Using (\ref{Eq8.23}) and then (\ref{Eq8.5}) iteratively, we have further 
\begin{align*}
a_{jk}^{(i)}\left( \xi \right) & =\int_{\mathbb{T}}e^{+ikN^{j}\!t}f\left(
t\right) \bar{m}_{0}\left( N^{j-1}t\right) \bar{m}_{i}\left( N^{j-2}t\right)
\cdots \bar{m}_{i}\left( Nt\right) \bar{m}_{i}\left( t\right) \func{PER}%
\left( \Bar{\Hat{\varphi}}\hat{\varphi}\right) \left( t\right) \,dt \\
& =\left( 2\pi \right) ^{-1}\int_{\mathbb{T}}e^{+ikN^{j}\!t}f\left( t\right) 
\bar{m}_{0}\left( N^{j-1}t\right) \cdots \bar{m}_{i}\left( t\right) \,dt,
\end{align*}
where the last step used the orthonormality of $\left\{ \varphi \left(
\,\cdot \,-k\right) \right\} $ in the form (\ref{Eq8.7}). Introducing $%
e_{k}\left( t\right) =e^{-ikt}$, or in complex form $e_{k}\left( z\right)
=z^{k}$, we furthermore compute 
\begin{align*}
a_{jk}^{(i)}\left( \xi \right) & =\frac{1}{2\pi }\int_{\mathbb{T}}\overline{%
S_{0}^{j-1}S_{i}^{{}}\left( e_{k}\right) \left( t\right) }f\left( t\right)
\,dt \\
& =\ip{S_0^{j-1}S_i^{}(e_k^{})}{f}_{L^{2}\left( \mathbb{T}\right) } \\
& =\ip{e_k^{}}{S_i^{*}S_0^{*\,j-1}f}_{L^{2}\left( \mathbb{T}\right) } \\
& =\frac{1}{2\pi }\int_{-\pi }^{\pi }e_{{}}^{ikt}\left(
S_{i}^{*}S_{0}^{*\,j-1}f\right) \left( t\right) \,dt,
\end{align*}
which is the desired conclusion.%
%TCIMACRO{\TeXButton{End Proof}{\end{proof}}}
%BeginExpansion
\end{proof}%
%EndExpansion

Generalizations of these results to non-orthogonal translates of the father
wavelet will be given in Section \ref{Fath}.

\section{\label{Scat}Scattering theory for scale $2$ versus scale $N$}

Let $\Phi \in L^{2}\left( \mathbb{R}\right) $ be a scale-$N$ father function
as introduced in (\ref{Eq8.1})--(\ref{Eq8.3}). By iteration of (\ref{Eq8.5})
we have 
\begin{equation}
\hat{\Phi}\left( t\right) =\prod_{k=1}^{n}\left( N^{-\frac{1}{2}}M_{0}\left(
tN^{-k}\right) \right) \hat{\Phi}\left( t/N^{n}\right) ,  \label{Eq9.1}
\end{equation}
where we now denote the function $m_{0}$ by $M_{0}$. It is known (by Remark
3 following Proposition 5.3.2 in \cite{Dau92}) that $\hat{\Phi}\left(
0\right) \ne 0$, and thus $M_{0}\left( 0\right) =N^{\frac{1}{2}}$ by (\ref
{Eq8.5}). Thus, if we assume that $\hat{\Phi}$ is continuous at $0$, it
follows from (\ref{Eq9.1}) and the normalization (\ref{Eq8.7}) that 
\begin{equation}
\hat{\Phi}\left( t\right) =\left( 2\pi \right) ^{-\frac{1}{2}%
}\prod_{k=1}^{\infty }\left( N^{-\frac{1}{2}}M_{0}\left( tN^{-k}\right)
\right) ,  \label{Eq9.2}
\end{equation}
at least up to a phase factor, but we choose the latter to be $1$.
Correspondingly, if $\varphi $ is a scale-$2$ father function, we have 
\begin{equation}
\hat{\varphi}\left( t\right) =\left( 2\pi \right) ^{-\frac{1}{2}%
}\prod_{k=1}^{\infty }\left( 2^{-\frac{1}{2}}m_{0}\left( t2^{-k}\right)
\right)  \label{Eq9.3}
\end{equation}
under slight regularity assumptions, where $m_{0}$ is now defined by (\ref
{Eq7.8}). Throughout this section we will assume that $\Phi $ and $\varphi $
are sufficiently regular that (\ref{Eq9.2}) and (\ref{Eq9.3}) are valid.
(See also the discussion at the end of Section \ref{Intr}.) Denote the
associated isometries defined by (\ref{Eq1.16}) by $T_{0}$ and $S_{0}$,
respectively, i.e., 
\begin{align}
\left( T_{0}\xi \right) \left( z\right) & =M_{0}\left( z\right) \xi \left(
z^{N}\right) ,  \label{Eq9.4} \\
\left( S_{0}\xi \right) \left( z\right) & =m_{0}\left( z\right) \xi \left(
z^{2}\right) .  \label{Eq9.5}
\end{align}

\begin{proposition}
\label{Pro9.1}Adopt the notation and assumptions above. The following three
conditions are equivalent. 
%TCIMACRO{
%\TeXButton{Custom Aligned Equation}{\begin{align}
%& \begin{minipage}[t]{\displayboxwidth}\raggedright $\varphi =\Phi 
%$.\end{minipage}  \label{Eq9.6} \\
%& \begin{minipage}[t]{\displayboxwidth}\raggedright $m_{0}\left( 
%Nt\right) M_{0}\left( t\right) =m_{0}\left( t\right) M_{0}\left(
%2t\right) $ for almost all $t\in \mathbb{R}$.\end{minipage}  \label{Eq9.7} \\
%& \begin{minipage}[t]{\displayboxwidth}\raggedright 
%$S_{0}T_{0}=T_{0}S_{0}$.\end{minipage}  \label{Eq9.8}
%\end{align}}}
%BeginExpansion
\begin{align}
& \begin{minipage}[t]{\displayboxwidth}\raggedright $\varphi =\Phi 
$.\end{minipage}  \label{Eq9.6} \\
& \begin{minipage}[t]{\displayboxwidth}\raggedright $m_{0}\left( 
Nt\right) M_{0}\left( t\right) =m_{0}\left( t\right) M_{0}\left(
2t\right) $ for almost all $t\in \mathbb{R}$.\end{minipage}  \label{Eq9.7} \\
& \begin{minipage}[t]{\displayboxwidth}\raggedright 
$S_{0}T_{0}=T_{0}S_{0}$.\end{minipage}  \label{Eq9.8}
\end{align}%
%EndExpansion
\end{proposition}

%TCIMACRO{\TeXButton{Begin Proof}{\begin{proof}}}
%BeginExpansion
\begin{proof}%
%EndExpansion
(\ref{Eq9.7})$\Leftrightarrow $(\ref{Eq9.8}): If $\xi \in L^{2}\left( %
\mathbb{T}\right) $ then 
\begin{align*}
\left( S_{0}T_{0}\xi \right) \left( z\right) & =m_{0}\left( z\right) \left(
T_{0}\xi \right) \left( z^{2}\right) \\
& =m_{0}\left( z\right) M_{0}\left( z^{2}\right) \xi \left( z^{2N}\right) \\
%TCIMACRO{\TeXButton{and}{\intertext{and}} }
%BeginExpansion
\intertext{and}%
%EndExpansion
\left( T_{0}S_{0}\xi \right) \left( z\right) & =M_{0}\left( z\right) \left(
S_{0}\xi \right) \left( z^{N}\right) \\
& =M_{0}\left( z\right) m_{0}\left( z^{N}\right) \xi \left( z^{2N}\right) ,
\end{align*}
so (\ref{Eq9.7})$\Leftrightarrow $(\ref{Eq9.8}) is immediate.

(\ref{Eq9.6})$\Rightarrow $(\ref{Eq9.7}): If $\Phi =\varphi $, it follows
from (\ref{Eq9.1}) with $n=1$ and $N=2,N$ that 
\begin{align*}
\hat{\varphi}\left( t\right) & =2^{-\frac{1}{2}}m_{0}\left( t/2\right) \hat{%
\varphi}\left( t/2\right) \\
& =2^{-\frac{1}{2}}m_{0}\left( t/2\right) N^{-\frac{1}{2}}M_{0}\left(
t/2N\right) \hat{\varphi}\left( t/2N\right) \\
%TCIMACRO{\TeXButton{and}{\intertext{and}} }
%BeginExpansion
\intertext{and}%
%EndExpansion
\hat{\varphi}\left( t\right) & =N^{-\frac{1}{2}}M_{0}\left( t/N\right) \hat{%
\varphi}\left( t/N\right) \\
& =N^{-\frac{1}{2}}M_{0}\left( t/N\right) 2^{-\frac{1}{2}}m_{0}\left(
t/2N\right) \hat{\varphi}\left( t/2N\right) ,
\end{align*}
and (\ref{Eq9.7}) is immediate.

(\ref{Eq9.7})$\Rightarrow $(\ref{Eq9.6}): Assuming (\ref{Eq9.7}) we have 
\begin{align*}
M_{0}\left( t\right) \hat{\varphi}\left( t\right) & =\left( 2\pi \right) ^{-%
\frac{1}{2}}M_{0}\left( t\right) 2^{-\frac{1}{2}}m_{0}\left( t/2\right)
\prod_{k=2}^{\infty }\left( 2^{-\frac{1}{2}}m_{0}\left( t2^{-k}\right)
\right) \\
& =\left( 2\pi \right) ^{-\frac{1}{2}}2^{-\frac{1}{2}}m_{0}\left(
Nt/2\right) M_{0}\left( t/2\right) \prod_{k=2}^{\infty }\left( 2^{-\frac{1}{2%
}}m_{0}\left( t2^{-k}\right) \right) \\
& =\left( 2\pi \right) ^{-\frac{1}{2}}2^{-\frac{2}{2}}m_{0}\left(
Nt/2\right) m_{0}\left( Nt/2^{2}\right) M_{0}\left( t/2^{2}\right)
\prod_{k=3}^{\infty }\left( 2^{-\frac{1}{2}}m_{0}\left( t2^{-k}\right)
\right) \\
& =\dots \\
& =\left( 2\pi \right) ^{-\frac{1}{2}}\prod_{k=1}^{n}\left( 2^{-\frac{1}{2}%
}m_{0}\left( Nt2^{-k}\right) \right) M_{0}\left( t2^{-n}\right)
\prod_{k=n+1}^{\infty }\left( 2^{-t}m_{0}\left( t2^{-k}\right) \right) \\
& \underset{n\rightarrow \infty }{\longrightarrow }\hat{\varphi}\left(
Nt\right) M_{0}\left( 0\right) =N^{\frac{1}{2}}\hat{\varphi}\left( Nt\right)
\end{align*}
Thus 
\begin{equation*}
N^{\frac{1}{2}}\hat{\varphi}\left( Nt\right) =M_{0}\left( t\right) \hat{%
\varphi}\left( t\right) .
\end{equation*}
But from this one deduces in the same way as (\ref{Eq9.2}) that 
\begin{equation*}
\hat{\varphi}\left( t\right) =\left( 2\pi \right) ^{-\frac{1}{2}%
}\prod_{k=1}^{\infty }\left( N^{-\frac{1}{2}}M_{0}\left( tN^{-k}\right)
\right) ,
\end{equation*}
and hence 
\begin{equation*}
\hat{\varphi}\left( t\right) =\hat{\Phi}\left( t\right) . 
%TCIMACRO{
%\TeXButton{qed}{\settowidth{\qedskip}{$\displaystyle \hat{\varphi}\left( 
%t\right) =\hat{\Phi}\left( t\right) .$}\addtolength{\qedskip
%}{-\textwidth}\makebox[0pt][l]{\makebox[-0.5\qedskip][r]{\qed}\hss}}}
%BeginExpansion
\settowidth{\qedskip}{$\displaystyle \hat{\varphi}\left( 
t\right) =\hat{\Phi}\left( t\right) .$}\addtolength{\qedskip
}{-\textwidth}\makebox[0pt][l]{\makebox[-0.5\qedskip][r]{\qed}\hss}%
%EndExpansion
\end{equation*}
\renewcommand{\qed}{\relax}%
%TCIMACRO{\TeXButton{End Proof}{\end{proof}}}
%BeginExpansion
\end{proof}%
%EndExpansion

Proposition \ref{Pro9.1} gives a characterization of the scale-$2$ father
functions $\varphi $ which are also of scale $N$. Next, assume that $\varphi 
$ is only a scale-$2$ father function satisfying (\ref{Eq7.1})--(\ref{Eq7.4}%
), define $m_{0}$ by (\ref{Eq7.8}), and next define $m_{1}$ by (\ref{Eq7.11}%
), i.e., 
\begin{equation}
m_{1}\left( z\right) =z\bar{m}_{0}\left( -z\right) .  \label{Eq9.9}
\end{equation}
Let $\psi $ be the corresponding mother function defined by (\ref{Eq7.11})
or (\ref{Eq7.13}). Then we have the orthogonal decomposition 
\begin{equation}
N^{-\frac{1}{2}}\varphi \left( x/N\right) =\sum_{k}A_{k}\varphi \left(
x-k\right) +\sum_{k}B_{k}\psi \left( x-k\right) +\xi _{-}\left( x\right) ,
\label{Eq9.10}
\end{equation}
where $\xi _{-}\in \mathcal{V}_{-1}^{\perp }$ and $\sum_{k}\left( \left|
A_{k}\right| ^{2}+\left| B_{k}\right| ^{2}\right) +\left\| \xi _{-}\right\|
_{2}^{2}=1$. Define 
\begin{equation}
A\left( t\right) =\sum_{k}A_{k}e^{-ikt},\quad B\left( t\right)
=\sum_{k}B_{k}e^{-ikt}.  \label{Eq9.11}
\end{equation}

\begin{proposition}
\label{Pro9.2}If $\varphi $ is a scale-$2$ father function we have, with the
notation introduced above, 
\begin{equation}
A\left( 2t\right) m_{0}\left( t\right) +B\left( 2t\right) m_{1}\left(
t\right) =A\left( t\right) m_{0}\left( Nt\right) ,  \label{Eq9.12}
\end{equation}
or, in terms of the representation $S_{0},S_{1}$ of $\mathcal{O}_{2}$
defined by $\varphi $, 
\begin{equation}
S_{0}\left( A\right) \left( z\right) +S_{1}\left( B\right) \left( z\right)
=m_{0}\left( z^{N}\right) A\left( z\right) .  \label{Eq9.13}
\end{equation}
In particular, the functions in the left sum are orthogonal, so 
%TCIMACRO{
%\TeXButton{Custom Aligned Equation}{\begin{align}
%\left\| A\right\| _{L^{2}\left( \mathbb{T}\right) }^{2}+\left\| B\right\|
%_{L^{2}\left( \mathbb{T}\right) }^{2} &=\int_{\mathbb{T}}\left| m_{0}\left(
%z^{N}\right) A\left( z\right) \right| ^{2}\frac{\left| dz\right| }{2\pi }
%\label{Eq9.14} \\
%&=N^{-1}\int_{\mathbb{T}}\left| m_{0}\left( z\right) \right| 
%^{2}\sum_{\substack{ w \\ \makebox[0pt]{\hss {$\scriptstyle w^{N}=z
%$}\hss }}}\left| A\left( w\right) \right| ^{2}\frac{\left| dz\right| 
%}{2\pi }.  \notag
%\end{align}}}
%BeginExpansion
\begin{align}
\left\| A\right\| _{L^{2}\left( \mathbb{T}\right) }^{2}+\left\| B\right\|
_{L^{2}\left( \mathbb{T}\right) }^{2} &=\int_{\mathbb{T}}\left| m_{0}\left(
z^{N}\right) A\left( z\right) \right| ^{2}\frac{\left| dz\right| }{2\pi }
\label{Eq9.14} \\
&=N^{-1}\int_{\mathbb{T}}\left| m_{0}\left( z\right) \right| 
^{2}\sum_{\substack{ w \\ \makebox[0pt]{\hss {$\scriptstyle w^{N}=z
$}\hss }}}\left| A\left( w\right) \right| ^{2}\frac{\left| dz\right| 
}{2\pi }.  \notag
\end{align}%
%EndExpansion
\end{proposition}

%TCIMACRO{\TeXButton{Begin Proof}{\begin{proof}}}
%BeginExpansion
\begin{proof}%
%EndExpansion
By Fourier transform of (\ref{Eq9.10}) we have 
\begin{equation}
N^{\frac{1}{2}}\hat{\varphi}\left( Nt\right) =A\left( t\right) \hat{\varphi}%
\left( t\right) +B\left( t\right) \hat{\psi}\left( t\right) +\hat{\xi}%
_{-}\left( t\right) .  \label{Eq9.15}
\end{equation}
Thus, by (\ref{Eq7.7}) and (\ref{Eq9.15}), 
\begin{align}
\left( 2N\right) ^{\frac{1}{2}}\hat{\varphi}\left( 2Nt\right) & =N^{\frac{1}{%
2}}m_{0}\left( Nt\right) \hat{\varphi}\left( Nt\right)  \label{Eq9.16} \\
& =m_{0}\left( Nt\right) \left( A\left( t\right) \hat{\varphi}\left(
t\right) +B\left( t\right) \hat{\psi}\left( t\right) +\hat{\xi}_{-}\left(
t\right) \right) .  \notag
\end{align}
On the other hand, by (\ref{Eq9.15}), (\ref{Eq7.7}), and (\ref{Eq7.11}), 
\begin{align}
\left( 2N\right) ^{\frac{1}{2}}\hat{\varphi}\left( 2Nt\right) & =2^{\frac{1}{%
2}}\left( A\left( 2t\right) \hat{\varphi}\left( 2t\right) +B\left( 2t\right) 
\hat{\psi}\left( 2t\right) +\hat{\xi}_{-}\left( 2t\right) \right)
\label{Eq9.17} \\
& =A\left( 2t\right) m_{0}\left( t\right) \hat{\varphi}\left( t\right)
+B\left( 2t\right) m_{1}\left( t\right) \hat{\varphi}\left( t\right) +2^{%
\frac{1}{2}}\hat{\xi}_{-}\left( 2t\right) .  \notag
\end{align}
Now, applying the orthogonal projection onto $\mathcal{\hat{V}}_{0}$ on (\ref
{Eq9.16}) and (\ref{Eq9.17}) and equating the two expressions, we obtain 
\begin{equation*}
m_{0}\left( Nt\right) A\left( t\right) \hat{\varphi}\left( t\right) =\left(
A\left( 2t\right) m_{0}\left( t\right) +B\left( 2t\right) m_{1}\left(
t\right) \right) \hat{\varphi}\left( t\right)
\end{equation*}
for almost all $t\in \mathbb{R}$. But multiplying both sides by $2\pi %
\Bar{\Hat{\varphi}}\left( t\right) $ and adding over all $t:=t+2\pi k$, $%
k\in \mathbb{Z}$, using (\ref{Eq8.7}), we obtain (\ref{Eq9.12}). The formula
(\ref{Eq9.13}) is just a transcription of (\ref{Eq9.12}) (using (\ref{Eq1.16}%
) with $N=2$), and since $S_{0}$ and $S_{1}$ are isometries with orthogonal
ranges, (\ref{Eq9.14}) follows.%
%TCIMACRO{\TeXButton{End Proof}{\end{proof}}}
%BeginExpansion
\end{proof}%
%EndExpansion

\begin{scholium}
\label{Sch9.3}Note in particular that $\mathcal{V}_{0}$ is invariant under
scaling by $N$ if and only if $B=0$ and $\xi _{-}=0$, and then $A\left(
t\right) =M_{0}\left( t\right) $, and (\ref{Eq9.13}) reduces to the relation
(\ref{Eq9.7}). Thus $B$ is a measure of the non-$N$-scale invariance of $%
\mathcal{V}_{0}$.

By Theorem \ref{Thm7.1}, to apply the projection onto $\mathcal{\hat{V}}_{0}$
is equivalent to projecting onto the vectors of the form $\sum_{n=1}^{\infty
}\xi _{n}z^{n}$ in the $z$-transformed Hilbert space. In this space $\xi _{-}
$ has the form 
\begin{equation}
\hat{\xi}_{-}\sim z^{-1}C_{1}+z^{-2}C_{2}+\cdots   \label{Eq9.18}
\end{equation}
while 
\begin{equation}
m_{0}\left( Nt\right) B\left( t\right) \hat{\psi}\left( t\right) \sim
m_{0}\left( Nt\right) B\left( t\right)   \label{Eq9.19}
\end{equation}
by (\ref{Eq7.27}). But (\ref{Eq7.27}) implies that 
\begin{align}
J\left( w\left( \,\cdot \,\right) \hat{\psi}_{n,k}\left( \,\cdot \,\right)
\right) &=J\left( w\left( 2^{-n}2^{n}\,\cdot \,\right) e^{-i2^{n}\!k\,\mathbf{%
\cdot }\,}\,2^{\frac{n}{2}}\hat{\psi}\left( 2^{n}\,\cdot \,\right) \right) 
\label{Eq9.20} \\
&=w\left( 2^{-n}t\right) e^{-ikt}z^{n}  \notag
\end{align}
if $w$ is $2\pi $-periodic and $n=0,-1,-2,\dots $. Thus (\ref{Eq9.18}) means 
\begin{equation}
\hat{\xi}_{-}\left( t\right) =\sum_{n=-1}^{-\infty }C_{-n}\left(
2^{n}t\right) 2^{\frac{n}{2}}\hat{\psi}\left( 2^{n}t\right) ,  \label{Eq9.21}
\end{equation}
and we have 
\begin{equation}
m_{0}\left( Nt\right) \hat{\xi}_{-}\left( t\right) \sim \sum_{n=-1}^{-\infty
}C_{-n}\left( t\right) m_{0}\left( N2^{-n}t\right) z^{n}.  \label{Eq9.22}
\end{equation}
Finally, 
\begin{equation}
2^{\frac{1}{2}}\hat{\xi}_{-}\left( 2t\right) \sim \sum_{n=0}^{-\infty
}C_{1-n}\left( t\right) z^{n}.  \label{Eq9.23}
\end{equation}
Thus it follows from (\ref{Eq9.16}), (\ref{Eq9.17}), (\ref{Eq9.19}), (\ref
{Eq9.22}), and (\ref{Eq9.23}) that 
\begin{equation*}
m_{0}\left( Nt\right) B\left( t\right) =C_{1}\left( t\right) 
\end{equation*}
and 
\begin{equation*}
C_{-n}\left( t\right) m_{0}\left( N2^{-n}t\right) =C_{1-n}\left( t\right) 
\end{equation*}
for $n=-1,-2,\dots $, i.e., 
\begin{equation}
C_{n}\left( t\right) =\prod_{k=0}^{n-1}m_{0}\left( N2^{k}t\right) \cdot
B\left( t\right)   \label{Eq9.24}
\end{equation}
for $n=1,2,3,\dots $. This specifies $\xi _{-}$ as an $A$-dependent operator
applied to $B$, and combining with (\ref{Eq9.3}) we obtain 
\begin{equation}
\lim_{n\rightarrow \infty }2^{-\frac{n}{2}}C_{n}\left( t/\left(
N2^{n}\right) \right) =\left( 2\pi \right) ^{\frac{1}{2}}B\left( 0\right) 
\hat{\varphi}\left( t\right) ,  \label{Eq9.25}
\end{equation}
with convergence in $L^{2}\left( 
%TCIMACRO{
%\TeXButton{RR hat}{\smash{\hat{\mathbb{R}}}\vphantom{\mathbb{R}}}}
%BeginExpansion
\smash{\hat{\mathbb{R}}}\vphantom{\mathbb{R}}%
%EndExpansion
\right) $.

Note also that, Fourier-transforming (\ref{Eq9.10}) using (\ref{Eq9.11}) and
(\ref{Eq9.21}), we obtain the following orthogonal expansion in $L^{2}\left( 
%TCIMACRO{
%\TeXButton{RR hat}{\smash{\hat{\mathbb{R}}}\vphantom{\mathbb{R}}}}
%BeginExpansion
\smash{\hat{\mathbb{R}}}\vphantom{\mathbb{R}}%
%EndExpansion
\right) $: 
\begin{equation}
N^{\frac{1}{2}}\hat{\varphi}\left( Nt\right) =A\left( t\right) \hat{\varphi}%
\left( t\right) +B\left( t\right) \hat{\psi}\left( t\right)
+\sum_{n=1}^{\infty }C_{n}\left( 2^{-n}t\right) 2^{-\frac{n}{2}}\hat{\psi}%
\left( 2^{-n}t\right) ,  \label{Eq9.26}
\end{equation}
and thus 
\begin{equation}
1=\left\| A\right\| _{L^{2}\left( \mathbb{T}\right) }^{2}+\left\| B\right\|
_{L^{2}\left( \mathbb{T}\right) }^{2}+\sum_{n=1}^{\infty }\left\|
C_{n}^{{}}\right\| _{L^{2}\left( \mathbb{T}\right) }^{2}.  \label{Eq9.27}
\end{equation}
\end{scholium}

\section{\label{Fath}\label{AppA}Father functions with non-orthogonal
translates}

It is known (see, e.g., \cite[Section 5.3 and Section 6.2]{Dau92}) that the
multiresolution analysis can be extended to cases where the translates of
the father function $\varphi $ are not exactly orthogonal. In this section
we will consider the case that (\ref{Eq8.1}) is replaced by the weaker
condition that there exists a constant $c>0$ such that 
\begin{equation}
\left\| \sum_{n\in \mathbb{Z}}\xi _{n}\varphi \left( \,\cdot \,-n\right)
\right\| _{L^{2}\left( \mathbb{R}\right) }^{2}\le c\left\| \xi \right\|
_{\ell ^{2}}^{2}  \label{Eq10.1}
\end{equation}
for any sequence $\left( \xi _{n}\right) _{n\in \mathbb{Z}}$ such that only
finitely many components are nonzero. The assumptions (\ref{Eq8.2}) (scale
invariance), (\ref{Eq8.3a}) (refinement), and (\ref{Eq8.3b}) (ergodicity)
are kept as before. By the same reasoning leading to (\ref{Eq8.7}),
condition (\ref{Eq10.1}) is equivalent to 
\begin{equation}
\func{PER}\left( \left| \hat{\varphi}\right| ^{2}\right) \left( t\right) \le 
\frac{c}{2\pi }.  \label{Eq10.2}
\end{equation}
Let us try to establish an analogue of the commutative diagram (\ref{Eq8.29}%
) in this more general setting. We first construct the left side of the
diagram.

\begin{lemma}
\label{Lem10.1}Adopt the assumptions \textup{(\ref{Eq10.1}), (\ref{Eq8.2}),}
and \textup{(\ref{Eq8.3}).} Let $\mu _{\varphi }$ be the measure on $%
\mathbb{T}=\mathbb{R}\diagup 2\pi \mathbb{Z}$ with Radon-Nikodym derivative 
\begin{equation}
\frac{d\mu _{\varphi }}{dt}=\func{PER}\left( \left| \hat{\varphi}\right|
^{2}\right) .  \label{Eq10.3}
\end{equation}
Then there is a one-to-one correspondence between $f\in \mathcal{V}_{0}$ and 
$m\in L^{2}\left( \mathbb{T},\mu _{\varphi }\right) $ given by 
\begin{equation}
\hat{f}\left( t\right) =m\left( e^{-it}\right) \hat{\varphi}\left( t\right) .
\label{Eq10.4}
\end{equation}
Moreover, 
\begin{equation}
\left\| f\right\| _{L^{2}\left( \mathbb{R}\right) }=\left\| m\right\|
_{L^{2}\left( \mathbb{T},\mu _{\varphi }\right) },  \label{Eq10.5}
\end{equation}
i.e., $f\rightarrow m_{f}$ is a unitary operator $\mathcal{V}_{0}\rightarrow
L^{2}\left( \mathbb{T},\mu _{\varphi }\right) $.
\end{lemma}

%TCIMACRO{\TeXButton{Begin Proof}{\begin{proof}}}
%BeginExpansion
\begin{proof}%
%EndExpansion
Assume first that $f$ is a finite linear combination 
\begin{equation}
f\left( \,\cdot \,\right) =\sum_{k\in \mathbb{Z}}a_{k}\varphi \left( \,\cdot
\,-k\right) ,  \label{Eq10.6}
\end{equation}
and put 
\begin{equation}
m_{f}\left( \,\cdot \,\right) =m\left( \,\cdot \,\right)
=\sum_{k}a_{k}e^{-ik\,\mathbf{\cdot }\,}.  \label{Eq10.7}
\end{equation}
Then (\ref{Eq10.4}) is valid, and 
\begin{align*}
\int_{\mathbb{R}}\left| \hat{f}\left( t\right) \right| ^{2}\,dt& =\int_{%
\mathbb{R}}\left| m\left( e^{-it}\right) \hat{\varphi}\left( t\right)
\right| ^{2}\,dt \\
& =\int_{\mathbb{T}}\left| m\right| ^{2}\func{PER}\left( \left| \hat{\varphi}%
\right| ^{2}\right) \,dt \\
& =\int_{\mathbb{T}}\left| m\right| ^{2}\,d\mu _{\varphi },
\end{align*}
so (\ref{Eq10.5}) holds. Since the set of $f$ of the form (\ref{Eq10.6}) is
dense in $\mathcal{V}_{0}$ by the definition of $\mathcal{V}_{0}$, and the
set of $m$ of the form (\ref{Eq10.7}) is dense in $L^{2}\left( \mathbb{T}%
,\mu _{\varphi }\right) $, the Lemma follows by closure.%
%TCIMACRO{\TeXButton{End Proof}{\end{proof}}}
%BeginExpansion
\end{proof}%
%EndExpansion

An immediate corollary is the following:

\begin{lemma}
\label{Lem10.2}Adopt the assumptions \textup{(\ref{Eq10.1}), (\ref{Eq8.2}),}
and \textup{(\ref{Eq8.3}).} Then there is a one-to-one unitary
correspondence between $\psi \in \mathcal{V}_{-n}^{{}}=U_{N}^{-n}\mathcal{V}%
_{0}^{{}}$ and $m=m_{\psi ,n}\in L^{2}\left( \mathbb{T},\mu _{\varphi
}\right) $ given by 
\begin{equation}
N^{\frac{n}{2}}\hat{\psi}\left( N^{n}t\right) =m\left( t\right) \hat{\varphi}%
\left( t\right) .  \label{Eq10.6bis}
\end{equation}
\end{lemma}

%TCIMACRO{\TeXButton{Begin Proof}{\begin{proof}}}
%BeginExpansion
\begin{proof}%
%EndExpansion
We have $\psi \in \mathcal{V}_{-n}^{{}}$ if and only if $U_{N}^{n}\psi \in 
\mathcal{V}_{0}^{{}}$, and $\left\| \psi \right\| _{2}^{{}}=\left\|
U_{N}^{n}\psi \right\| _{2}^{{}}$. Apply Lemma \ref{Lem10.1} on $%
f=U_{N}^{n}\psi $.%
%TCIMACRO{\TeXButton{End Proof}{\end{proof}}}
%BeginExpansion
\end{proof}%
%EndExpansion

Note that (\ref{Eq10.4}) of Lemma \ref{Lem10.1} precisely says that the
diagram 
\begin{equation}
\begin{array}{ccc}
\mathcal{V}_{0} & 
%TCIMACRO{
%\TeXButton{F phi rightleft}{\mkern18mu \smash[b]{\raisebox{0.5ex
%}{\makebox[0pt]{\hss $\displaystyle \overset{\mathcal{F}_\varphi 
%}{\longrightarrow }$\hss }}\raisebox{-0.5ex}{\makebox[0pt]{\hss 
%$\displaystyle \underset{\mathcal{F}_{\varphi }^{-1}}{\longleftarrow 
%}$\hss }}}\mkern18mu } }
%BeginExpansion
\mkern18mu \smash[b]{\raisebox{0.5ex
}{\makebox[0pt]{\hss $\displaystyle \overset{\mathcal{F}_\varphi 
}{\longrightarrow }$\hss }}\raisebox{-0.5ex}{\makebox[0pt]{\hss 
$\displaystyle \underset{\mathcal{F}_{\varphi }^{-1}}{\longleftarrow 
}$\hss }}}\mkern18mu %
%EndExpansion
& L^{2}\left( \mathbb{T},\mu _{\varphi }\right) \\ 
%TCIMACRO{\TeXButton{hook down}{\hookdownarrow } }
%BeginExpansion
\hookdownarrow %
%EndExpansion
&  & 
%TCIMACRO{
%\TeXButton{M phi hat down}{\makebox[0pt]{\hss \rule{0.3pt}{6pt}\hss 
%}\setbox\frogdown=\hbox to 0pt{\hss $\displaystyle \downarrow $\hss 
%}\lower\ht\frogdown\box\frogdown\raisebox{-0.5ex
%}{\makebox[0pt][l]{$\scriptstyle \mkern6mu M_{\hat{\varphi }}$\hss }}} }
%BeginExpansion
\makebox[0pt]{\hss \rule{0.3pt}{6pt}\hss 
}\setbox\frogdown=\hbox to 0pt{\hss $\displaystyle \downarrow $\hss 
}\lower\ht\frogdown\box\frogdown\raisebox{-0.5ex
}{\makebox[0pt][l]{$\scriptstyle \mkern6mu M_{\hat{\varphi }}$\hss }}%
%EndExpansion
\vphantom{\rule[-20pt]{0pt}{36pt}} \\ 
L^{2}\left( \mathbb{R}\right) & 
%TCIMACRO{
%\TeXButton{F rightleft}{\mkern18mu \smash[t]{\raisebox{0.5ex
%}{\makebox[0pt]{\hss $\displaystyle \overset{\mathcal{F}}{\longrightarrow 
%}$\hss }}\raisebox{-0.5ex}{\makebox[0pt]{\hss $\displaystyle 
%\underset{\mathcal{F}^{-1}}{\longleftarrow }$\hss }}}\mkern18mu } }
%BeginExpansion
\mkern18mu \smash[t]{\raisebox{0.5ex
}{\makebox[0pt]{\hss $\displaystyle \overset{\mathcal{F}}{\longrightarrow 
}$\hss }}\raisebox{-0.5ex}{\makebox[0pt]{\hss $\displaystyle 
\underset{\mathcal{F}^{-1}}{\longleftarrow }$\hss }}}\mkern18mu %
%EndExpansion
& L^{2}\left( 
%TCIMACRO{
%\TeXButton{RR hat}{\smash{\hat{\mathbb{R}}}\vphantom{\mathbb{R}}} }
%BeginExpansion
\smash{\hat{\mathbb{R}}}\vphantom{\mathbb{R}}%
%EndExpansion
\right)
\end{array}
\label{Eq10.7bis}
\end{equation}
is commutative, where $\mathcal{F}_{\varphi }$ now is the map $f\rightarrow
m_{f}$, and $M_{\hat{\varphi}}$ still is the map of multiplying the
periodized function by $\hat{\varphi}$. $\mathcal{F}_{\varphi }$ is still
unitary.

If $\psi _{1},\psi _{2}\in \mathcal{V}_{-1}$, and $m_{i}=\mathcal{F}%
_{\varphi }\left( U\psi _{i}\right) $, we compute 
\begin{align}
\ip{\psi _1}{\psi _2}& =\int_{\hat{\mathbb{R}}}\Bar{\Hat{\psi}}_{1}\left(
t\right) \hat{\psi}_{2}\left( t\right) \,dt  \label{Eq10.8} \\
& =\frac{1}{N}\int_{\hat{\mathbb{R}}}\bar{m}_{1}\left( t/N\right) \bar{m}%
_{2}\left( t/N\right) \left| \hat{\varphi}\left( t/N\right) \right| ^{2}\,dt
\notag \\
& =\int_{\hat{\mathbb{R}}}\bar{m}_{1}\left( t\right) m_{2}\left( t\right)
\left| \hat{\varphi}\left( t\right) \right| ^{2}\,dt  \notag \\
& =\int_{\mathbb{T}}\bar{m}_{1}\left( z\right) m_{2}\left( z\right) \,d\mu
_{\varphi }\left( z\right) .  \notag
\end{align}
Correspondingly, 
\begin{align}
\ip{\psi _1}{\smash{T^k\psi _2}}& =\int_{\hat{\mathbb{R}}}\Bar{\Hat{\psi}}%
_{1}\left( t\right) e^{-ikt}\hat{\psi}_{2}\left( t\right) \,dt
\label{Eq10.9} \\
& =\frac{1}{N}\int_{\hat{\mathbb{R}}}\bar{m}_{1}\left( t/N\right) \bar{m}%
_{2}\left( t/N\right) e^{-ikt}\left| \hat{\varphi}\left( t/N\right) \right|
^{2}\,dt  \notag \\
& =\int_{\mathbb{T}}\bar{m}_{1}\left( z\right) m_{2}\left( z\right)
z^{kN}\,d\mu _{\varphi }\left( z\right) ,  \notag
\end{align}
and hence $\ip{\psi _1}{\smash{T^k\psi _2}}=0$ for all $k\in \mathbb{Z}$ if
and only if 
\begin{equation*}
\int_{\mathbb{T}}\bar{m}_{1}\left( z\right) m_{2}\left( z\right) f\left(
z^{N}\right) \,d\mu _{\varphi }\left( z\right) =0
\end{equation*}
for all $f\in L^{\infty }\left( \mathbb{T}\right) $. This is equivalent to 
\begin{equation}
\sum_{k\in \mathbb{Z}_{N}}\bar{m}_{1}\left( \rho ^{k}z\right) m_{2}\left(
\rho ^{k}z\right) \func{PER}\left( \left| \hat{\varphi}\right| ^{2}\right)
\left( \rho ^{k}z\right) =0  \label{Eq10.10}
\end{equation}
for almost all $z$.

From this point one could make a similar theory as in Section \ref{Scal},
replacing (\ref{Eq8.18}) by (\ref{Eq10.10}) and using a selection theorem to
find $m_{1},\dots ,m_{N-1}$, and thus $\psi _{1},\dots ,\psi _{N-1}$. See (%
\ref{Eq12.31})--(\ref{Eq12.32}) below. However, in this case the matrix (\ref
{Eq1.11}) will not be unitary, and thus the connection with representations
of $\mathcal{O}_{N}$ is less direct. Let us rather sketch a completely
different approach, where one starts with functions $m_{0},m_{1},\dots
,m_{N-1}:\mathbb{T}\rightarrow \mathbb{C}$ such that the matrix (\ref{Eq1.11}%
) is assumed to be unitary at the outset. In addition we will assume that $%
m_{0}\left( 0\right) =\sqrt{N}$ and that $m_{0}$ is Lipschitz continuous at $%
0$, or merely that the infinite product 
\begin{equation}
\hat{\varphi}\left( t\right) =\left( 2\pi \right) ^{-\frac{1}{2}%
}\prod_{k=1}^{\infty }\left( m_{0}\left( N^{-k}t\right) /N^{\frac{1}{2}%
}\right)   \label{Eq12.13}
\end{equation}
converges pointwise almost everywhere. By \cite{Mal89}, or \cite[Lemma 6.2.1]
{Dau92}, it follows from the condition $\sum_{k\in \mathbb{Z}_{N}}\left|
m_{0}\left( t+2\pi k/N\right) \right| ^{2}=N$ that $\hat{\varphi}\in
L^{2}\left( \mathbb{R}\right) $ and $\left\| \varphi \right\| _{2}\le 1$. We
will also still assume (\ref{Eq10.1}). If we now \emph{define} $\psi
_{1},\dots ,\psi _{N-1}$ by 
\begin{equation}
\sqrt{N}\hat{\psi}_{k}\left( Nt\right) =m_{k}\left( t\right) \hat{\varphi}%
\left( t\right)   \label{Eq12.14}
\end{equation}
then 
\begin{equation}
U^{n}T^{k}\psi _{m}\left( \,\cdot \,\right) =N^{-\frac{n}{2}}\psi _{m}\left(
N^{-n}\,\cdot \,-k\right) ,  \label{Eq12.15}
\end{equation}
$m=1,\dots ,N-1$, $n,k\in \mathbb{Z}$, no longer forms an orthonormal basis
for $L^{2}\left( \mathbb{R}\right) $, but a \emph{tight frame} in the sense
that 
\begin{equation}
\sum_{n,k,m}\left| \ip{U^{n}T^{k}\psi _{m}}{f}\right| ^{2}=\left\| f\right\|
_{2}^{2}  \label{Eq12.16}
\end{equation}
for all $f\in L^{2}\left( \mathbb{R}\right) $; see \cite[Proposition 6.2.3]
{Dau92}. It is known that a tight frame is an orthonormal basis precisely
when $\left\| \psi _{m}\right\| _{2}=1$ for $m=1,\dots ,N-1$, and in general 
\begin{equation}
f=\sum_{n,k,m}\ip{U^{n}T^{k}\psi _{m}}{f}U^{n}T^{k}\psi _{m};
\label{Eq12.17}
\end{equation}
see \cite[Section 3.2]{Dau92}.

The crucial property used in proving (\ref{Eq12.16}) as in 
\cite[Proposition 6.2.3]{Dau92} is the identity 
\begin{equation}
\sum_{k\in \mathbb{Z}}\left| \ip{U^{n}T^{k}\varphi }{f}\right| ^{2}+\sum_{m=1%
\vphantom{k\in \mathbb{Z}}}^{N-1}\sum_{k\in \mathbb{Z}\vphantom{m=1}}\left| %
\ip{U^{n}T^{k}\psi _{m}}{f}\right| ^{2}=\sum_{k\in \mathbb{Z}\vphantom{m=1}%
}\left| \ip{U^{n-1}T^{k}\varphi }{f}\right| ^{2},  \label{Eq12.17bis}
\end{equation}
which is verified from (\ref{Eq12.14}) and the unitarity of (\ref{Eq1.11}),
and is valid for all $f\in L^{2}\left( \mathbb{R}\right) $. Let us check the
details in the case $N=2$, where (\ref{Eq12.17bis}) takes the form 
\begin{equation}
\sum_{k\in \mathbb{Z}}\left| \ip{\varphi _{n,k}}{f}\right| ^{2}+\sum_{k\in %
\mathbb{Z}}\left| \ip{\smash{\psi _{n,k}}\vphantom{\varphi _{n,k}}}{f}%
\right| ^{2}=\sum_{k\in \mathbb{Z}}\left| \ip{\varphi _{n-1,k}}{f}\right|
^{2},  \label{Eq46.7}
\end{equation}
valid for all $f\in L^{2}\left( \mathbb{R}\right) $. Using the argument from
Theorem \ref{Thm7.1} adjusted as in Lemmas \ref{Lem10.1}--\ref{Lem10.2}, we
note that it is enough to check (\ref{Eq46.7}) for vectors in $U^{-j}%
\mathcal{V}_{0}$ for all $j\in \mathbb{Z}$. Note that by (\ref{Eq12.13}),
the spaces $U^{-j}\mathcal{V}_{0}$ increase as $j\rightarrow \infty $. The
vectors $f\in U^{-j}\mathcal{V}_{0}$ have representations as $\hat{f}\left(
t\right) =2^{-\frac{j}{2}}\left( \xi \hat{\varphi}\right) \left(
t/2^{j}\right) $, where $\xi \in L^{2}\left( \mathbb{T},\mu _{\varphi
}\right) $, according to (\ref{Eq10.4}) or (\ref{Eq10.6bis}). On the
Fourier-transform side the terms in (\ref{Eq46.7}) then take the following
form (we may assume $j\ge n$, and for simplicity, we shall do the
calculation only for $n=0$, and omit the subindex $n$ when $n=0$): 
\begin{align*}
\ip{\varphi _{k}}{f}& =\ip{\smash{\hat{\varphi}_{k}}\vphantom{\varphi
_{k}}}{\smash{\hat{f}}\vphantom{f}} \\
& =\int_{-\infty }^{\infty }\overline{\left( e_{k}\hat{\varphi}\right)
\left( t\right) }2^{-\frac{j}{2}}\left( \xi \hat{\varphi}\right) \left(
t/2^{j}\right) \,dt \\
& =2^{\frac{j}{2}}\int_{-\infty }^{\infty }\overline{\left( e_{k}\hat{\varphi%
}\right) \left( 2^{j}t\right) }\left( \xi \hat{\varphi}\right) \left(
t\right) \,dt \\
& =\int_{-\infty }^{\infty }\overline{m_{0}^{(j)}\left( t\right)
e_{k}^{{}}\left( 2^{j}t\right) }\left| \hat{\varphi}\left( t\right) \right|
^{2}\xi \left( t\right) \,dt \\
& =\ip{S_{0}^{j}e_{k}^{}}{P\xi }_{L^{2}\left( \mathbb{T}\right) },
\end{align*}
where $P=2\pi \func{PER}\left( \left| \hat{\varphi}\right| ^{2}\right) $ and 
$e_{k}\left( t\right) =e^{-ikt}$. By a similar calculation, 
\begin{align*}
\ip{\psi _{k}}{f}& =\ip{S_{0}^{j-1}S_{1}^{}e_{k}^{}}{P\xi }_{L^{2}\left( %
\mathbb{T}\right) } \\
%TCIMACRO{\TeXButton{and}{\intertext{and}}}
%BeginExpansion
\intertext{and}%
%EndExpansion
\ip{\varphi _{-1,k}}{f}& =\ip{S_{0}^{j-1}e_{k}^{}}{P\xi }_{L^{2}\left( %
\mathbb{T}\right) }.
\end{align*}
Substituting back into (\ref{Eq46.7}) and using the fact that $\left\{
e_{k}\right\} $ is an orthonormal basis for $L^{2}\left( \mathbb{T}\right) $%
, we see that (\ref{Eq46.7}) just says that 
\begin{equation*}
\left\| S_{0}^{*\,j}P\xi \right\| _{L^{2}\left( \mathbb{T}\right)
}^{2}+\left\| S_{1}^{*}S_{0}^{*\,j-1}P\xi \right\| _{L^{2}\left( \mathbb{T}%
\right) }^{2}=\left\| S_{0}^{*\,j-1}P\xi \right\| _{L^{2}\left( \mathbb{T}%
\right) }^{2},
\end{equation*}
which in turn takes the form 
\begin{equation*}
PS_{0}^{j-1}\left( S_{0}^{{}}S_{0}^{*}+S_{1}^{{}}S_{1}^{*}\right)
S_{0}^{*\,j-1}P=PS_{0}^{j-1}S_{0}^{*\,j-1}P,
\end{equation*}
and this follows immediately from $%
S_{0}^{{}}S_{0}^{*}+S_{1}^{{}}S_{1}^{*}=I_{L^{2}\left( \mathbb{T}\right) }$.
This proves (\ref{Eq12.17bis}), and thus $\left\{ U^{n}T^{k}\psi
_{m}\right\} $ forms a tight frame. (Compare the present argument to the one
of the proof of Corollary \ref{Cor8.4}, and to (\ref{Eq1.35}) and (\ref
{Eq1.37}) in Section \ref{Intr}.)

\begin{remark}
\label{Rem12.3}An alternative way of defining a commutative diagram like (%
\ref{Eq10.7bis}) is the following. Define a map $\mathcal{F}_{\varphi }:%
\mathcal{V}_{0}\rightarrow L^{2}\left( \mathbb{T},\mu _{\varphi }\right) $, 
\emph{different} from the $\mathcal{F}_{\varphi }$ defined after (\ref
{Eq10.7bis}), by 
\begin{equation}
\left( \mathcal{F}_{\varphi }f\right) \left( e^{-it}\right) =\sum_{k\in %
\mathbb{Z}}\ip{\varphi\left( \,\cdot \,-k\right) }{f}e^{-ikt}.
\label{Eq12.18}
\end{equation}
Then 
\begin{align}
\left\| \mathcal{F}_{\varphi }f\right\| _{2}^{2}&=\sum_{k}\left| %
\ip{\varphi\left( \,\cdot \,-k\right) }{f}\right| ^{2}  \label{Eq12.19} \\
&=\sum_{k}\left| \ip{e^{-ik\,\hat{\mathbf{\cdot }}\,}\hat{\varphi }\left(
\,\hat{\cdot }\,\right) }{m_{f}\left( e^{-i\,\hat{\mathbf{\cdot }}\,}\right)
\hat{\varphi }\left( \,\hat{\cdot }\,\right) }\right| ^{2}  \notag \\
&=\sum_{k}\left| \int_{\mathbb{T}}e^{ikt}m_{f}\left( e^{-ikt}\right) \func{PER%
}\left( \left| \hat{\varphi}\right| ^{2}\right) \left( t\right) \,dt\right|
^{2}  \notag \\
&=\left( 2\pi \right) ^{2}\int_{\mathbb{T}}\left| m_{f}\func{PER}\left(
\left| \hat{\varphi}\right| ^{2}\right) \right| ^{2},  \notag
\end{align}
where the Haar measure $dt/2\pi $ on $\mathbb{T}$ is implicit. On the other
hand, 
\begin{align}
\left\| f\right\| _{2}^{2}&=\int_{\mathbb{R}}\left| m_{f}\hat{\varphi}\right|
^{2}\left( t\right) \,dt  \label{Eq12.20} \\
&=2\pi \int_{\mathbb{T}}\left| m_{f}\right| ^{2}\func{PER}\left( \left| \hat{%
\varphi}\right| ^{2}\right) .  \notag
\end{align}
Thus, using (\ref{Eq10.1}) in the form (\ref{Eq10.2}) we have 
\begin{equation}
\left\| \mathcal{F}_{\varphi }f\right\| _{2}^{2}\le c\left\| f\right\|
_{2}^{2},  \label{Eq12.21}
\end{equation}
so $\mathcal{F}_{\varphi }:\mathcal{V}_{0}\rightarrow L^{2}\left( \mathbb{T}%
,\mu _{\varphi }\right) $ is bounded. Next, define a map $\mathcal{F}%
:L^{2}\left( \mathbb{R}\right) \rightarrow L^{2}\left( 
%TCIMACRO{
%\TeXButton{RR hat}{\smash{\hat{\mathbb{R}}}\vphantom{\mathbb{R}}}}
%BeginExpansion
\smash{\hat{\mathbb{R}}}\vphantom{\mathbb{R}}%
%EndExpansion
\right) $ as a modified Fourier transform: 
\begin{equation}
\mathcal{F}\left( f\right) =2\pi \func{PER}\left( \left| \hat{\varphi}%
\right| ^{2}\right) \hat{f}=P\hat{f},  \label{Eq12.22}
\end{equation}
where $P=2\pi \func{PER}\left( \left| \hat{\varphi}\right| ^{2}\right) $.
The reason for this definition is the following computation, valid for $f\in 
\mathcal{V}_{0}$, i.e., $\hat{f}=m_{f}\hat{\varphi}$: 
\begin{align}
\left( M_{\hat{\varphi}}\mathcal{F}_{\varphi }f\right) \left( \,\hat{\cdot}%
\,\right) &=\sum_{k\in \mathbb{Z}}\ip{\varphi\left( \,\cdot \,-k\right) }{f}%
\hat{\varphi}\left( \,\hat{\cdot}\,\right) e^{-ik\,\hat{\mathbf{\cdot }}\,}
\label{Eq12.23} \\
&=\sum_{k\in \mathbb{Z}}\ip{e^{-ik\,\hat{\mathbf{\cdot }}\,}\hat{\varphi
}}{m_{f}\hat{\varphi }}\hat{\varphi}\left( \,\hat{\cdot}\,\right) e^{-ik\,%
\hat{\mathbf{\cdot }}\,}  \notag \\
&=\sum_{k\in \mathbb{Z}}\left( \int_{\mathbb{R}}e^{ikt}m_{f}\left( t\right)
\left| \hat{\varphi}\right| ^{2}\left( t\right) \,dt\right) \hat{\varphi}%
\left( \,\hat{\cdot}\,\right) e^{-ik\,\hat{\mathbf{\cdot }}\,}  \notag \\
&=\sum_{k\in \mathbb{Z}}\left( \int_{\mathbb{T}}e^{ikt}m_{f}\left( t\right) 
\func{PER}\left( \left| \hat{\varphi}\right| ^{2}\right) \left( t\right)
\,dt\right) e^{-ik\,\hat{\mathbf{\cdot }}\,}\hat{\varphi}\left( \,\hat{\cdot}%
\,\right)   \notag \\
&=m_{f}\left( \,\hat{\cdot}\,\right) P\left( \,\hat{\cdot}\,\right) \hat{%
\varphi}\left( \,\hat{\cdot}\,\right)   \notag \\
&=P\left( \,\hat{\cdot}\,\right) \hat{f}\left( \,\hat{\cdot}\,\right)   \notag
\\
&=\left( \mathcal{F}f\right) \left( \,\hat{\cdot}\,\right) .  \notag
\end{align}
Thus, the following diagram is commutative: 
\begin{equation}
\begin{array}{ccc}
\mathcal{V}_{0} & 
%TCIMACRO{
%\TeXButton{F phi right}{\smash[b]{\overset{\mathcal{F}_\varphi 
%}{\longrightarrow }}} }
%BeginExpansion
\smash[b]{\overset{\mathcal{F}_\varphi }{\longrightarrow }}%
%EndExpansion
& L^{2}\left( \mathbb{T},\mu _{\varphi }\right)  \\ 
%TCIMACRO{\TeXButton{hook down}{\hookdownarrow } }
%BeginExpansion
\hookdownarrow %
%EndExpansion
&  & 
%TCIMACRO{
%\TeXButton{M phi hat down}{\makebox[0pt]{\hss \rule{0.3pt}{6pt}\hss 
%}\setbox\frogdown=\hbox to 0pt{\hss $\displaystyle \downarrow $\hss 
%}\lower\ht\frogdown\box\frogdown\raisebox{-0.5ex
%}{\makebox[0pt][l]{$\scriptstyle \mkern6mu M_{\hat{\varphi }}$\hss }}}}
%BeginExpansion
\makebox[0pt]{\hss \rule{0.3pt}{6pt}\hss 
}\setbox\frogdown=\hbox to 0pt{\hss $\displaystyle \downarrow $\hss 
}\lower\ht\frogdown\box\frogdown\raisebox{-0.5ex
}{\makebox[0pt][l]{$\scriptstyle \mkern6mu M_{\hat{\varphi }}$\hss }}%
%EndExpansion
\vphantom{\rule[-20pt]{0pt}{36pt}} \\ 
L^{2}\left( \mathbb{R}\right)  & 
%TCIMACRO{
%\TeXButton{F right}{\smash[t]{\overset{\mathcal{F}}{\longrightarrow }}} }
%BeginExpansion
\smash[t]{\overset{\mathcal{F}}{\longrightarrow }}%
%EndExpansion
& L^{2}\left( 
%TCIMACRO{
%\TeXButton{RR hat}{\smash{\hat{\mathbb{R}}}\vphantom{\mathbb{R}}}}
%BeginExpansion
\smash{\hat{\mathbb{R}}}\vphantom{\mathbb{R}}%
%EndExpansion
\right) 
\end{array}
\label{Eq12.24}
\end{equation}
This new diagram should not be confused with (\ref{Eq10.7bis}), as the maps
are defined differently. The new maps $\mathcal{F}_{\varphi }$ and $\mathcal{%
F}$ are no longer isometries, but merely continuous, and they are invertible
if and only if there is a lower estimate 
\begin{equation}
b\left\| \xi \right\| _{\ell ^{2}}^{2}\le \left\| \sum_{n\in \mathbb{Z}}\xi
_{n}\varphi \left( \,\cdot \,-n\right) \right\| _{L^{2}\left( \mathbb{R}%
\right) }^{2}  \label{Eq12.25}
\end{equation}
or, equivalently, 
\begin{equation}
\frac{b}{2\pi }\le \func{PER}\left( \left| \hat{\varphi}\right| ^{2}\right)
\left( t\right) .  \label{Eq12.26}
\end{equation}
\end{remark}

\bigskip After this digression, we let $\mathcal{F}$ and $\mathcal{F}%
_{\varphi }$ have the same meaning as in (\ref{Eq10.7bis}) in the rest of
the discussion. For the same reason as above, the map 
\begin{equation}
\func{id}\nolimits_{\varphi }:L^{2}\left( \mathbb{T}\right) \rightarrow
L^{2}\left( \mathbb{T},\mu _{\varphi }\right)   \label{Eq12.27}
\end{equation}
given by $f\rightarrow f$ is bounded and of norm at most $c$ if and only if (%
\ref{Eq10.2}) is valid, and then $\func{id}_{\varphi }^{-1}$ exists as a
bounded operator if and only if (\ref{Eq12.26}) holds.

We now define 
\begin{equation}
J:L^{2}\left( 
%TCIMACRO{
%\TeXButton{RR hat}{\smash{\hat{\mathbb{R}}}\vphantom{\mathbb{R}}} }
%BeginExpansion
\smash{\hat{\mathbb{R}}}\vphantom{\mathbb{R}}%
%EndExpansion
\right) \longrightarrow \mathbb{C}^{N-1}\otimes \mathcal{K}\otimes
L^{2}\left( \mathbb{T}\right)  \label{Eq12.28}
\end{equation}
differently from the $J$ in (\ref{Eq8.28}), by 
\begin{equation}
J\left( \check{f}\right) \left( e^{-it},z\right) =\sum_{m=1\vphantom{n,k\in
\mathbb{Z}}}^{N-1}\sum_{n,k\in \mathbb{Z}\vphantom{m=1}}\ip{U^{n}T^{k}\psi
_{m}}{f}e_{m}\otimes e^{-ikt}\otimes z^{n},  \label{Eq12.29}
\end{equation}
where $\check{f}=\mathcal{F}^{-1}f$ is the inverse Fourier transform of $f$.
The map $J$, so defined, is an isometry because of (\ref{Eq12.16}), but it
is not necessarily surjective, and (\ref{Eq12.29}) coincides with (\ref
{Eq8.28}) if and only if $\left\{ U^{n}T^{k}\psi _{m}\right\} $ is an
orthonormal basis. The intertwining property (\ref{Eq7.28}) also carries
over to the present more general setting, and it takes the form 
\begin{equation}
J\mathcal{F}U=M_{z}J\mathcal{F}.  \label{EqX}
\end{equation}
Let us now do the simple computation of this, omitted in (\ref{Eq7.28}), in
the case $N=2$. Putting $\hat{J}=J\mathcal{F}$, we have 
\begin{equation}
\hat{J}f=\sum_{n\vphantom{k}\in \mathbb{Z}}\sum_{k\vphantom{n}\in \mathbb{Z}}%
\ip{\psi _{n,k}}{f}_{L^{2}\left( \mathbb{R}\right) }e_{k}z^{n}.
\label{Eq46.6}
\end{equation}
We must show 
\begin{equation}
M_{z}\hat{J}=\hat{J}U,  \label{Eq46.8}
\end{equation}
where $M_{z}$ again is the operator on $\mathcal{K}\otimes L^{2}\left( %
\mathbb{T}\right) $ given by $\left( M_{z}\xi \right) \left( z\right) =z\xi
\left( z\right) $, $z\in \mathbb{T}$.

We now check (\ref{Eq46.8}): let $f\in L^{2}\left( \mathbb{R}\right) $. Then
by (\ref{Eq46.6}), using temporarily the notation $e_{k}\left( t\right)
=e^{-ikt}$, 
\begin{align*}
\left( \hat{J}Uf\right) \left( \,\cdot \,,z\right) & =\sum_{n\vphantom{k}%
}\sum_{k\vphantom{n}}\ip{\psi _{n,k}}{Uf}e_{k}z^{n} \\
& =\sum_{n\vphantom{k}}\sum_{k\vphantom{n}}2^{\frac{1+n}{2}}\int_{-\infty
}^{\infty }\overline{\left( e_{k}\hat{\psi}\right) \left( 2^{n}t\right) }%
\hat{f}\left( 2t\right) \,dt\,e_{k}z^{n} \\
& =\sum_{n\vphantom{k}}\sum_{k\vphantom{n}}2^{\frac{n-1}{2}}\int_{-\infty
}^{\infty }\overline{\left( e_{k}\hat{\psi}\right) \left( 2^{n-1}t\right) }%
\hat{f}\left( t\right) \,dt\,e_{k}z^{n} \\
& =\sum_{n\vphantom{k}}\sum_{k\vphantom{n}}\ip{\psi _{n-1,k}}{f}e_{k}z^{n} \\
& =\sum_{n\vphantom{k}}\sum_{k\vphantom{n}}\ip{\psi _{n,k}}{f}e_{k}z^{n+1} \\
& =\left( M_{z}\hat{J}f\right) \left( \,\cdot \,,z\right) ,
\end{align*}
and this completes the proof of (\ref{Eq46.8}).

If the operator $V:\mathbb{C}_{{}}^{N-1}\otimes \mathcal{K}\otimes
H_{+}^{2}\left( \mathbb{T}\right) \rightarrow L_{{}}^{2}\left( \mathbb{T}%
\right) $ is defined exactly as before in (\ref{Eq4.14}) and (\ref{Eq8.26}),
it is still unitary. The diagram corresponding to the right-hand side of (%
\ref{Eq8.29}) is 
\begin{equation}
\begin{array}{ccc}
L^{2}\left( \mathbb{T},\mu _{\varphi }\right)  & 
%TCIMACRO{
%\TeXButton{id phi left}{\smash[b]{\xleftarrow{\func{id}_{\varphi }}}}}
%BeginExpansion
\smash[b]{\xleftarrow{\func{id}_{\varphi }}}%
%EndExpansion
\quad \mathcal{K}=L^{2}\left( \mathbb{T}\right) \quad 
%TCIMACRO{\TeXButton{V left}{\smash[b]{\overset{V}{\longleftarrow }}} }
%BeginExpansion
\smash[b]{\overset{V}{\longleftarrow }}%
%EndExpansion
& \mathbb{C}_{{}}^{N-1}\otimes \mathcal{K}\otimes H_{+}^{2}\left( \mathbb{T}%
\right)  \\ 
%TCIMACRO{
%\TeXButton{M phi hat down}{\makebox[0pt]{\hss \rule{0.3pt}{6pt}\hss 
%}\setbox\frogdown=\hbox to 0pt{\hss $\displaystyle \downarrow $\hss 
%}\lower\ht\frogdown\box\frogdown\raisebox{-0.5ex
%}{\makebox[0pt][l]{$\scriptstyle \mkern6mu M_{\hat{\varphi }}$\hss }}}}
%BeginExpansion
\makebox[0pt]{\hss \rule{0.3pt}{6pt}\hss 
}\setbox\frogdown=\hbox to 0pt{\hss $\displaystyle \downarrow $\hss 
}\lower\ht\frogdown\box\frogdown\raisebox{-0.5ex
}{\makebox[0pt][l]{$\scriptstyle \mkern6mu M_{\hat{\varphi }}$\hss }}%
%EndExpansion
\vphantom{\rule[-20pt]{0pt}{36pt}} &  & 
%TCIMACRO{\TeXButton{hook down}{\hookdownarrow } }
%BeginExpansion
\hookdownarrow %
%EndExpansion
\\ 
L^{2}\left( 
%TCIMACRO{
%\TeXButton{RR hat}{\smash{\hat{\mathbb{R}}}\vphantom{\mathbb{R}}}}
%BeginExpansion
\smash{\hat{\mathbb{R}}}\vphantom{\mathbb{R}}%
%EndExpansion
\right)  & 
%TCIMACRO{
%\TeXButton{long J right}{\settowidth{\qedskip}{$\displaystyle 
%\smash[b]{\xleftarrow{\func{id}_{\varphi }}}\quad \mathcal{K}=L^{2}\left( 
%\mathbb{T}\right) \quad \smash[b]{\overset{V}{\longleftarrow 
%}}$}\settowidth{\Jwidth}{$\scriptstyle J$}\addtolength{\qedskip
%}{-\Jwidth}\smash[t]{\xrightarrow{\kern0.5\qedskip J\kern0.5\qedskip }}} }
%BeginExpansion
\settowidth{\qedskip}{$\displaystyle 
\smash[b]{\xleftarrow{\func{id}_{\varphi }}}\quad \mathcal{K}=L^{2}\left( 
\mathbb{T}\right) \quad \smash[b]{\overset{V}{\longleftarrow 
}}$}\settowidth{\Jwidth}{$\scriptstyle J$}\addtolength{\qedskip
}{-\Jwidth}\smash[t]{\xrightarrow{\kern0.5\qedskip J\kern0.5\qedskip }}%
%EndExpansion
& \mathbb{C}^{N-1}\otimes \mathcal{K}\otimes L^{2}\left( \mathbb{T}\right) 
\end{array}
\label{Eq12.30}
\end{equation}
This diagram is necessarily not commutative, however, unless $\func{PER}%
\left( \left| \hat{\varphi}\right| ^{2}\right) =\frac{1}{2\pi }$, i.e., (\ref
{Eq8.1}) is fulfilled. The reason is that the maps $V$, $M_{\hat{\varphi}}$, 
$J$, and the inclusion map are all isometries, while $\func{id}_{\varphi }$
is not unless $\func{PER}\left( \left| \hat{\varphi}\right| ^{2}\right) =%
\frac{1}{2\pi }$. In order to make the diagram commutative, $V$ would have
to be redefined. One way would be to choose $m_{1},\dots ,m_{N-1}\in
L^{\infty }\left( \mathbb{T}\right) $ such that the relations 
\begin{equation}
\sum_{k\in \mathbb{Z}_{N}}\bar{m}_{i}\left( \rho ^{k}z\right) m_{j}\left(
\rho ^{k}z\right) \func{PER}\left( \left| \hat{\varphi}\right| ^{2}\right)
\left( \rho ^{k}z\right) 2\pi =N\delta _{ij}  \label{Eq12.31}
\end{equation}
are valid for almost all $z$ (see (\ref{Eq10.8})--(\ref{Eq10.10})), and then
define $S_{k}$ on $L^{2}\left( \mathbb{T};\mu _{\varphi }\right) $ by 
\begin{equation}
\left( S_{k}\xi \right) \left( z\right) =m_{k}\left( z\right) \xi \left(
z^{N}\right) .  \label{Eq12.32}
\end{equation}
One verifies that this defines a representation of $\mathcal{O}_{N}$, and
hence the corresponding $V$ is a unitary. These remaining details for making
a commutative variant of the diagram will be published in a forthcoming
paper.

Let us now consider further the orthogonality properties of the $\mathbb{Z}$%
-translates $\left\{ \varphi \left( \,\cdot \,-k\right) \right\} _{k\in %
\mathbb{Z}}$ in $L^{2}\left( \mathbb{R}\right) $.

Let $\varphi \in L^{2}\left( \mathbb{R}\right) $, and assume that $\varphi $
can be expanded like 
\begin{equation}
U\varphi =\sum_{k}a_{k}\varphi \left( \,\cdot \,-k\right) ,  \label{EqA.1}
\end{equation}
where $\sum_{k}\left| a_{k}\right| ^{2}<\infty $, and where 
\begin{equation}
\left( U\varphi \right) \left( x\right) =N^{-\frac{1}{2}}\varphi \left(
x/N\right) .  \label{EqA.2}
\end{equation}
Thus 
\begin{equation}
\sqrt{N}\hat{\varphi}\left( Nt\right) =m_{0}\left( t\right) \hat{\varphi}%
\left( t\right) ,  \label{EqA.3}
\end{equation}
where 
\begin{equation}
m_{0}\left( t\right) =\sum_{k}a_{k}e^{-ikt},  \label{EqA.4}
\end{equation}
so $m_{0}\in L^{2}\left( \mathbb{T}\right) $.

From now and through the rest of this section, we will make the overall
assumption that $m_{0}$ is uniformly Lipschitz continuous, i.e., there
exists a $K>0$ such that $\left| m_{0}\left( t\right) -m_{0}\left( s\right)
\right| \le K\left| t-s\right| $ for all $t,s\in \mathbb{R}$. This condition
is for example implied by the stronger condition $\sum_{k}\left|
ka_{k}\right| <\infty $ which is much used in \cite{Dau92}. We may then
define an operator $R:C\left( \mathbb{T}\right) \rightarrow C\left( %
\mathbb{T}\right) $ by 
\begin{equation}
\left( R\xi \right) \left( z\right) =\frac{1}{N}\sum_{%
\substack{ w \\ \makebox[0pt]{\hss {$\scriptstyle w^{N}=z
$}\hss }}}\left| m_{0}\left( w\right) \right| ^{2}\xi \left( w\right) .
\label{EqA.5}
\end{equation}

\begin{proposition}
\label{ProA.1}If $m_{0}$ is uniformly Lipschitz continuous, $\hat{\varphi}$
is continuous at $0$, and $\varphi $ and $m_{0}$ are normalized by $\hat{%
\varphi}\left( 0\right) =\left( 2\pi \right) ^{-\frac{1}{2}}$ and $%
m_{0}\left( 0\right) =N^{\frac{1}{2}}$ so that $\hat{\varphi}\left( t\right)
=\left( 2\pi \right) ^{-\frac{1}{2}}\prod_{k=1}^{\infty }\left( N^{-\frac{1}{%
2}}m_{0}\left( tN^{-k}\right) \right) $, then the following conditions are
equivalent. 
%TCIMACRO{
%\TeXButton{Custom Aligned Equation}{\begin{align}
%\begin{minipage}[t]{\displayboxwidth}\raggedright $\left\{ \varphi 
%\left( \,\cdot \,-k\right) \right\} $ is an orthonormal
%set.\end{minipage}  \label{EqA.6} \\
%\begin{minipage}[t]{\displayboxwidth}\raggedright $\func{PER}\left( 
%\left| \hat{\varphi}\right| ^{2}\right) =\left( 2\pi
%\right) ^{-1}\openone$. \end{minipage}  \label{EqA.7} \\
%\begin{minipage}[t]{\displayboxwidth}\raggedright Up to a scalar, $\openone
%$ is the unique eigenvector of $R$ of eigenvalue $1$.\end{minipage}
%\label{EqA.8}
%\end{align}}}
%BeginExpansion
\begin{align}
\begin{minipage}[t]{\displayboxwidth}\raggedright $\left\{ \varphi 
\left( \,\cdot \,-k\right) \right\} $ is an orthonormal
set.\end{minipage}  \label{EqA.6} \\
\begin{minipage}[t]{\displayboxwidth}\raggedright $\func{PER}\left( 
\left| \hat{\varphi}\right| ^{2}\right) =\left( 2\pi
\right) ^{-1}\openone$. \end{minipage}  \label{EqA.7} \\
\begin{minipage}[t]{\displayboxwidth}\raggedright Up to a scalar, $\openone
$ is the unique eigenvector of $R$ of eigenvalue $1$.\end{minipage}
\label{EqA.8}
\end{align}%
%EndExpansion
\end{proposition}

%TCIMACRO{\TeXButton{Begin Proof}{\begin{proof}}}
%BeginExpansion
\begin{proof}%
%EndExpansion
We already proved the implications (\ref{EqA.6})$\Leftrightarrow $(\ref
{EqA.7})$\Leftarrow $(\ref{EqA.8}) in the remarks around (\ref{Eq8.6})--(\ref
{Eq8.10}). In particular, (\ref{EqA.8}) and (\ref{Eq8.8}) imply that $\func{%
PER}\left( \left| \hat{\varphi}\right| ^{2}\right) $ is a scalar multiple of 
$1$, and hence the $\left\{ \varphi \left( \,\cdot \,-k\right) \right\} $
are orthogonal, but then, as 
\begin{equation*}
\frac{1}{N}\sum_{%
\substack{ w \\ \makebox[0pt]{\hss {$\scriptstyle w^{N}=z
$}\hss }}}\left| m_{0}\left( w\right) \right| ^{2}=1
\end{equation*}
as a consequence of (\ref{EqA.8}) and $m_{0}\left( 0\right) =N^{\frac{1}{2}}$%
, we have $m_{0}\left( 2\pi \frac{k}{N}\right) =0$ for $k=1,\dots ,N-1$, and
it follows from the product expansion of $\hat{\varphi}$ that $\hat{\varphi}%
\left( 2\pi k\right) =0$ for all $k\in \mathbb{Z}\setminus \left\{ 0\right\} 
$. Thus $\func{PER}\left( \left| \hat{\varphi}\right| ^{2}\right) \left(
0\right) =\left( 2\pi \right) ^{-1}$, and as $\func{PER}\left( \left| \hat{%
\varphi}\right| ^{2}\right) $ is a scalar multiple of $\openone$, (\ref
{EqA.7}) follows. Implicit in this reasoning is that $\func{PER}\left(
\left| \hat{\varphi}\right| ^{2}\right) $ is continuous on $\mathbb{T}$, but
this follows from the uniform Lipschitz condition on $m_{0}$ and the product
expansion. It remains to prove (\ref{EqA.6})$\Rightarrow $(\ref{EqA.8}). It
follows from (\ref{EqA.7}) ($\Leftrightarrow $(\ref{EqA.6})) and (\ref{Eq8.8}%
) that $\openone$ is indeed an eigenvector of $R$ of eigenvalue $1$, and it
remains to show that it is the only one. (See Remark \ref{RemAA.3}.)

\begin{lemma}
\label{LemAA.1}Let $N\in \mathbb{N}$, $N\ge 2$, and $m_{0}\in L^{\infty
}\left( \mathbb{T}\right) $ be given, satisfying 
\begin{equation}
\sum_{\substack{ w \\ \makebox[0pt]{\hss {$\scriptstyle w^{N}=z
$}\hss }}}\left| m_{0}\left( w\right) \right| ^{2}=N,  \label{EqAA.1}
\end{equation}
for almost all $\,z\in \mathbb{T}$. Let $S_{0}$ be the corresponding
isometry of $L^{2}\left( \mathbb{T}\right) $, 
\begin{equation*}
\left( S_{0}f\right) \left( z\right) =m_{0}\left( z\right) f\left(
z^{N}\right) ,\quad f\in L^{2}\left( \mathbb{T}\right) .
\end{equation*}
Let $m_{0}$ and $\varphi $ further satisfy the general conditions in
Proposition \textup{\ref{ProA.1}.} Then the orthogonality condition \textup{(%
\ref{EqA.6})} for $\varphi $ in $L^{2}\left( \mathbb{R}\right) $ is
equivalent to 
\begin{equation}
\lim_{n\rightarrow \infty }\ip{\openone }{S_0^{*\,n}M_fS_0^n\openone }%
_{L^{2}\left( \mathbb{T}\right) }=f\left( 0\right)   \label{EqAA.2}
\end{equation}
for all $f\in C\left( \mathbb{T}\right) =C\left( \mathbb{R}\diagup 2\pi %
\mathbb{Z}\right) $. Thus the two conditions \textup{(\ref{EqAA.1})} and 
\textup{(\ref{EqAA.2})} together are equivalent to the other conditions in
Proposition \textup{\ref{ProA.1}.}
\end{lemma}

In general, if we do not assume the orthogonality (\ref{EqA.6}) but merely
its consequence (\ref{EqAA.1}), the left-hand limit in (\ref{EqAA.2}) exists
and defines a probability measure $D$ on $\mathbb{T}$. The fact that a Borel
measure $D$ is defined as 
\begin{align}
D\left( f\right) & =\lim_{n\rightarrow \infty }\ip{\openone
}{S_0^{*\,n}M_fS_0^n\openone }  \label{EqAA.2bis} \\
& =\lim_{n\rightarrow \infty }\int_{\mathbb{T}}f\,d\nu _{n}  \notag
\end{align}
is justified in the discussion of Remark \ref{RemAA.2} below. It is a
compactness argument, referring to the Hausdorff metric on the Borel
probability measures on $\mathbb{T}$, and it requires the Lipschitz
assumption on the function $m_{0}$; see also \cite{Hut81} for definitions.
If $\mu $ and $\nu $ are Borel measures on $\mathbb{T}=\mathbb{R}\diagup
2\pi \mathbb{Z}$, then the Hausdorff metric $d_{H}$ is 
\begin{equation*}
d_{H}\left( \mu ,\nu \right) :=\sup_{f}\left\{ \left| \int_{\mathbb{T}%
}f\,d\mu -\int_{\mathbb{T}}f\,d\nu \right| \biggm|f\in C^{1}\left( \mathbb{T}%
\right) ,\sup_{t}\left| f^{\prime }\left( t\right) \right| \le 1\right\}
\end{equation*}
with $C^{1}$-functions on $\mathbb{T}$ identified with differentiable $2\pi $%
-periodic functions on $\mathbb{R}$. The approximation $\nu _{n}\rightarrow
D $ in (\ref{EqAA.2bis}) refers to $\lim_{n\rightarrow \infty }d_{H}\left(
\nu _{n},D\right) =0$.

%TCIMACRO{
%\TeXButton{Begin Proof of Lemma AA.1}
%{\begin{proof}[Proof of Lemma \upshape\ref{LemAA.1}]}}
%BeginExpansion
\begin{proof}[Proof of Lemma \upshape\ref{LemAA.1}]%
%EndExpansion
We use the result of Meyer and Paiva \cite{MePa93} mentioned in (\ref{Eq1.45}%
) to the effect that (\textup{\ref{EqA.6}}) is equivalent in turn to 
\begin{equation}
\int_{\delta \le \left| t\right| \le \pi }P_{n}\left( t\right) \,dt\underset{%
n\rightarrow \infty }{\longrightarrow }0,  \label{EqAA.3}
\end{equation}
for all positive $\delta $, where (by identification) $m_{0}\left( t\right)
\sim m_{0}\left( e^{-it}\right) $, and 
\begin{equation*}
P_{n}\left( t\right) :=\left| m_{0}\left( t\right) m_{0}\left( Nt\right)
\cdots m_{0}\left( N^{n-1}t\right) \right| ^{2}.
\end{equation*}
The lemma follows from Meyer--Paiva using 
\begin{align*}
\left( S_{0}^{n}\openone\right) \left( z\right) & =m_{0}^{{}}\left( z\right)
m_{0}^{{}}\left( z^{N}\right) \cdots m_{0}^{{}}\left( z_{{}}^{N^{n-1}}\right)
\\
& =:m_{0}^{(n)}\left( z\right) , \\
P_{n}^{{}}\left( t\right) & =\left| m_{0}^{(n)}\left( e_{{}}^{-it}\right)
\right| ^{2}, \\
%TCIMACRO{\TeXButton{thus,}{\intertext{thus,}} }
%BeginExpansion
\intertext{thus,}%
%EndExpansion
\ip{S_{0}^{n}\openone }{fS_{0}^{n}\openone }& =\ip{\openone
}{S_{0}^{*\,n}M_{f}^{}S_{0}^{n}\openone } \\
& =\int_{\mathbb{T}}\left( R^{n}f\right) \left( z\right) \frac{\left|
dz\right| }{2\pi } \\
& =\int_{\mathbb{T}}\left| m_{0}^{(n)}\left( z\right) \right| ^{2}f\left(
z\right) \frac{\left| dz\right| }{2\pi } \\
& =\int_{\mathbb{T}}P_{n}\left( z\right) f\left( z\right) \frac{\left|
dz\right| }{2\pi },
\end{align*}
and an elementary characterization of the Dirac mass at $z=1$.%
%TCIMACRO{\TeXButton{End Proof}{\end{proof}}}
%BeginExpansion
\end{proof}%
%EndExpansion

To prove (\ref{EqA.6})$\Rightarrow $(\ref{EqA.8}), it thus suffices to show
that (\ref{EqAA.1}) and (\ref{EqAA.2}) imply (\ref{EqA.8}). To this end, one
easily deduces from (\ref{Eq3.19}) that 
\begin{equation}
S_{0}^{*\,n}M_{f}^{{}}S_{0}^{n}=M_{R^{n}f}^{{}},  \label{EqAA.4}
\end{equation}
where still 
\begin{equation*}
\left( Rf\right) \left( z\right) =N^{-1}\sum_{%
\substack{ w \\ \makebox[0pt]{\hss {$\scriptstyle w^{N}=z
$}\hss }}}\left| m_{0}\left( w\right) \right| ^{2}f\left( w\right) .
\end{equation*}
We conclude from (\ref{EqAA.2}) that 
\begin{equation}
f\left( 0\right) =\lim_{n\rightarrow \infty }\int_{\mathbb{T}}\left(
R^{n}f\right) \left( z\right) \frac{\left| dz\right| }{2\pi }  \label{EqAA.5}
\end{equation}
for all $f\in C\left( \mathbb{T}\right) $. Note the formula 
\begin{equation}
\left( R^{n}f\right) \left( z\right) =N^{-n}\sum_{%
\substack{ w \\ \makebox[0pt]{\hss {$\scriptstyle w^{N^n}=z
$}\hss }}}P_{n}\left( w\right) f\left( w\right) .  \label{Rnfz}
\end{equation}

It follows from \textup{(\ref{EqAA.5})} that 
\begin{equation}
f\left( 0\right) =\int_{\mathbb{T}}f\left( z\right) \frac{\left| dz\right| }{%
2\pi }  \label{f0if}
\end{equation}
if $f$ satisfies 
\begin{equation}
Rf=f.  \label{EqAA.8}
\end{equation}
From (\ref{EqAA.8}), we obtain by the Schwarz inequality (see 
\cite[Notes and Remarks to Section 5.3.1]{BrRoII}) applied to $R$: 
\begin{equation*}
\left| f\right| ^{2}=\left| Rf\right| ^{2}\le R\left( \left| f\right|
^{2}\right) .
\end{equation*}
By induction, then, 
\begin{equation*}
\left| f\right| ^{2}\le R^{n}\left( \left| f\right| ^{2}\right) \le
R^{n+1}\left( \left| f\right| ^{2}\right) ,
\end{equation*}
and therefore, from (\ref{EqAA.5}), 
\begin{equation*}
\left| f\left( 0\right) \right| ^{2}=\lim_{n\rightarrow \infty }\int_{%
\mathbb{T}}R^{n}\left( \left| f\right| ^{2}\right) \ge \int_{\mathbb{T}%
}\left| f\right| ^{2},
\end{equation*}
where the Haar measure $\frac{\left| dz\right| }{2\pi }$ is implicitly
understood. Using Cauchy-Schwarz and (\ref{f0if}), we then conclude that 
\begin{equation*}
\left| f\left( 0\right) \right| ^{2}=\left| \int_{\mathbb{T}}f\right|
^{2}\le \int_{\mathbb{T}}\left| f\right| ^{2}\le \left| f\left( 0\right)
\right| ^{2}\text{.\quad (Recall }z=e^{-it}\text{.)}
\end{equation*}
We conclude that the Cauchy--Schwarz inequality is an equality when applied
to the two functions $f$ and $\openone$. Hence $f$ is a constant multiple of 
$\openone
$, and we have proved that the eigenspace of $R$ corresponding to the
eigenvalue $1$ is one-dimensional. This ends the proof of Proposition \ref
{ProA.1}.%
%TCIMACRO{\TeXButton{End Proof}{\end{proof}}}
%BeginExpansion
\end{proof}%
%EndExpansion

\begin{remark}
\label{RemAA.2}Using (\ref{EqAA.1}), we conclude that the measure $D$
defined in general by (\ref{EqAA.2bis}) satisfies the two invariance
properties below (\ref{EqAA.6})--(\ref{EqAA.7}), even when (\ref{EqA.6}) is
not assumed: 
\begin{equation}
D\left( R\left( f\right) \right) =D\left( f\right)   \label{EqAA.6}
\end{equation}
for all $f\in C\left( \mathbb{T}\right) $, and 
\begin{equation}
D\left( \sigma \left( f\right) \right) =D\left( f\right)   \label{EqAA.7}
\end{equation}
where $\sigma \left( f\right) \left( z\right) =f\left( z^{N}\right) $. Since 
$\sigma $ is mixing (see \cite{Kea72}), the measure $D$ (in the wavelet
examples) must be singular, but with support on $\mathbb{T}$ invariant under 
$\sigma $.

The measure $D$ defined by (\ref{EqAA.2bis}) exists by the
Ruelle-Perron-Frobenius theorem in the form of \cite[p.\ 21]{PaPo90}, at
least if $m_{0}\left( z\right) \ne 0$ for all $z$, except for a finite
number of zeroes of $m_{0}$; see \cite{Kea72} or \cite{Hut81}. The required
regularity condition on $m_{0}=\sum_{k\in \mathbb{Z}}a_{k}z^{k}$ is $%
\sum_{k\in \mathbb{Z}}\left| ka_{k}\right| <\infty $, which also guarantees
convergence of the infinite product formula for $\hat{\varphi}$. Since $%
\left| m_{0}\left( z\right) \right| ^{2}P_{n}\left( z^{N}\right)
=P_{n+1}\left( z\right) $, the invariance (\ref{EqAA.7}) of the measure $D$
follows: specifically, 
\begin{align*}
\int_{\mathbb{T}}P_{n+1}\left( z\right) f\left( z^{N}\right) \frac{\left|
dz\right| }{2\pi }& =\int_{\mathbb{T}}\left| m_{0}\left( z\right) \right|
^{2}P_{n}\left( z^{N}\right) f\left( z^{N}\right) \frac{\left| dz\right| }{%
2\pi } \\
& =\frac{1}{N}\int_{\mathbb{T}}\sum_{%
\substack{ w \\ \makebox[0pt]{\hss {$\scriptstyle w^{N}=z
$}\hss }}}\left| m_{0}\left( w\right) \right| ^{2}P_{n}\left( z\right)
f\left( z\right) \frac{\left| dz\right| }{2\pi } \\
& =\int_{\mathbb{T}}P_{n}\left( z\right) f\left( z\right) \frac{\left|
dz\right| }{2\pi },
\end{align*}
so (\ref{EqAA.7}) follows upon taking the $n\rightarrow \infty $ limit.

Using the estimate 
\begin{equation}
N^{n}\left| \hat{\varphi}\left( N^{n}t\right) \right| ^{2}\le \frac{1}{2\pi }%
P_{n}\left( e^{-it}\right) ,  \label{Nnphihatnormestimate}
\end{equation}
which follows from (\ref{Eq12.13}), (\ref{EqAA.1}), and the definition of $%
P_{n}$ after (\ref{EqAA.3}), we will argue that 
\begin{equation}
\left\| \varphi \right\| _{L^{2}\left( \mathbb{R}\right) }^{2}\le D\left(
\left\{ 0\right\} \right) .  \label{kphinormbound}
\end{equation}
Integrating (\ref{Nnphihatnormestimate}) over $\left\langle -\eta ,\eta
\right\rangle $, where $\eta >0$, we obtain
\begin{equation}
\int_{-\eta N^{n}}^{\eta N^{n}}\left| \hat{\varphi}\left( t\right) \right|
^{2}\,dt\le \int_{\left| t\right| <\eta }P_{n}\left( t\right) \,dt.
\label{phihatintbound}
\end{equation}
The limit on the left-hand side (as $n\rightarrow \infty $) is 
\begin{equation*}
\int_{\mathbb{R}}\left| \hat{\varphi}\right| _{{}}^{2}\,dt=\left\| \varphi
\right\| _{L^{2}\left( \mathbb{R}\right) }^{2},
\end{equation*}
and on the right it is $D\left( \left\langle -\eta ,\eta \right\rangle
\right) $ by formula (\ref{EqAA.2bis}) applied to $f=\chi _{\left\langle
-\eta ,\eta \right\rangle }$. Letting $\eta \rightarrow 0$, (\ref
{kphinormbound}) follows. Using (\ref{EqAA.6}) and (\ref{EqAA.7}), we note
that when $\func{supp}\left( D\right) $ is finite, then there are cycles
corresponding to roots $a\in \mathbb{T}$ of $a^{N^{k}}=a$, $k=1,2,\dots $ ($k
$ chosen minimal), such that $D$ is a convex combination of associated
measures $D_{a}$ defined as 
\begin{equation*}
D_{a}:=\frac{1}{k}\sum_{j=0}^{k-1}\delta _{a^{N^{j}}},
\end{equation*}
where $\delta _{a^{N^{j}}}$ denotes the Dirac mass at the point $z=a^{N^{j}}$
on $\mathbb{T}$. The case $D_{1}=\delta _{1}$ may occur in the convex
combination because $\left| m_{0}\left( 1\right) \right| =\sqrt{N}$ in our
case, referring to the $z$-parameter on $\mathbb{T}$. By (\ref{EqAA.6}), a
cycle $D_{a}$ may in general occur in the convex combination for $D$ iff $%
\left| m_{0}\left( a^{N^{j}}\right) \right| =\sqrt{N}$ for $j=0,1,\dots
,N^{k-1}$.

Recall that we have encountered the functions $P_{n}\left( t\right) $ before
in a situation where the normalization $m_{0}\left( 0\right) =N^{\frac{1}{2}}
$ is not fulfilled, in the proof of Lemma \ref{Lem3.2}.

The finite-orbit picture for the $z\mapsto z^{N}$ action on $\limfunc{supp}%
\left( D\right) $ ($\subset \mathbb{T}$), and its connection to the Cohen
cycles (see \cite{Coh90} and \cite[Theorem 6.3.3, p.\ 188]{Dau92}), will be
taken up in a subsequent paper. This decomposition is also closely connected
(in a special case) to one which arises in an earlier paper of ours 
\cite[Proposition 8.2]{BrJo96b}.

In a forthcoming paper, we plan to study the other possibilities for $%
\limfunc{supp}\left( D\right) $: possibly infinite, possibly allowing
infinite orbits, or an infinite number of finite orbits, under the
restricted action of $z\mapsto z^{N}$ on $\limfunc{supp}\left( D\right) $.
\end{remark}

\begin{remark}
\label{RemAA.3}Let us check how much mileage we can get towards the proof of
(\ref{EqA.6})$\Rightarrow $(\ref{EqA.8}) in Proposition \ref{ProA.1} without
using Lemma \ref{LemAA.1}. For this, let $\xi \in L^{\infty }\left( %
\mathbb{T}\right) $ and $\eta \in L^{2}\left( \mathbb{T}\right) $, and
compute 
%TCIMACRO{
%\TeXButton{Annotated Aligned Equation}{\begin{align*}
%\int_{\mathbb{T}}\left( R\xi \right) \left( z\right) \eta \left( z\right) 
%\frac{\left| dz\right| }{2\pi }& =\int_{\mathbb{T}}\frac{1}{N}
%\sum_{\substack{ w \\ \makebox[0pt]{\hss {$\scriptstyle w^{N}=z
%$}\hss }}}\left| m_{0}\left( w\right) \right| ^{2}\xi \left( w\right) \eta
%\left( w^{N}\right) \frac{\left| dz\right| }{2\pi } & & \\
%& =\int_{\mathbb{T}}\left| m_{0}\left( z\right) \right| ^{2}\xi \left(
%z\right) \eta \left( z^{N}\right) \frac{\left| dz\right| }{2\pi } 
%& &\text{(by (\ref{Eq3.0}))} \\
%& =\int_{\mathbb{R}}\left| m_{0}\left( t\right) \right| ^{2}\left| 
%\hat{\varphi}\left( t\right) \right| ^{2}\xi \left( t\right) \eta \left(
%Nt\right) \,dt & &\text{(by (\ref{EqA.7}))} \\
%& =N\int_{\mathbb{R}}\xi \left( t\right) \left| \hat{\varphi}\left(
%Nt\right) \right| ^{2}\eta \left( Nt\right) \,dt 
%& &\text{(by (\ref{EqA.3}))} \\
%& =\int_{\mathbb{R}}\xi \left( t/N\right) \left| \hat{\varphi}\left(
%t\right) \right| ^{2}\eta \left( t\right) \,dt. & &
%\end{align*}}}
%BeginExpansion
\begin{align*}
\int_{\mathbb{T}}\left( R\xi \right) \left( z\right) \eta \left( z\right) 
\frac{\left| dz\right| }{2\pi }& =\int_{\mathbb{T}}\frac{1}{N}
\sum_{\substack{ w \\ \makebox[0pt]{\hss {$\scriptstyle w^{N}=z
$}\hss }}}\left| m_{0}\left( w\right) \right| ^{2}\xi \left( w\right) \eta
\left( w^{N}\right) \frac{\left| dz\right| }{2\pi } & & \\
& =\int_{\mathbb{T}}\left| m_{0}\left( z\right) \right| ^{2}\xi \left(
z\right) \eta \left( z^{N}\right) \frac{\left| dz\right| }{2\pi } 
& &\text{(by (\ref{Eq3.0}))} \\
& =\int_{\mathbb{R}}\left| m_{0}\left( t\right) \right| ^{2}\left| 
\hat{\varphi}\left( t\right) \right| ^{2}\xi \left( t\right) \eta \left(
Nt\right) \,dt & &\text{(by (\ref{EqA.7}))} \\
& =N\int_{\mathbb{R}}\xi \left( t\right) \left| \hat{\varphi}\left(
Nt\right) \right| ^{2}\eta \left( Nt\right) \,dt 
& &\text{(by (\ref{EqA.3}))} \\
& =\int_{\mathbb{R}}\xi \left( t/N\right) \left| \hat{\varphi}\left(
t\right) \right| ^{2}\eta \left( t\right) \,dt. & &
\end{align*}%
%EndExpansion
If $\xi $ is an eigenvector for $R$ with eigenvalue $1$, this computation
gives 
\begin{equation*}
\int_{\mathbb{R}}\left( \xi \left( t\right) -\xi \left( t/N\right) \right)
\left| \hat{\varphi}\left( t\right) \right| ^{2}\eta \left( t\right) \,dt=0
\end{equation*}
for all $\eta \in L^{2}\left( \mathbb{T}\right) $. Conversely, one checks
that this $N$-scale condition on $\xi $ implies that $\xi $ is an
eigenvector for $R$ of eigenvalue $1$.
\end{remark}

\section{\label{Conc}Concluding remarks}

Operators of the form (\ref{Eq1.8}) or (\ref{Eq1.16}) occur in a variety of
contexts: for example, in Ruelle's work on dynamical systems \cite{Rue94}, 
\cite{BaRu96}; as operators in spaces of analytic functions \cite{CoMa95}, 
\cite{Ho96}, \cite{Lam86}, \cite{HoJa96}, \cite{LaSt91} under the names
``weighted translation operators'', ``composition operators'', or ``slash
Toeplitz operators''; and in ergodic theory \cite{Kea72}, \cite{Wal96} (in
the positive case). Our present approach is different from those mentioned
in that we ask the questions in a geometric Hilbert-space setting in $%
L^{2}\left( \mathbb{T}\right) $, and in that we make the connections between
wavelets and the theory of representations of the $C^{*}$-algebras $\mathcal{%
O}_{N}$. Our analysis of the mother functions $\psi _{1},\dots ,\psi _{N-1}$
is motivated by results in \cite{GrMa92}, \cite{GrHa94}, and \cite
{Mey93,Mey92b}, while our study of the correspondence between the $\mathcal{O%
}_{N}$-representations and $\left\{ \varphi ,\psi _{i}\right\} $, in (\ref
{Eq1.41})--(\ref{Eq1.42}), is motivated by our desire to put the results of 
\cite{CoRy95} and \cite{CoDa96} in a more general geometric and
operator-theoretic framework. Our viewpoint here, and in particular in
Section \ref{Scat}, is that of Lax and Phillips \cite{LaPh89} in their
approach to scattering on obstacles for the classical wave equation.

\begin{acknowledgements}
Most of this work was done when P.E.T.J. visited the University of Oslo with
support from the university and from NFR. We have benefitted from
discussions with Nils \O vrelid, Edwin Beggs, Ra\'{u}l Curto, Mark C. Ho,
and others. Expert typesetting by Brian Treadway is gratefully acknowledged.
\end{acknowledgements}

\bibliographystyle{bftalpha}
\bibliography{jorgen}

\end{document}